\documentclass[11pt,a4paper]{article}
\usepackage[T1]{fontenc}
\usepackage[utf8]{inputenc}
\usepackage{authblk}
\usepackage{geometry}
\usepackage{color}

\usepackage{amsmath,amsthm}
\usepackage{charter}
\usepackage[charter]{mathdesign}
\usepackage{xcolor}
\usepackage{tikz-cd}

\newcommand{\nn}{\nonumber}
\newcommand{\eps}{\varepsilon}
\newcommand{\bx}{{\bf x}}

\newtheorem{remark}{Remark}

\graphicspath{{figures/}}

\author[1]{Weizhu Bao}
\affil[1]{Department of Mathematics, National University of Singapore, Singapore 119076, Singapore}

\author[2]{Yue Feng}
\affil[2]{Laboratoire Jacques-Louis Lions, Sorbonne Universit\'e, Paris 75005, France }

\author[3]{Ying Ma}
\affil[3]{Department of Mathematics, Faculty of Science, Beijing University of Technology, Beijing 100124, China}
\affil[*]{Corresponding author: Yue Feng, yue.feng@sorbonne-universite.fr}

\date{}

\begin{document}

\title{Regularized numerical methods for the nonlinear Schr\"odinger equation with singular nonlinearity}

\maketitle
\begin{abstract}
We present different regularizations and numerical methods for the nonlinear Schr\"odinger equation with singular nonlinearity (sNLSE) including the regularized Lie-Trotter time-splitting (LTTS) methods and regularized Lawson-type exponential integrator (LTEI) methods. Due to the blowup of the singular nonlinearity, i.e., $f(\rho)=\rho^{\alpha}$ with a fixed exponent $\alpha<0$ goes to infinity when $\rho \to 0^+$ ($\rho = |\psi|^2$ represents the density with $\psi$ being the complex-valued wave function or order parameter), there are significant difficulties in designing accurate and efficient numerical schemes to solve the sNLSE. In order to suppress the round-off error and avoid blowup near $\rho = 0^+$, two types of regularizations for the sNLSE are proposed with a small regularization parameter $0 < \eps \ll 1$. One is based on the local energy regularization (LER) for the sNLSE via regularizing the energy density $F(\rho) = \frac{1}{\alpha+1}\rho^{\alpha+1}$ locally near $\rho = 0^+$ with a polynomial approximation and then obtaining a local energy regularized nonlinear Schr\"odinger equation via energy variation. The other one is the global nonlinearity regularization which directly regularizes the singular nonlinearity $f(\rho)=\rho^{\alpha}$ to avoid  blowup near $\rho = 0^+$. For the regularized models, we apply the first-order Lie-Trotter time-splitting method and Lawson-type exponential integrator method for temporal discretization and combine with the Fourier pseudospectral method in space to numerically solve them. Numerical examples are provided to show the convergence of the regularized models to the sNLSE and they suggest that the local energy regularization performs better than directly regularizing the singular nonlinearity globally.
\end{abstract}


\maketitle


\section{Introduction}
The nonlinear Schr\"odinger equation (NLSE) is a prototypical dispersive partial differential equation (PDE) playing an important role in different areas of physics, chemistry and engineering. The relevant applications vary from Bose-Einstein condensate (BEC) \cite{BC,PS}, nonlinear optics \cite{ABB,ADK} to plasma and particle physics \cite{BC,SS}.

In general, the time-dependent NLSE is in the following form \cite{CT,Sch,SS}
\begin{equation}
\left\{
\begin{aligned}
&i\partial_t \psi(\bx, t) = -\Delta \psi(\bx, t) +\lambda\, f\left(|\psi(\bx, t)|^2\right)\psi(\bx, t), \quad \bx \in \Omega, \ t>0,	\\
&\psi(\bx, 0) = \psi_0(\bx), \quad \bx \in \overline{\Omega},
\end{aligned}\right.
\label{eq:NLSE}
\end{equation}
where $i = \sqrt{-1}$ is the complex unit, $t$ is time, $\bx \in \mathbb{R}^d$ $(d = 1, 2, 3)$ is the spatial coordinate, $\psi : = \psi(\bx, t) \in \mathbb{C}$ is the dimensionless wave function or order parameter, $\psi_0 := \psi_0(\bx)$ is a given complex-valued initial data, $\lambda\ne0$ is a given real constant with $\lambda>0$ for repulsive or defocusing  interaction and $\lambda<0$ for attractive or focusing  interaction,  and $\Omega =\mathbb{R}^d$ or $\Omega \subset \mathbb{R}^d$ is a bounded domain with periodic boundary condition or homogeneous Dirichlet boundary or homogeneous Neumann boundary condition. The nonlinearity is given as \cite{PS,CT,SS}
\begin{equation}\label{nonlinear}
f(\rho) := \rho^{\alpha}, \qquad \rho\ge0,
\end{equation}
where $\rho := |\psi|^2$ is the density and the exponent $\alpha>-1$ is a real constant, which is different in diverse applications. Specifically, when $\alpha = 1$, i.e., $f(\rho) = \rho$, it is the most popular NLSE with cubic nonlinearity and  also called Gross--Pitaevskii equation (GPE), especially in BEC \cite{PS,Sch,SS}; and when $\alpha = 2$, it is related to  the quintic Schr\"odinger equation, which is regarded as the mean field limit of a Boson gas with three-body interactions and also widely used in the study of optical lattices \cite{CP,RA}. When $0<\alpha<1$ or $1<\alpha<2$, it is usually stated that the NLSE with semi-smooth (or fractional) nonlinearity, which has been adapted in different applications \cite{Lee,Bao04,Cai,Ig}. For the NLSE with smooth or semi-smooth nonlinearity, i.e., $\alpha>0$,  the existence and uniqueness of the Cauchy problem as well as the finite time blow-up have been widely studied \cite{CT,SS}. Recently, interests have been surged for the study of the NLSE \eqref{eq:NLSE} with singular nonlinearity \eqref{nonlinear}, i.e., $\alpha \in (-1,0)$. In this case, the NLSE \eqref{eq:NLSE} can be formally obtained as the nonrelativistic limit of the nonlinear Dirac equation with singular (or fractional) nonlinearity \cite{M1,Ma2}, which was proposed as a model of strong interaction of particles and it recovered the MIT bag model \cite{Cho,Jo}. When $\alpha<0$ in \eqref{nonlinear}, the nonlinearity $f(\rho)$ has a singularity at the origin and it is also called sublinear Schr\"odinger equation for the case $\alpha \in (-1/2, 0)$ in the mathematical literature \cite{BOB,KRI}. The study of the NLSE with singular nonlinearity is much more complicated in both analytical and numerical aspects. In recent years, dispersive PDEs with singular nonlinearity have attracted much attention, e.g., the existence of standing waves for nonlinear Dirac fields has been proven and the solution is of class $C^1$ when $-1/3 < \alpha < 0$, while $|\nabla \psi|$ is infinite on some sphere $\{|x| = R\}$ for $-1 < \alpha < -1/3$ \cite{BCV}. Since the nonlinear Schr\"odinger equation is the nonrelativistic limit of the nonlinear Dirac equation, it is also an interesting and challenging problem to study the nonlinear Schr\"odinger equation with such a singular nonlinearity. The existence and multiplicity of solutions for the nonlinear Schr\"odinger equation with singular nonlinearity have been studied in the literature and references therein \cite{BOB,KRI,ZW}.

In this paper, for the sake of simplicity, we focus on the following singular nonlinear Schr\"odinger equation (sNLSE) with $\alpha \in (-1/3, 0)$
\begin{equation}
\left\{
\begin{aligned}
&i\partial_t \psi(\bx, t) = -\Delta \psi(\bx, t) +\lambda |\psi(\bx, t)|^{2\alpha}\psi(\bx, t), \quad \bx \in \Omega, \ t>0,	\\
&\psi(\bx, 0) = \psi_0(\bx), \quad \bx \in \overline{\Omega}.
\end{aligned}\right.
\label{eq:NLSE1}
\end{equation}
Similar to the cubic nonlinear Schr\"odinger equation, the sNLSE \eqref{eq:NLSE1} conserves the {\em mass}, {\em momentum} and {\em energy} as \cite{BC,CT}:
\begin{align}
M(t) & := \left\|\psi(\cdot, t)\right\|^2 = \int_{\Omega}|\psi(\bx, t)|^2 d\bx \equiv  	\int_{\Omega}|\psi_0(\bx)|^2 d\bx := M(0), \quad t \geq 0,\\
P(t) & := {\rm Im}\int_{\Omega}\overline{\psi}(\bx, t)\nabla \psi(\bx, t) d\bx \equiv  	{\rm Im}\int_{\Omega}\overline{\psi_0}(\bx)\nabla\psi_0(\bx) d\bx := P(0), \quad t \geq 0,\\
E(t) & := \int_{\Omega}\Bigg[|\nabla\psi(\bx, t)|^2 + \lambda\, F(|\psi(\bx, t)|^2)\Bigg] d\bx \nn\\
& \equiv \int_{\Omega}\Bigg[|\nabla\psi_0(\bx)|^2 + \lambda\, F(|\psi_0(\bx)|^2)\Bigg] d\bx := E(0), \quad t \geq 0,
\end{align}
where ${\rm Im}f$ and $\overline{f}$ denote the imaginary part and complex conjugate of $f$, respectively, and
\begin{equation}
F(\rho) = \int^{\rho}_0 f(s) ds = \frac{1}{\alpha+1} \rho^{\alpha+1}, \quad \rho \geq 0.
\end{equation}

As a consequence of their importance in various applications,  a large number of numerical schemes have been proposed and analyzed for the NLSE with smooth nonlinearity in the past decades \cite{ABB,BC,BBD,BJM,GL}, including the finite difference methods \cite{DFP,SV}, finite element methods \cite{AD,KAD}, exponential integrator methods \cite{CCO,HO,OS}, time-splitting methods \cite{BBD,TM,WH}. However, these methods for the smooth nonlinearity cannot be directly applied to the singular nonlinearity due to the blowup of the nonlinear term, i.e., $f(\rho) =\rho^{\alpha} \to \infty$ when $\rho \to 0^+$. As far as we know, there is few work on the numerical simulation and analysis on the NLSE with singular nonlinearity $f(\rho) = \rho^{\alpha}$ for $\alpha < 0$. The singularity of the nonlinear term makes it much more challenging to design accurate and efficient numerical schemes and establish their error estimates. Very recently, to deal with the logarithmic nonlinearity, different regularized methods were proposed and analyzed to study the logarithmic Schr\"odinger equation (LogSE) \cite{BCST0,BCST1,BCST2}. One is to introduce a small parameter $0 < \eps \ll 1$ to regularize the logarithmic nonlinearity in the LogSE globally and the convergence was established between the solutions of the LogSE and the regularized equation \cite{BCST0,BCST1}. The other strategy is to regularize the energy density locally in the region $\{0 \le \rho < \eps^2\}$ by a sequence of polynomials and remain it unchanged in the region $\{\rho \ge \eps^2\}$ and the convergence was also established \cite{BCST2}. In addition, different numerical schemes are applied to solve the regularized models including the finite difference method and time-splitting method. The main aim of this paper is to introduce two different types of regularizations to regularize the sNLSE locally or globally and solve the regularized models by the first-order Lie-Trotter time-splitting Fourier pseudospectral method and Lawson-type exponential integrator Fourier pseudospectral method as well as compare the performance of different regularizations and numerical schemes.

The rest of this paper is organized as follows. In section 2, we introduce two types of regularized models to regularize the singular nonlinear  Schr\"odinger equation (sNLSE) and compare their approximations of the sNLSE. In section 3, we apply the first-order Lie-Trotter time-slitting method and Lawson-type exponential integrator to discretize the regularized nonlinear Schr\"odinger equations in time and combine with the Fourier pseudospectral method for spatial discretization. Numerical results are shown in section 4 to validate our regularized models and compare the performance of different numerical schemes. Finally, some conclusions are drawn in section 5.

\section{Regularized nonlinear Schr\"odinger equation}
In this section, we introduce two types of regularized nonlinear Schr\"odinger equations (rNLSEs) to regularize the nonlinear Schr\"odinger equation with singular nonlinearity (sNLSE) and analyze their properties. One is based on the local energy regularization (LER) for the energy density function $F(\rho) = \frac{1}{\alpha + 1} \rho^{\alpha+1}$ in the region $\{0\le \rho < \eps^2\}$ and remains it unchanged in the region  $\{\rho \ge \eps^2\}$.  The other type is directly regularizing the singular nonlinearity $f(\rho) = \rho^{\alpha}$ globally via merging the small regularization parameter $0 < \eps \ll 1$ by different ways.

\subsection{The rNLSEs with local energy regularization}
We consider a local regularization by approximating  the energy density $F(\rho) = \frac{1}{\alpha+1}\rho^{\alpha+1}$ only in the region $\{0\le \rho <\eps^2\}$ and keeping it unchanged in the region $\{\rho \ge \eps^2\}$, which is a local energy regularization (LER) \cite{BCST2}. In the region $\{0\le \rho < \eps^2\}$, we use a sequence of polynomials to approximate the energy density $F(\rho)$ as
\begin{equation}
F^{\eps}_n(\rho) = F(\rho) \chi_{\{\rho \geq \eps^2\}} + P^{\eps}_{n+1}(\rho)\chi_{\{\rho < \eps^2\}}, \qquad \rho\ge0,
\label{eq:Feps}
\end{equation}
where $0 < \eps \ll 1$ is a small regularization parameter, $n \geq 1$ is an arbitrary integer, $\chi_A$ is the characteristic function of the set $A$, and $P^{\eps}_{n+1} $ is a polynomial of degree $n+1$. We require the piecewise smooth function $F^{\eps}_n(\rho) \in C^n([0, +\infty])$ and $F^{\eps}_n(0) = F(0) = 0$. Then, we are going to derive the explicit expression of $P^{\eps}_{n+1}(\rho)$. Noticing $P^{\eps}_{n+1}(0) = 0$, we can write
\begin{equation}
P^{\eps}_{n+1}(\rho) = \rho 	Q^{\eps}_{n}(\rho),\qquad \rho\ge0,
\end{equation}
where $Q^{\eps}_{n}(\rho)$ is a polynomial of degree of $n$. Accordingly, denote
 \begin{equation}
 F(\rho) = \rho Q(\rho), \quad Q(\rho) = \frac{1}{\alpha+1}\rho^{\alpha},\qquad \rho\ge0.	
 \end{equation}
The continuity conditions at $\rho = \eps^2$ are
\begin{equation*}
P^{\eps}_{n+1}(\eps) = F(\eps^2), \quad  (P^{\eps}_{n+1})'(\eps) = F'(\eps^2),\quad \ldots, \quad (P^{\eps}_{n+1})^{(n)}(\eps) = F^{(n)}(\eps^2),
\end{equation*}
which imply
\begin{equation*}
Q^{\eps}_{n+1}(\eps) = Q(\eps^2), \quad  (Q^{\eps}_{n+1})'(\eps) = Q'(\eps^2),\quad \ldots, \quad (Q^{\eps}_{n+1})^{(n)}(\eps) = Q^{(n)}(\eps^2).
\end{equation*}
As a consequence, $Q^{\eps}_{n}(\rho)$ is the $n$-th degree Taylor polynomial of $Q$ at $\rho = \eps^2$, i.e.,
\begin{align}
Q^{\eps}_{n}(\rho) & = Q(\eps^2) + \sum_{k=1}^n \frac{Q^{(k)}(\eps^2)}{k!}\left(\rho-\eps^2\right)^k \nn \\
& = 	\frac{\eps^{2\alpha}}{\alpha+1}\left[1+\sum_{k=1}^n \prod_{j=1}^k \frac{-\alpha+j-1}{j}\left(1-\frac{\rho}{\eps^2}\right)^k\right],\quad \rho\ge0.
\label{eq:Qeps}
\end{align}
Differentiating \eqref{eq:Feps} with respect to $\rho$ and combining \eqref{eq:Qeps}, we get
\begin{equation}
f^{\eps}_n(\rho) = (F^{\eps}_n)'(\rho) = \rho^{\alpha}\chi_{\{\rho\geq\eps^2\}} + q^{\eps}_n(\rho)  \chi_{\{\rho<\eps^2\}}, \quad \rho \geq 0,
\label{eq:feps}
\end{equation}
where
\begin{align}
q^{\eps}_n(\rho) & =  (P^{\eps}_{n+1})'(\rho)	= Q^{\eps}_n(\rho) + \rho (Q^{\eps}_n)'(\rho)\nn\\
& =	\frac{\eps^{2\alpha}}{\alpha+1}\left[1+\sum_{k=1}^n \prod_{j=1}^k \frac{-\alpha+j-1}{j}\left(1-\frac{\rho}{\eps^2}\right)^{k-1}
\left(1-\frac{\rho(k+1)}{\eps^2}\right)\right],\quad \rho\ge0.
\end{align}
Therefore, we obtain the following energy regularized nonlinear Schr\"odinger equation (erNLSE) with a small regularized parameter $0 < \eps \ll 1$ as
\begin{equation}
\left\{
\begin{aligned}
&i\partial_t \psi^{\eps}(\bx, t) = -\Delta \psi^{\eps}(\bx, t) +\lambda\,f^{\eps}_n(|\psi(\bx, t)|^2)\psi^{\eps}(\bx, t), \quad \bx \in \Omega, \ t>0,	\\
&\psi^{\eps}(\bx, 0) = \psi_0(\bx), \quad \bx \in \overline{\Omega},
\end{aligned}\right.
\label{eq:ERNLSE}
\end{equation}
with $f^{\eps}_n$ defined in \eqref{eq:feps}. In addition, it is easy to check the erNLSE \eqref{eq:ERNLSE} conserves the {\em mass}, {\em momentum} and {\em energy} as:
\begin{align}\label{massl}
M(t) & := \left\|\psi^\eps(\cdot, t)\right\|^2 = \int_{\Omega}|\psi^\eps(\bx, t)|^2 d\bx \equiv  	\int_{\Omega}|\psi_0(\bx)|^2 d\bx := M(0), \quad t \geq 0,\\ \label{momentuml}
P(t) & := {\rm Im}\int_{\Omega}\overline{\psi}^\eps(\bx, t)\nabla \psi^\eps(\bx, t) d\bx \equiv  	{\rm Im}\int_{\Omega}\overline{\psi_0}(\bx)\nabla\psi_0(\bx) d\bx := P(0), \quad t \geq 0,\\
E^{\eps}_n(t) & := \int_{\Omega}\Bigg[|\nabla\psi^{\eps}(\bx, t)|^2 + \lambda\,F^{\eps}_n(|\psi^{\eps}(\bx, t)|^2)\Bigg] d\bx \nn\\
& \equiv \int_{\Omega}\Bigg[|\nabla\psi_0(\bx)|^2 +\lambda\,F^{\eps}_n(|\psi_0(\bx)|^2)\Bigg] d\bx = E^{\eps}_n(0), \quad t \geq 0,
\end{align}
where $F^{\eps}_n$ is defined in \eqref{eq:Feps}.

\subsection{The rNLSEs with global regularization}
In contrast to the previous local energy regularization, we consider the global regularization by directly regularizing the singular nonlinearity $f(\rho) = \rho^{\alpha}$ via the following regularization:
\begin{equation}\label{globalr1}
f^{\eps}(\rho)= (\rho+\eps^2)^{\alpha}, \qquad \rho\ge0,
\end{equation}
where $0 < \eps \ll 1$ is  a small regularization parameter. Then, we obtain the following rNLSE
\begin{equation}
\left\{
\begin{aligned}
&i\partial_t \psi^{\eps}(\bx, t) = -\Delta \psi^{\eps}(\bx, t) +\lambda\left(|\psi^{\eps}(\bx, t)|^2 +\eps^2\right)^{\alpha}\psi^{\eps}(\bx, t), \quad \bx \in \Omega, \ t>0,	\\
&\psi^{\eps}(\bx, 0) = \psi_0(\bx), \quad \bx \in \overline{\Omega}.
\end{aligned}\right.
\label{eq:RNLSE1}
\end{equation}
Similarly, the rNLSE \eqref{eq:RNLSE1} also conserves the {\em mass} in \eqref{massl}, {\em momentum} in \eqref{momentuml} and the following {\em energy} as:
\begin{align}
E^{\eps}(t) & := \int_{\Omega}\Bigg[|\nabla\psi^{\eps}(\bx, t)|^2 + \lambda\, F^{\eps}(|\psi^{\eps}(\bx, t)|^2)\Bigg] d\bx \nn\\
& \equiv \int_{\Omega}\Bigg[|\nabla\psi_0(\bx)|^2 + \lambda\, F^{\eps}(|\psi_0(\bx)|^2)\Bigg] d\bx =  E^{\eps}(0), \quad t \geq 0,
\end{align}
where
\begin{equation}
F^{\eps}(\rho) = \int^{\rho}_0 (s+\eps^2)^{\alpha} ds = \frac{1}{\alpha+1}\left[(\rho+\eps^2)^{\alpha+1}-\eps^{2(\alpha+1)}\right], \quad \rho \geq 0.
\end{equation}

For comparison, we can also use  another function
\begin{equation}\label{globalr2}
\widetilde{f}^{\eps}(\rho)= 1/(\rho^{-\alpha}+\eps),\qquad \rho\ge0,
\end{equation}
to regularize the singular nonlinearity $f(\rho)$ in \eqref{nonlinear} and obtain the following rNLSE
\begin{equation}
\left\{
\begin{aligned}
&i\partial_t \psi^{\eps}(\bx, t) = -\Delta \psi^{\eps}(\bx, t) + \frac{\lambda}{|\psi^{\eps}(\bx, t)|^{-2\alpha}+\eps}\psi^{\eps}(\bx, t), \quad \bx \in \Omega, \ t>0,	\\
&\psi^{\eps}(\bx, 0) = \psi_0(\bx), \quad \bx \in \overline{\Omega}.
\end{aligned}\right.
\label{eq:RNLSE2}
\end{equation}
Again, the rNLSE \eqref{eq:RNLSE2} conserves the {\em mass} in \eqref{massl}, {\em momentum} in \eqref{momentuml} and the following {\em energy} as:
\begin{align}
\widetilde E^{\eps}(t) & := \int_{\Omega}\Bigg[|\nabla\psi^{\eps}(\bx, t)|^2 + \lambda\,\widetilde F_{\eps}(|\psi^{\eps}(\bx, t)|^2)\Bigg] d\bx \nn\\
& \equiv \int_{\Omega}\Bigg[|\nabla\psi_0(\bx)|^2 + \lambda\,\widetilde F_{\eps}(|\psi_0(\bx)|^2)\Bigg] d\bx = \widetilde E^{\eps}(0), \quad t \geq 0,
\end{align}
where
\begin{equation}
\widetilde F^{\eps}(\rho) = \int^{\rho}_0 \frac{1}{s^{-\alpha}+\eps} ds, \quad \rho \geq 0.
\end{equation}

The nonlinear term $f(\rho) = \rho^{\alpha}$ with $\alpha < 0$ is singular only at the origin, but the above two regularized models \eqref{eq:RNLSE1} and \eqref{eq:RNLSE2} regularize the nonlinear term globally, i.e., it not only changes the nonlinear term on the regime $\{0\le \rho < \eps^2\}$ (where $\rho = |\psi^{\eps}|^2$), but also has small effect on the regime $\{\rho \geq \eps^2\}$. However, the convergence between the solutions  of the regularized model \eqref{eq:RNLSE1} and the sNLSE \eqref{eq:NLSE1} is sublinear and dependent on $\alpha$ in term of $\eps$, which is different from that for the regularized model \eqref{eq:RNLSE2}. Based on our extensive numerical results reported later, we demonstrate linear convergence between the solutions of the sNLSE \eqref{eq:NLSE1}  and the regularized model \eqref{eq:RNLSE2}, i.e.,
\begin{equation}
\sup_{t \in [0, T]} \left\|\psi^{\eps} - \psi(t)\right\|_{L^2(\Omega)} = O(\eps), \quad T > 0.
\end{equation}
In fact, compared with the rNLSE \eqref{eq:RNLSE1}, the rNLSE \eqref{eq:RNLSE2} converges much faster to the sNLSE \eqref{eq:NLSE1} and
the convergent rate is independent of $\alpha\in(-1/3,0)$.

\begin{figure}[h!]
\begin{minipage}{0.5\textwidth}
\centerline{\includegraphics[width=7.5cm,height=5cm]{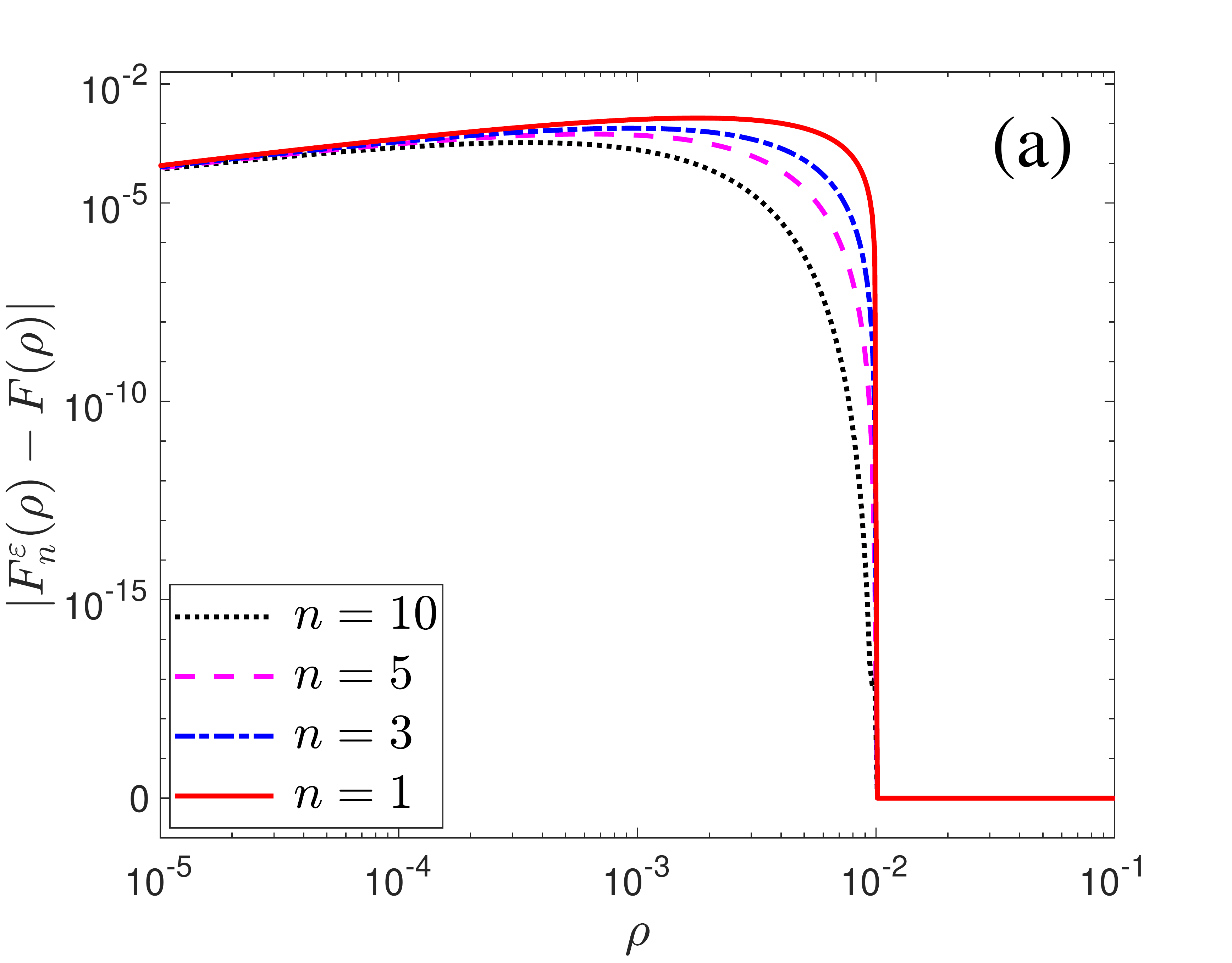}}
\end{minipage}
\begin{minipage}{0.5\textwidth}
\centerline{\includegraphics[width=7.5cm,height=5cm]{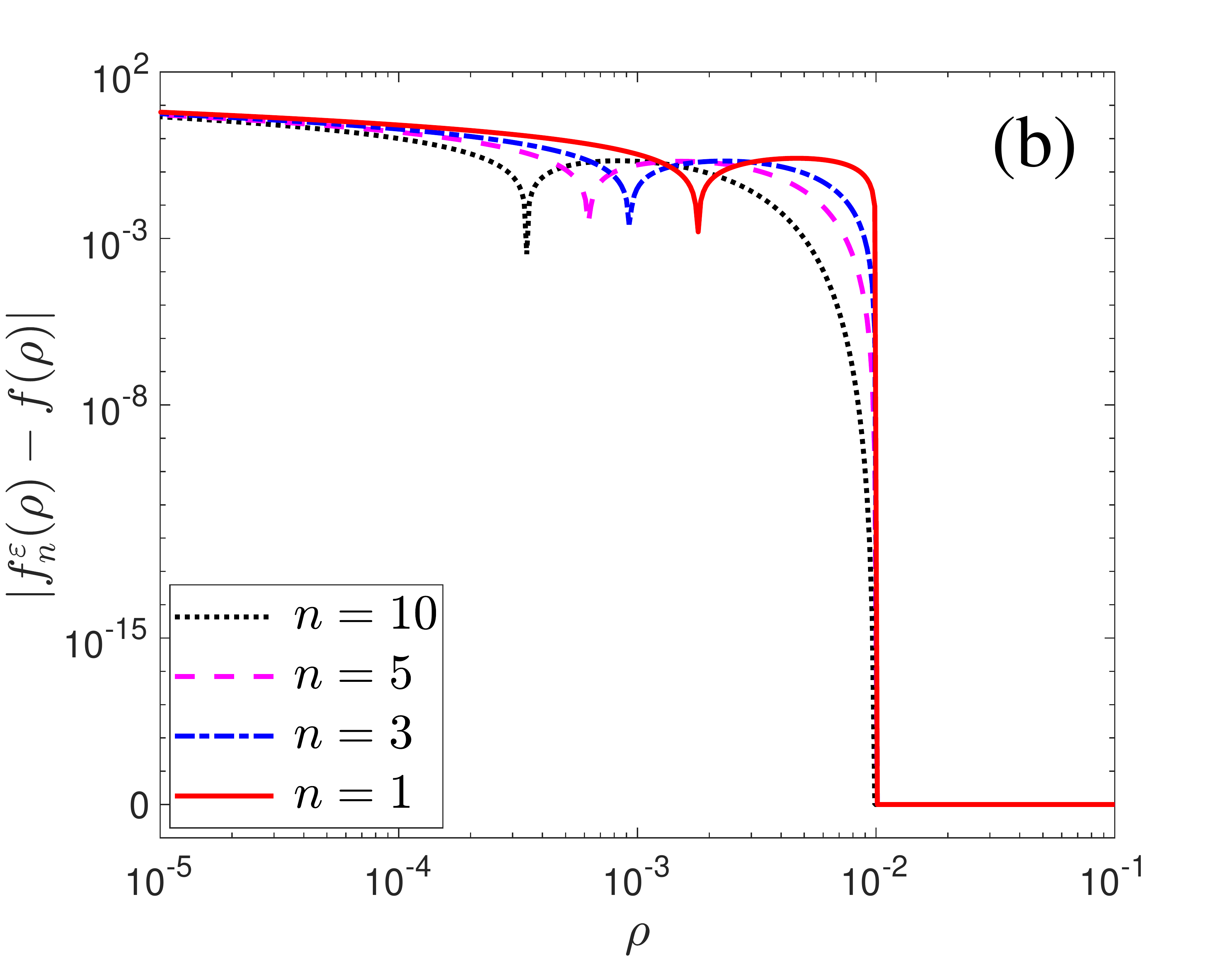}}
\end{minipage}
\caption{Errors of local energy regularization \eqref{eq:Feps} with different degrees $n$ for  (a) density $F(\rho) = \frac{1}{\alpha+1}\rho^{\alpha+1}$; and (b) singular nonlinearity $f(\rho) = \rho^{\alpha}$ with $\alpha=-0.2$.}
\label{fig:dif_n}
\end{figure}
\subsection{Comparisons of different regularizations}
The regularized function $q^{\eps}_n$ is decreasing in $[0, \eps^2]$, $f^{\eps}$ and $\widetilde{f}^{\eps}$ are decreasing on $[0, \infty)$, thus all these regularizations preserve the convexity of $F$. In addition, for the approximation of the energy density $F(\rho) = \frac{1}{\alpha+1}\rho^{\alpha+1} \in C^0([0, \infty)) \cap C^{\infty}((0, \infty))$, we have $F^{\eps}_n \in C^n([0, \infty))$ for any $n \geq 1$, while $F^{\eps} \in C^{\infty}([0, \infty))$ and $\widetilde{F}^{\eps}  \in C^1([0, \infty)) \cap C^{\infty}((0, \infty))$. Similarly, for the approximation of the singular nonlinearity $f(\rho) = \rho^{\alpha} \in C^{\infty}((0, \infty))$, we have $f^{\eps}_n \in C^{n-1}([0, \infty))$ for any $n \geq 1$, while $f^{\eps} \in C^{\infty}([0, \infty))$ and $\widetilde{f}^{\eps}  \in C^0([0, \infty)) \cap C^{\infty}((0, \infty))$.

\begin{figure}[t!]
\begin{minipage}{0.5\textwidth}
\centerline{\includegraphics[width=7.5cm,height=5cm]{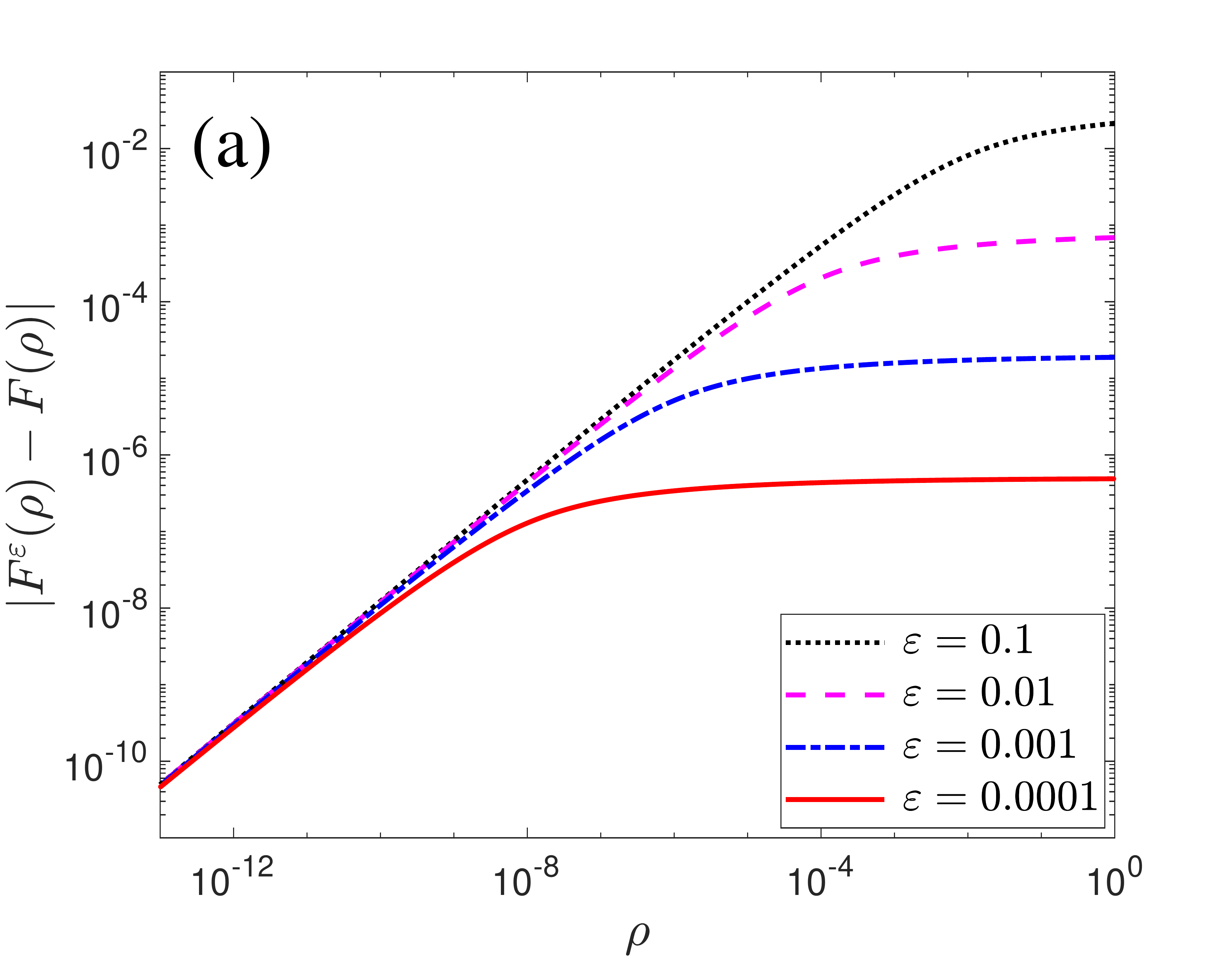}}
\end{minipage}
\begin{minipage}{0.5\textwidth}
\centerline{\includegraphics[width=7.5cm,height=5cm]{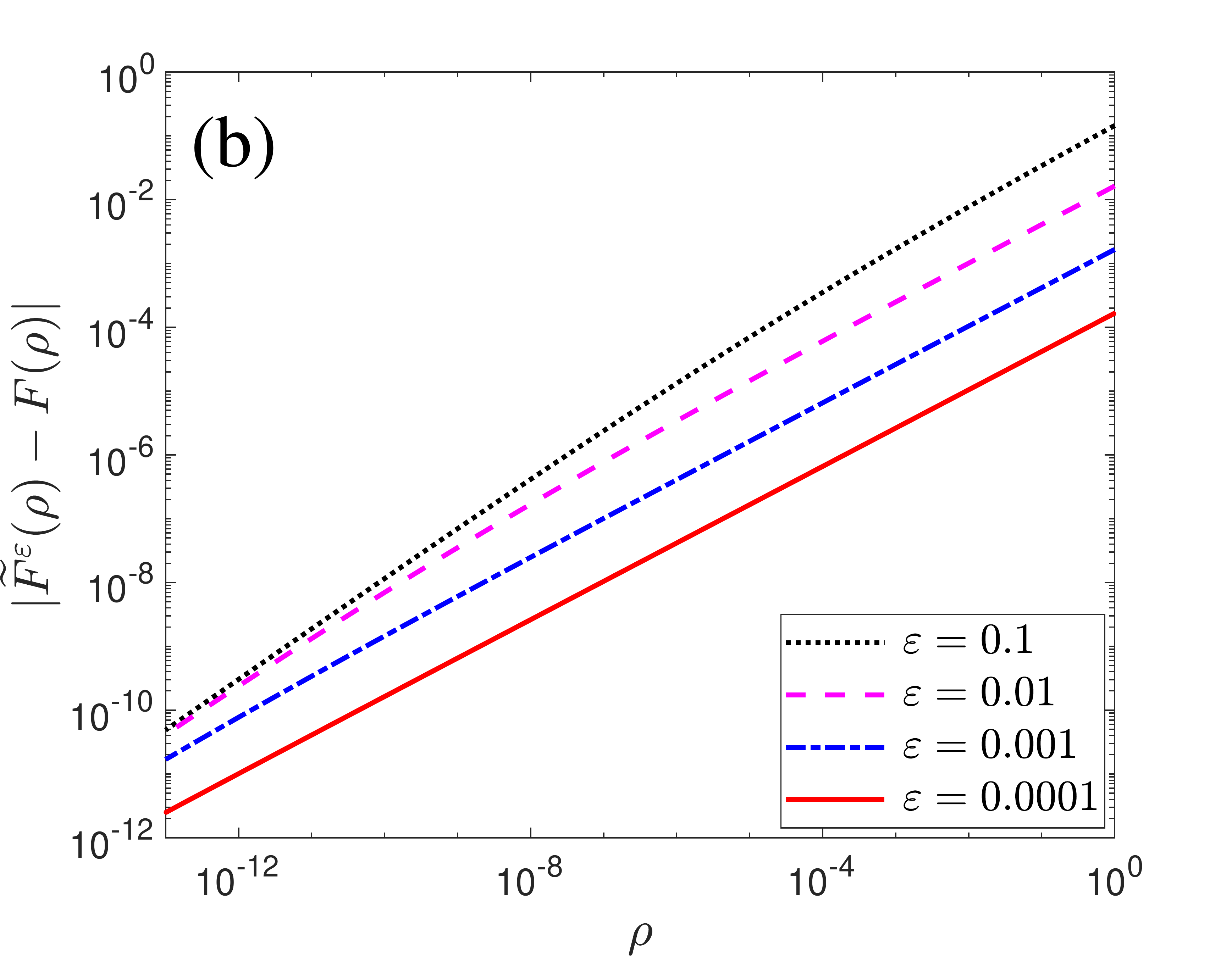}}
\end{minipage}
\begin{minipage}{0.5\textwidth}
\centerline{\includegraphics[width=7.5cm,height=5cm]{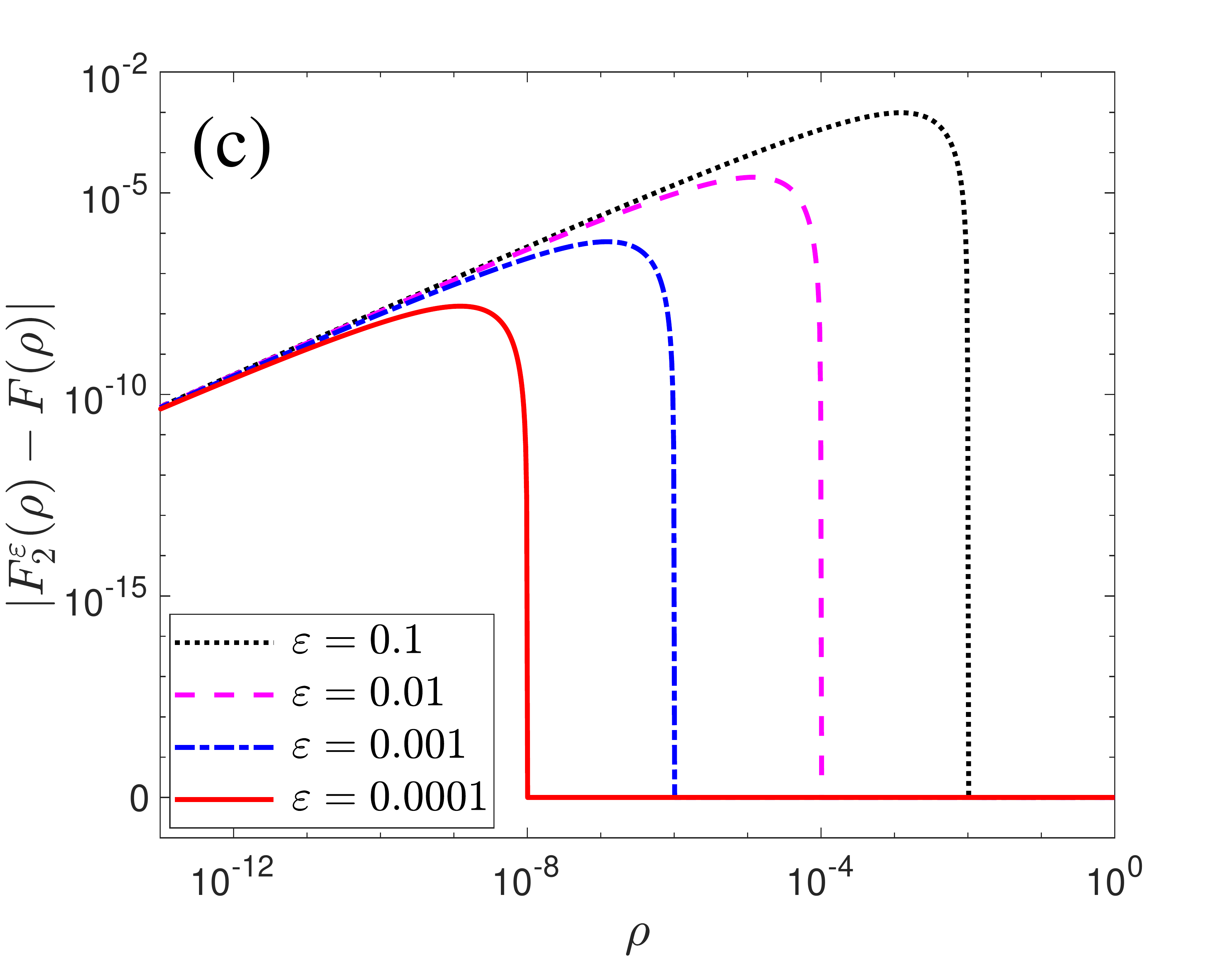}}
\end{minipage}
\begin{minipage}{0.5\textwidth}
\centerline{\includegraphics[width=7.5cm,height=5cm]{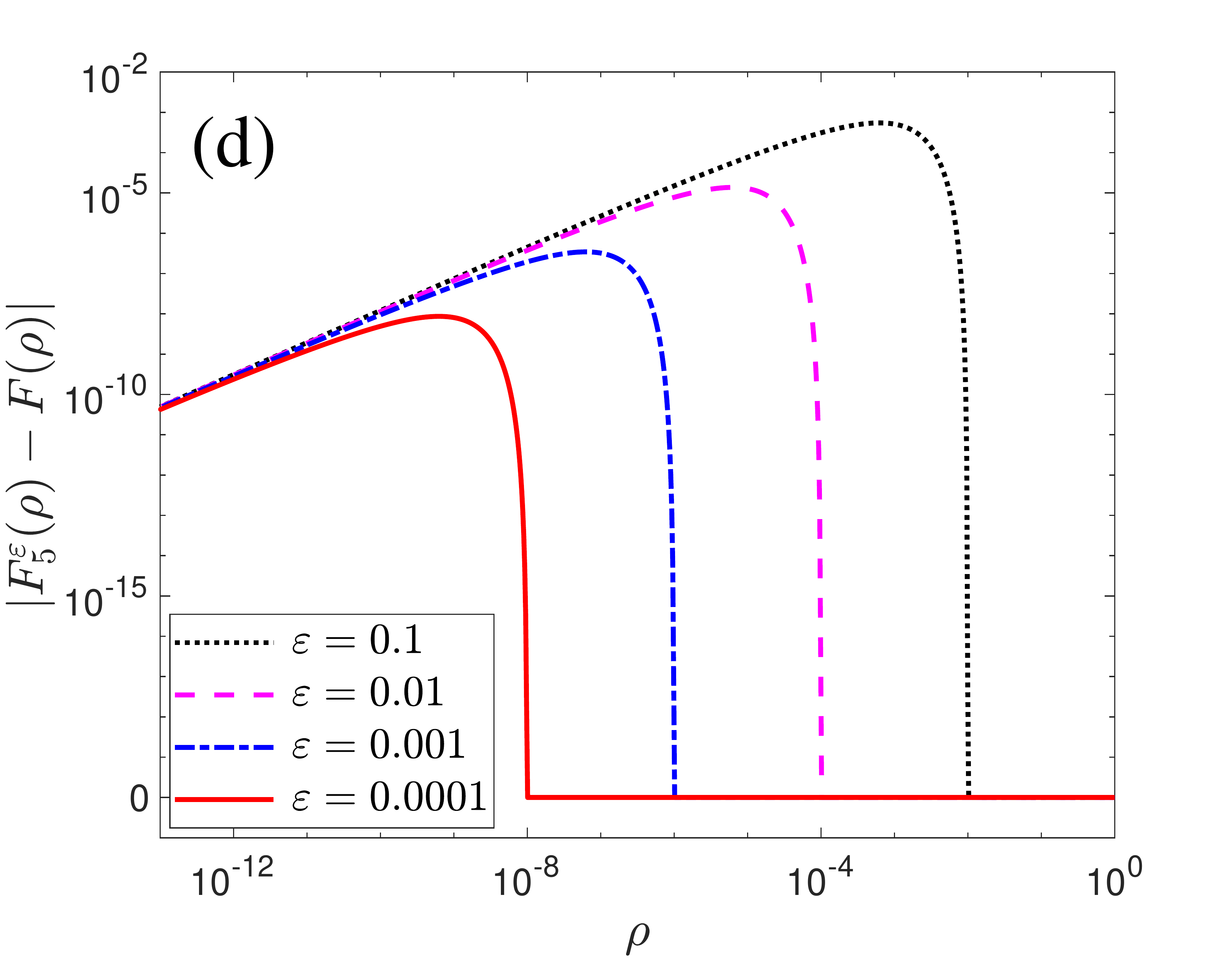}}
\end{minipage}
\caption{Comparison of different regularizations for $F(\rho) = \frac{1}{\alpha+1}\rho^{\alpha+1}$ with $\alpha = -0.2$ for different $\eps$: (a) $F^{\eps}(\rho)$; (b) $\widetilde{F}^{\eps}(\rho)$; (c) $F^{\eps}_2(\rho)$; and (d) $F^{\eps}_5(\rho)$.}
\label{fig:F_al02}
\end{figure}

\begin{figure}[ht!]
\begin{minipage}{0.5\textwidth}
\centerline{\includegraphics[width=7.5cm,height=5cm]{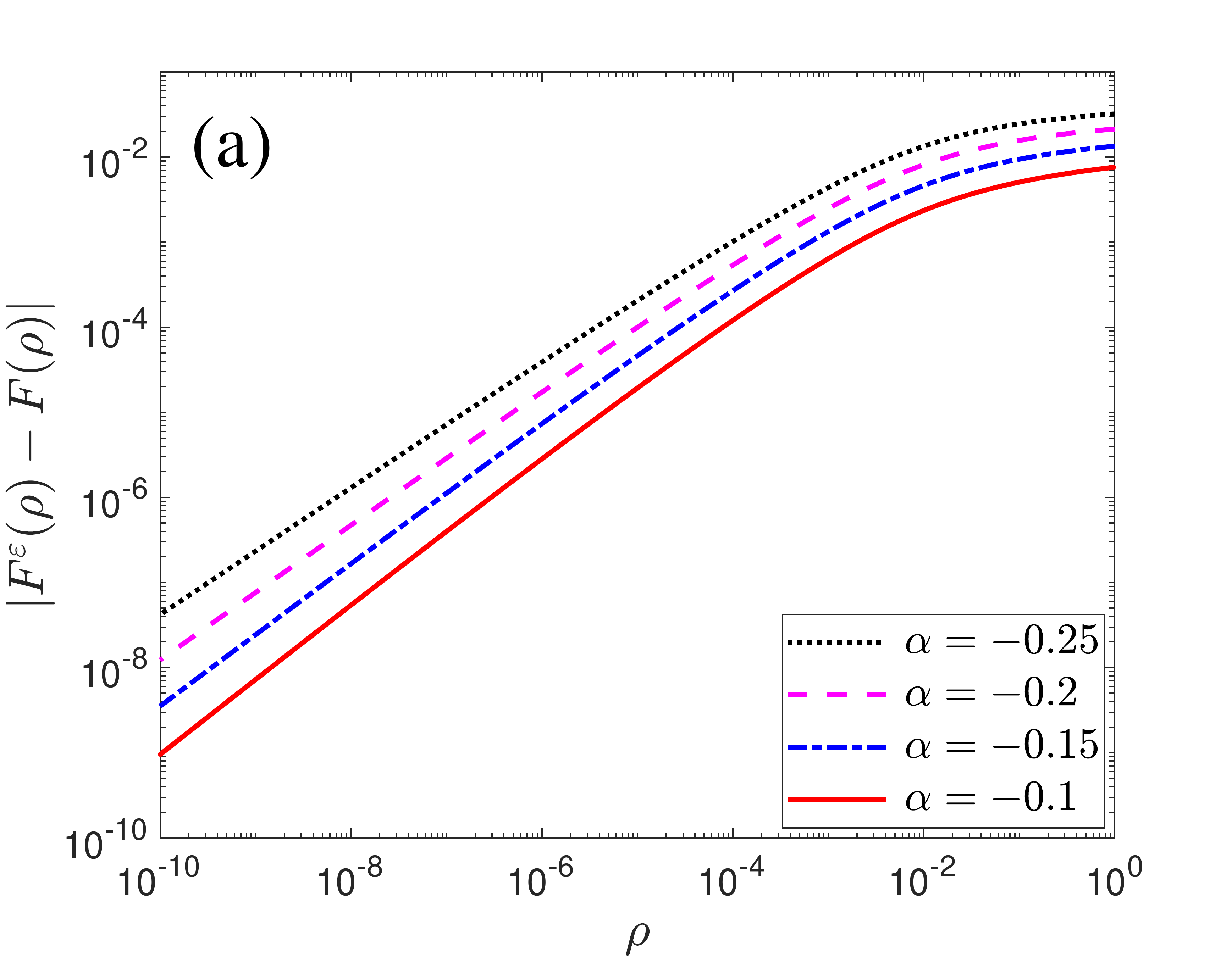}}
\end{minipage}
\begin{minipage}{0.5\textwidth}
\centerline{\includegraphics[width=7.5cm,height=5cm]{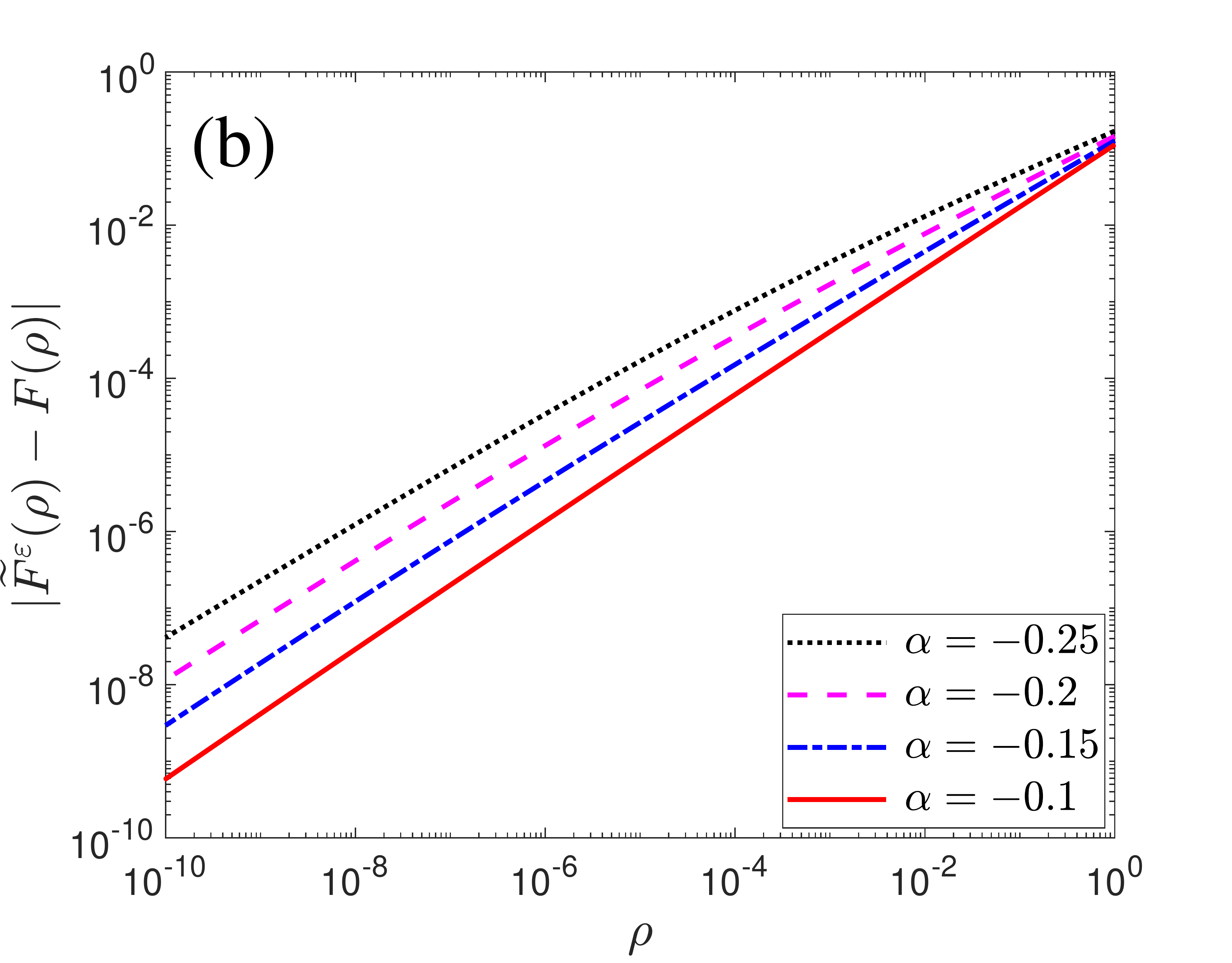}}
\end{minipage}
\begin{minipage}{0.5\textwidth}
\centerline{\includegraphics[width=7.5cm,height=5cm]{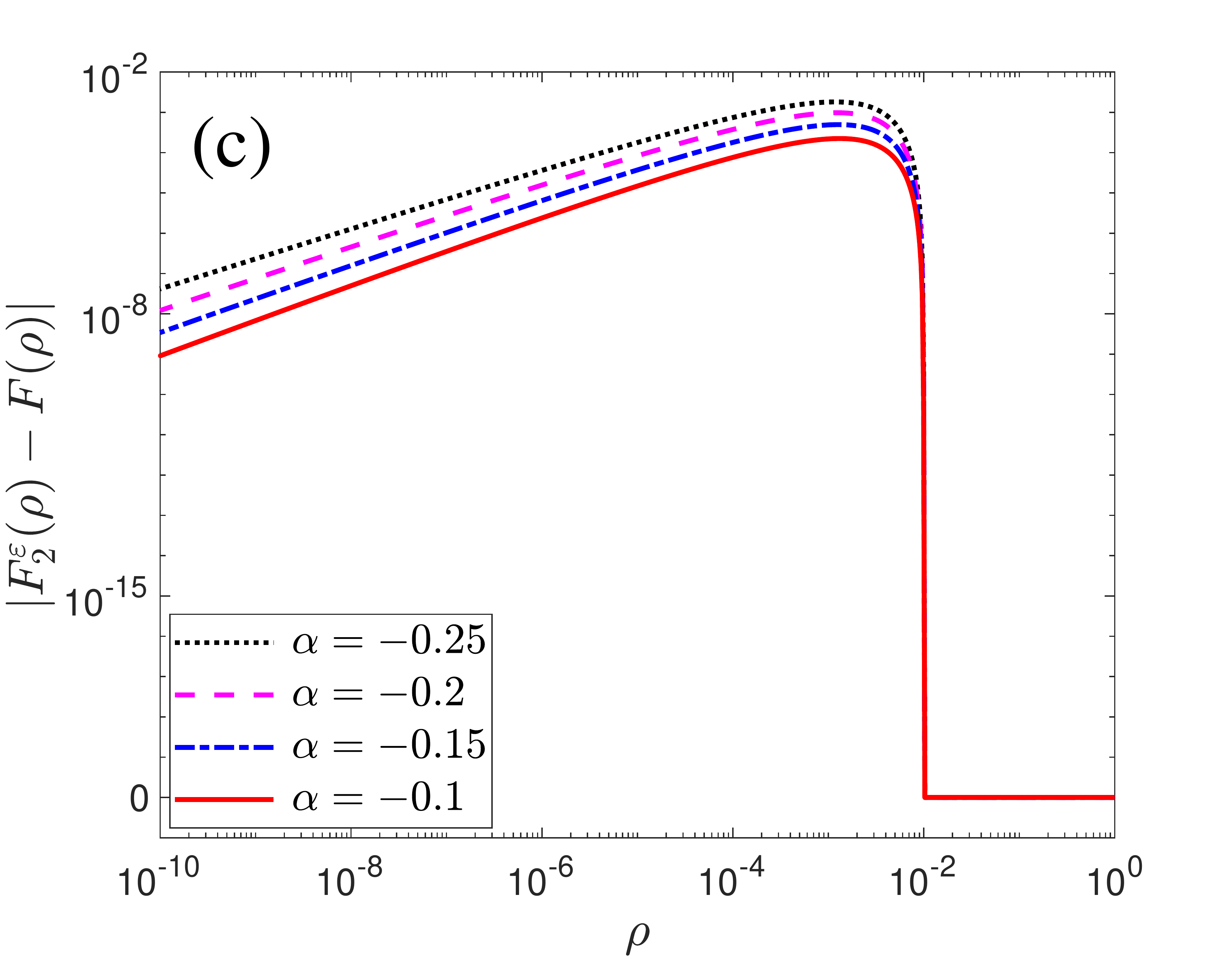}}
\end{minipage}
\begin{minipage}{0.5\textwidth}
\centerline{\includegraphics[width=7.5cm,height=5cm]{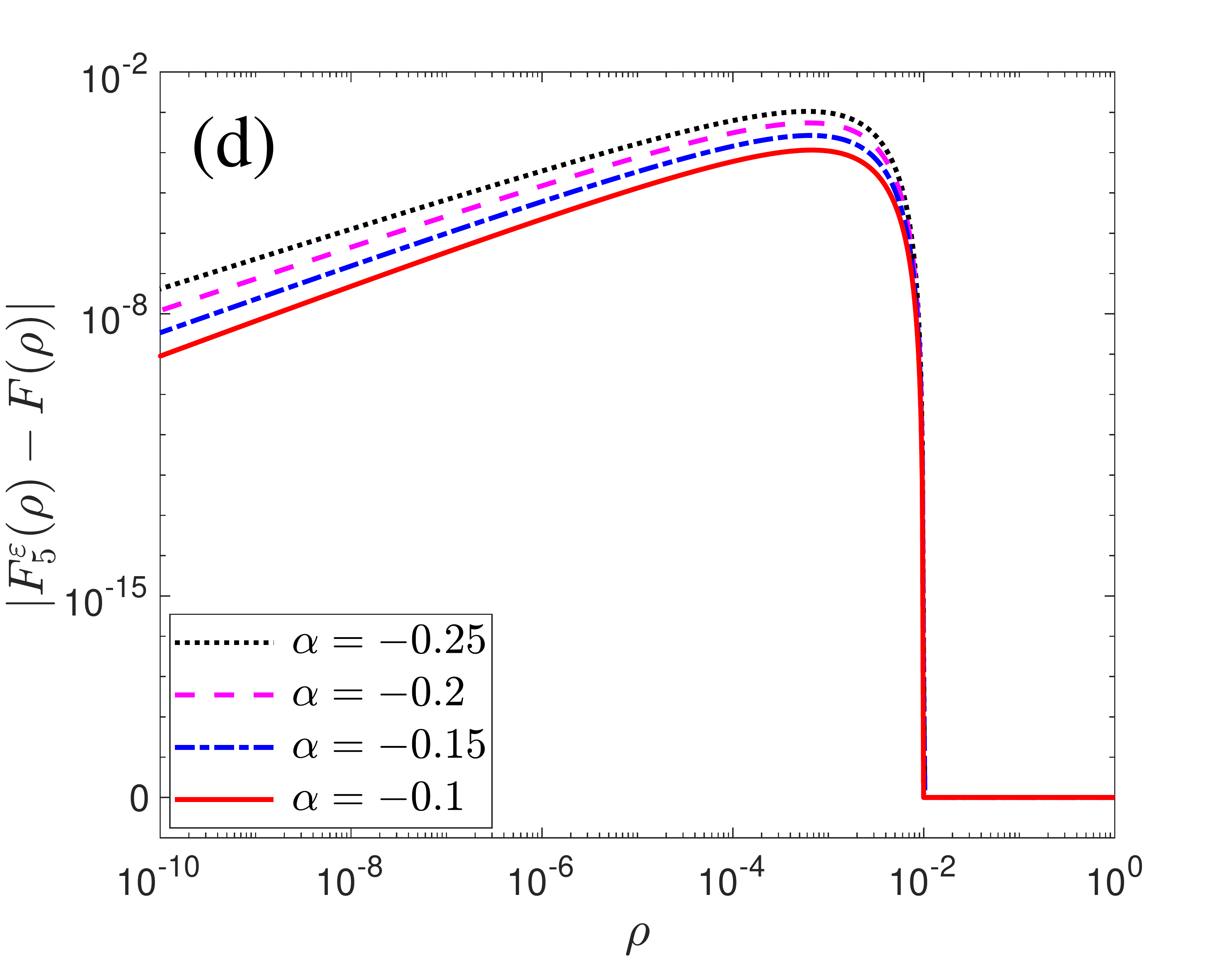}}
\end{minipage}
\caption{Comparison of different regularizations for $F(\rho) = \frac{1}{\alpha+1}\rho^{\alpha+1}$ with $\eps = 0.1$ for different $\alpha$: (a) $F^{\eps}(\rho)$; (b) $\widetilde{F}^{\eps}(\rho)$; (c) $F^{\eps}_2(\rho)$; and (d) $F^{\eps}_5(\rho)$.}
\label{fig:F_difalpha}
\end{figure}

For the local energy regularization, we use a polynomial of degree $n+1$ to approximate the energy density function $F(\rho) = \frac{1}{\alpha+1}\rho^{\alpha+1}$ in the region $\{0\le \rho < \eps^2\}$. Figure \ref{fig:dif_n} displays the approximations for $F(\rho)$ and $f(\rho)$ with $\alpha = -0.2$ for different degrees, which indicates that higher degree approximates $F(\rho)$ (and thus for $f(\rho)$) better in the regime close to the origin.  Figures \ref{fig:F_al02} and \ref{fig:F_difalpha} show the regularizations $F^{\eps}$, $\widetilde{F}^{\eps}$ and $F^{\eps}_n$ $(n = 2, 5)$ for different $\eps$ and $\alpha$, respectively. From these two figures, we find that the global regularizations approximate the function $F(\rho)$ well around the origin and introduce some additional errors when $\rho$ is far way from the origin, while the local energy regularization just change the values in the region $\{0\le \rho < \eps^2\}$, which indicate that the local regularization $F^{\eps}_n(\rho)$ approximates the energy density function $F(\rho)$ more accurately than the global regularization.  From Figure \ref{fig:F_al02}, we observe that for the fixed density $\rho$, the errors of the approximations become smaller when $\eps$ is smaller. Figure \ref{fig:F_difalpha} indicates that the errors of these three approximations become larger when $\alpha$ is smaller, and local energy regularization is the best choice in these three regularizations.

\begin{figure}[ht!]
\begin{minipage}{0.5\textwidth}
\centerline{\includegraphics[width=7.5cm,height=5cm]{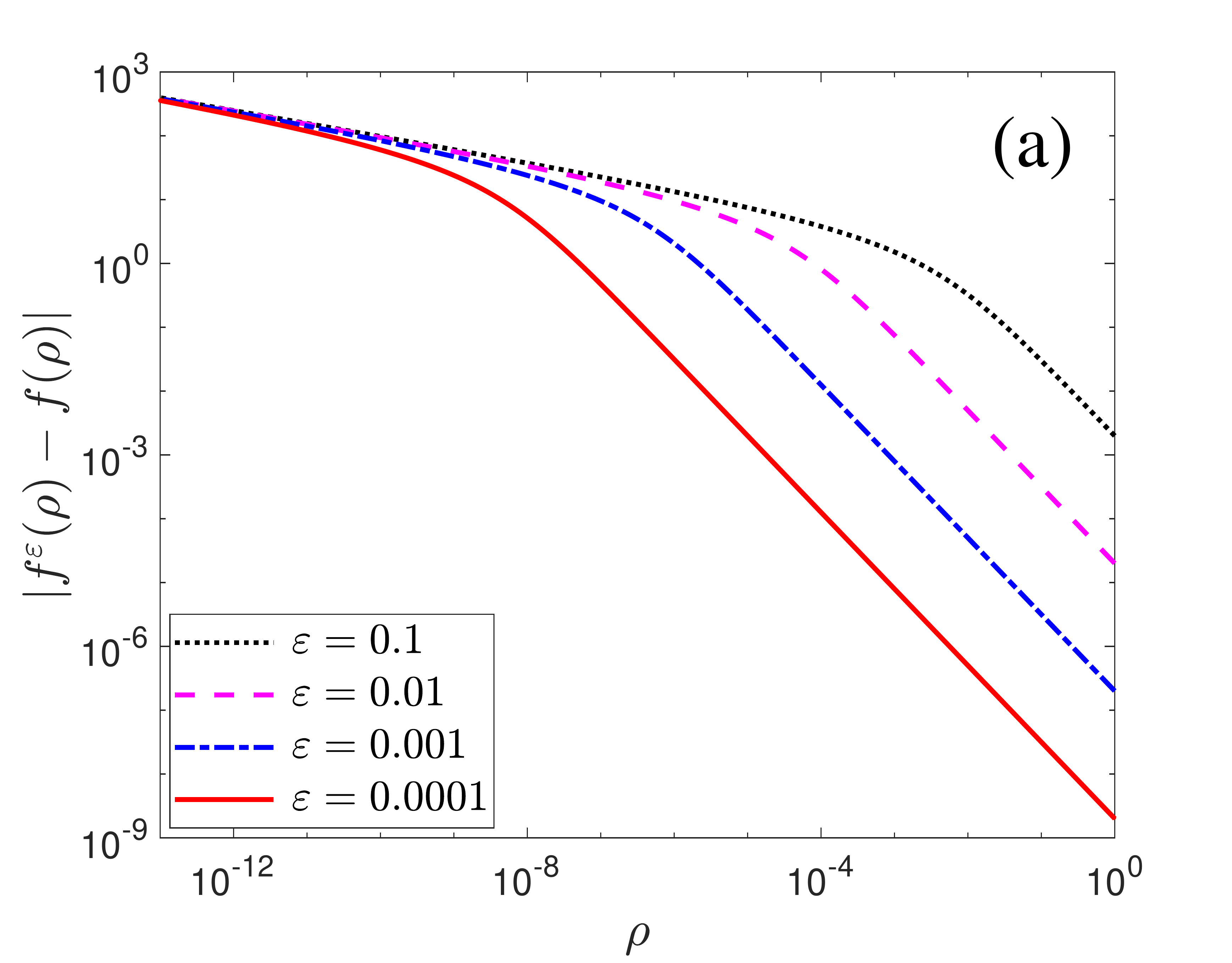}}
\end{minipage}
\begin{minipage}{0.5\textwidth}
\centerline{\includegraphics[width=7.5cm,height=5cm]{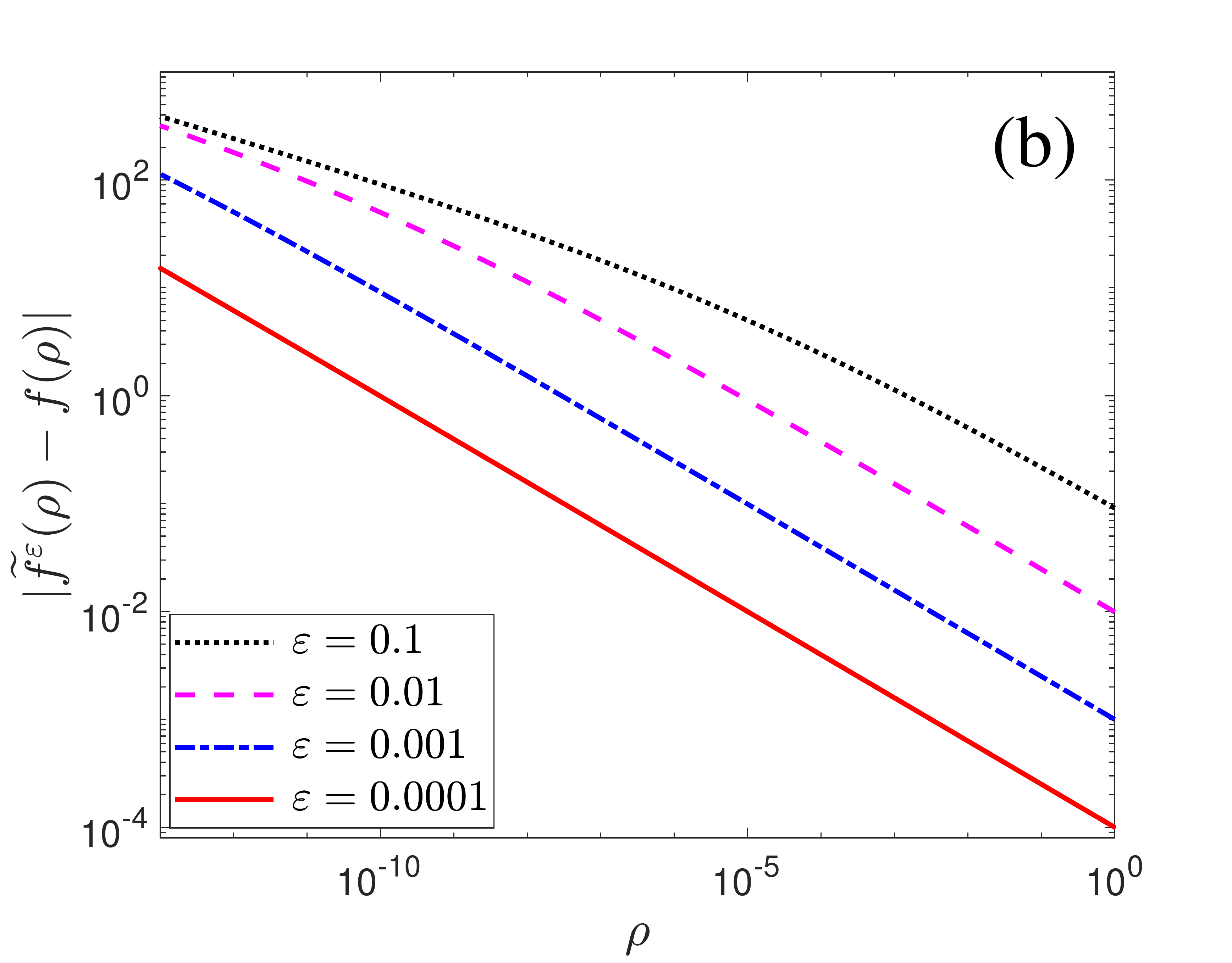}}
\end{minipage}
\begin{minipage}{0.5\textwidth}
\centerline{\includegraphics[width=7.5cm,height=5cm]{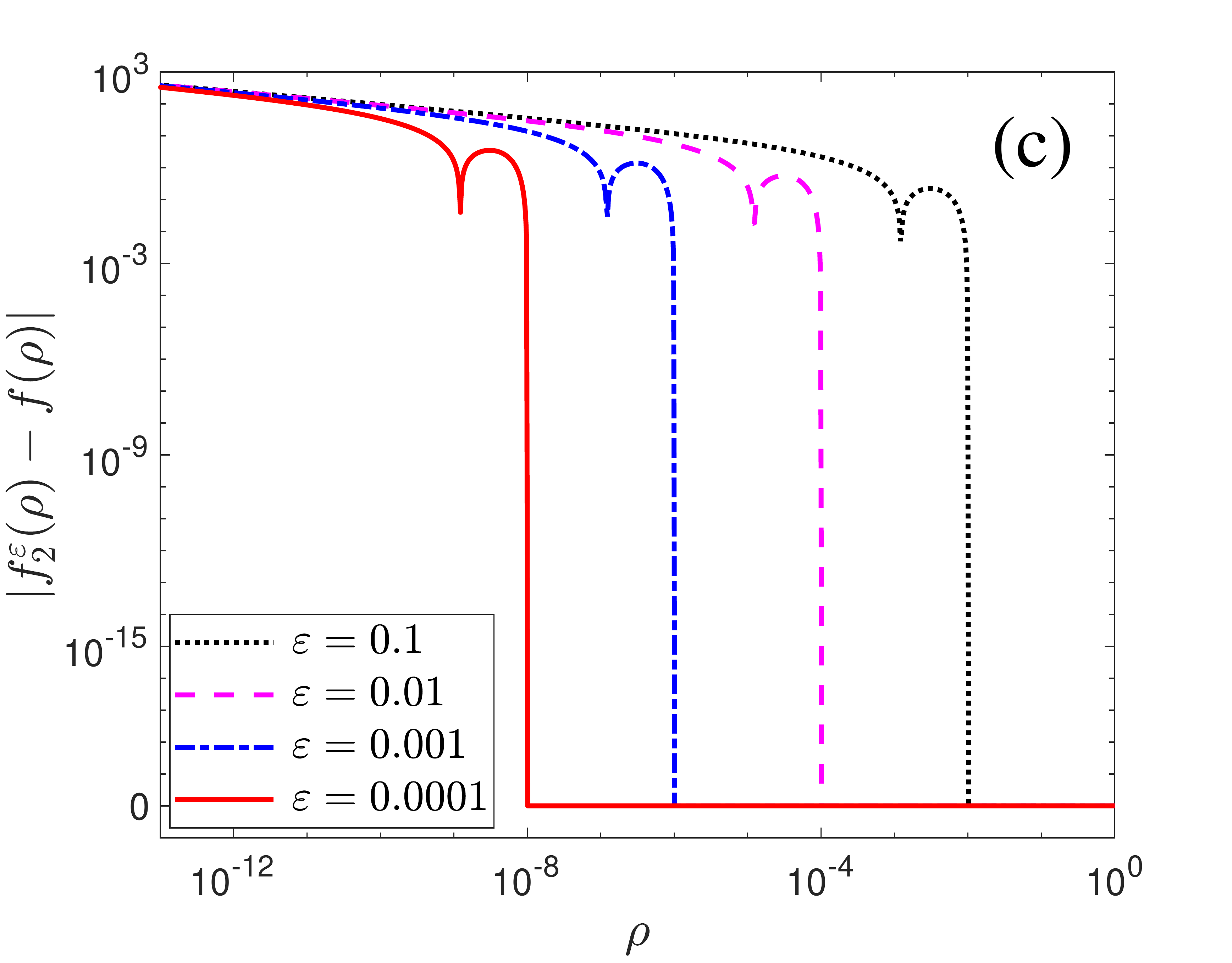}}
\end{minipage}
\begin{minipage}{0.5\textwidth}
\centerline{\includegraphics[width=7.5cm,height=5cm]{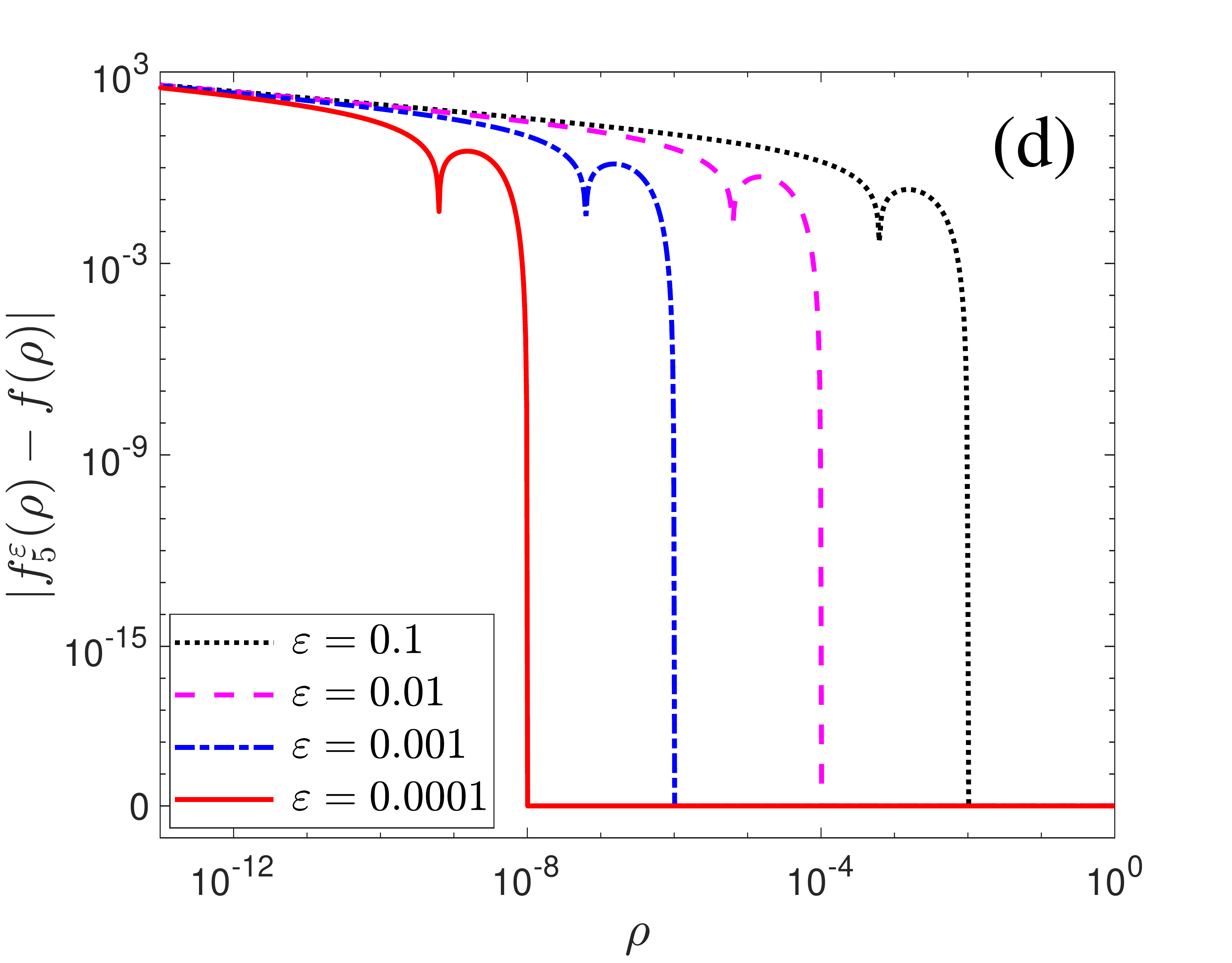}}
\end{minipage}
\caption{Comparison of different regularizations for $f(\rho) = \rho^{\alpha}$ with $\alpha = -0.2$ for different $\eps$: (a) $f^{\eps}(\rho)$; (b) $\widetilde{f}^{\eps}(\rho)$; (c) $f^{\eps}_2(\rho)$; and  (d) $f^{\eps}_5(\rho)$.}
\label{fig:sf_al02}
\end{figure}

\begin{figure}[ht!]
\begin{minipage}{0.5\textwidth}
\centerline{\includegraphics[width=7.5cm,height=5cm]{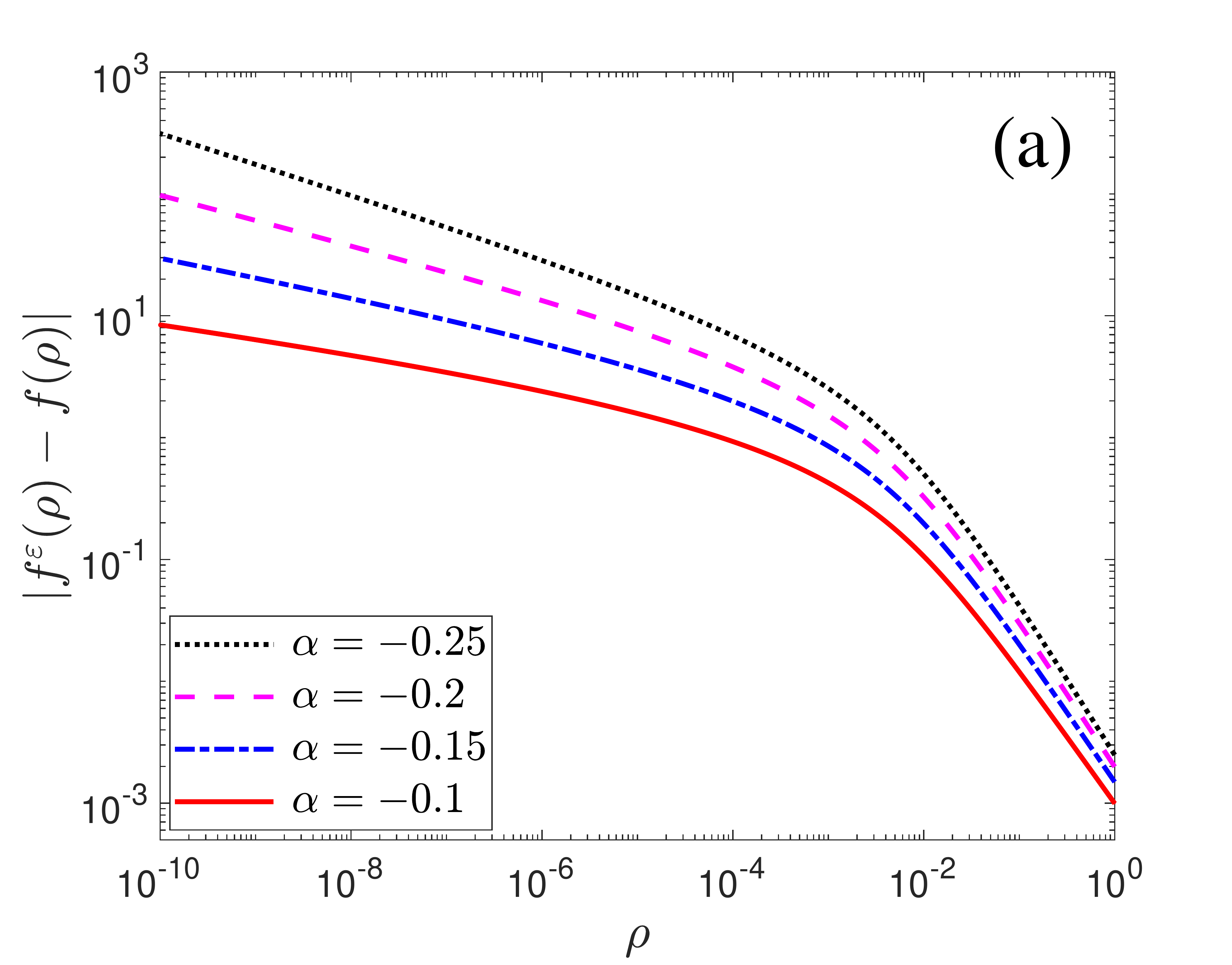}}
\end{minipage}
\begin{minipage}{0.5\textwidth}
\centerline{\includegraphics[width=7.5cm,height=5cm]{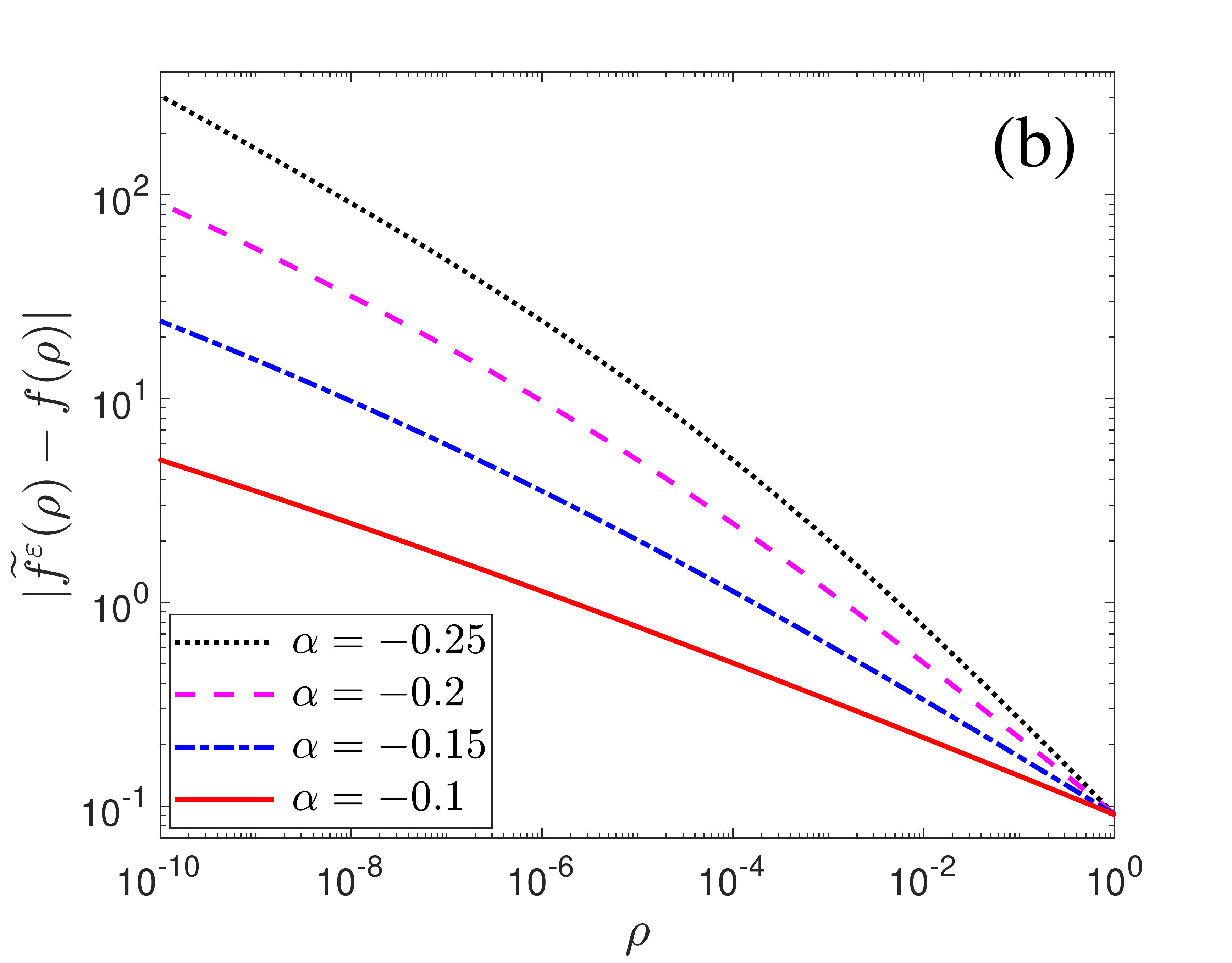}}
\end{minipage}
\begin{minipage}{0.5\textwidth}
\centerline{\includegraphics[width=7.5cm,height=5cm]{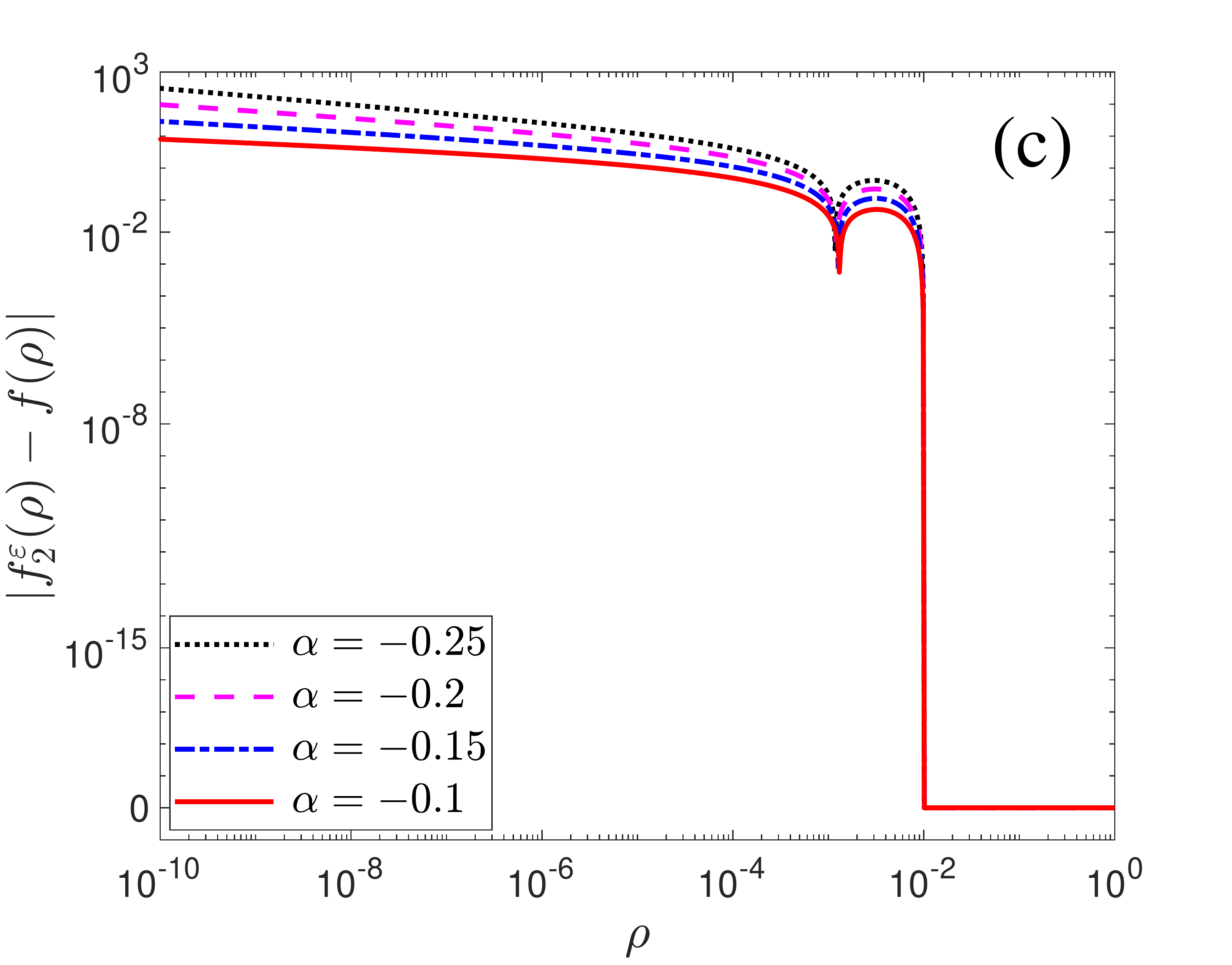}}
\end{minipage}
\begin{minipage}{0.5\textwidth}
\centerline{\includegraphics[width=7.5cm,height=5cm]{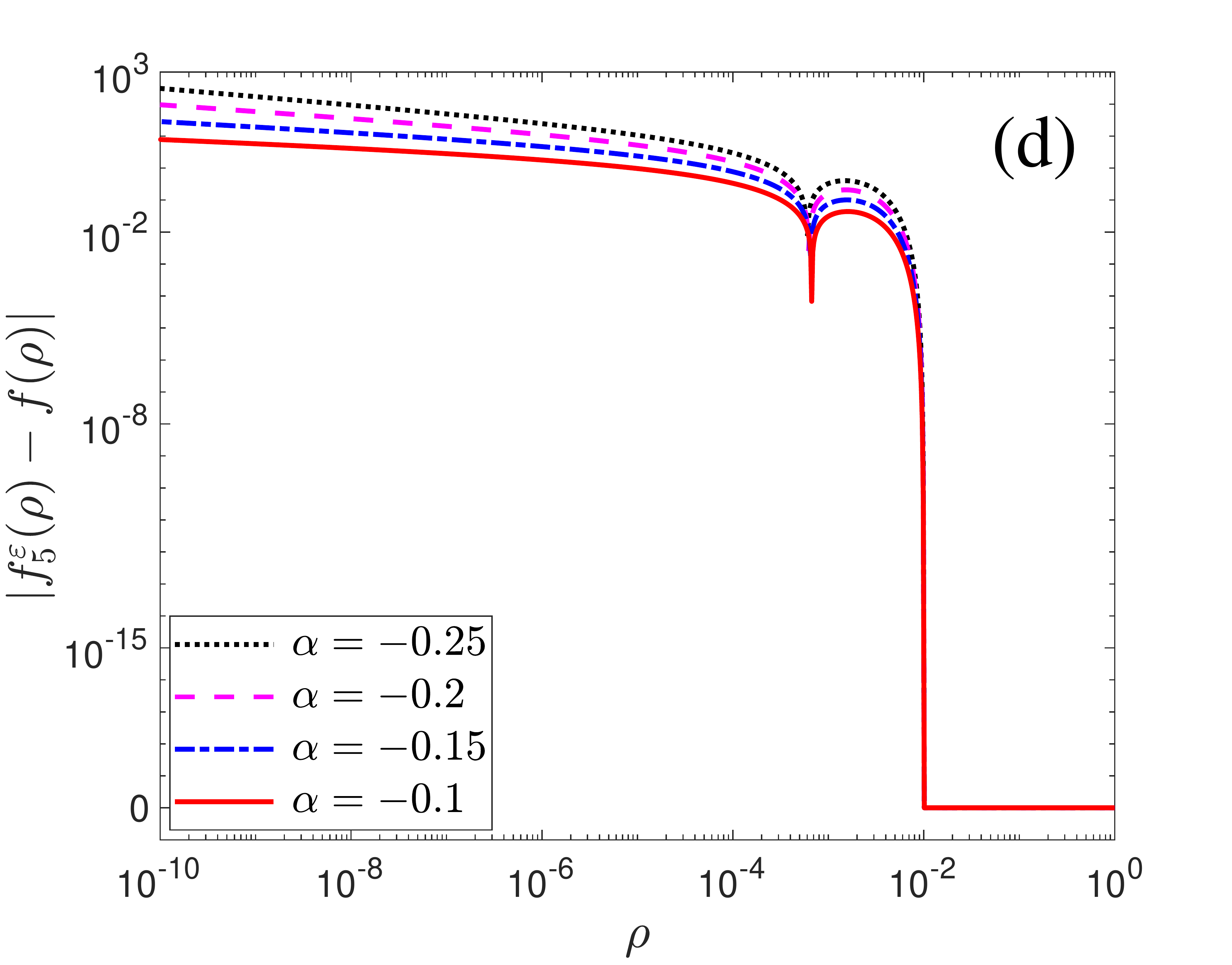}}
\end{minipage}
\caption{Comparison of different regularizations for $f(\rho) = \rho^{\alpha}$ with different $\alpha$  with $\eps = 0.1$ for different $\alpha$: (a) $f^{\eps}(\rho)$; (b) $\widetilde{f}^{\eps}(\rho)$; (c) $f^{\eps}_2(\rho)$; and (d) $f^{\eps}_5(\rho)$.}
\label{fig:sf_difalpha}
\end{figure}

Figures \ref{fig:sf_al02} and \ref{fig:sf_difalpha} display the regularizations $f^{\eps}$, $\widetilde{f}^{\eps}$ and $f^{\eps}_n$ $(n = 2, 5)$ for different $\eps$ and $\alpha$, respectively. From these figures, we find the local regularization $f^{\eps}_n$ approximates the singular nonlinearity $f(\rho) = \rho^{\alpha}$ more accurately than the other two global regularizations.

\section{Numerical methods for the regularized models}
In this section, we discretize the regularized models by the first-order Lie-Trotter time-splitting (LTTS) \cite{GL,MQ,HFT} and Lawson-type exponential integrator (LTEI) \cite{HO,JDL,OSU} in time. We impose the periodic boundary condition on the domain $\Omega$ and apply the Fourier pseudospectral method for spatial discretization \cite{BCST2,BJM,BJM2}. For simplicity, we only present the numerical schemes in one dimension (1D) and generalization to higher dimensions is straightforward. In 1D, the regularized models can be written in the general form as
\begin{equation}
\left\{
\begin{aligned}
&i\partial_t \psi^{\eps}(x, t) = -\Delta \psi^{\eps}(x, t) +\lambda\,  f^{\eps}_{\rm Reg}(|\psi^{\eps}(x, t)|^2)\psi^{\eps}(x, t), \quad x \in \Omega = (a, b), \ t>0,	\\
&\psi^{\eps}(x, 0) = \psi_0(x), \quad x \in \overline{\Omega},
\end{aligned}\right.
\label{eq:GNLSE}
\end{equation}
where $f^{\eps}_{\rm Reg}$ is the regularized nonlinear function, i.e., $f^{\eps}_{\rm Reg}(\rho) =(\rho+\eps^2)^{\alpha}$ for the rNLSE \eqref{eq:RNLSE1},   $f^{\eps}_{\rm Reg}(\rho) =1/(\rho^{-\alpha}+\eps)$ for the rNLSE \eqref{eq:RNLSE2}, and $f^{\eps}_{\rm Reg}(\rho) = f^{\eps}_n(\rho)$ for the erNLSE \eqref{eq:ERNLSE}. Then, we will derive the semi-discretization for \eqref{eq:GNLSE} in time via the first-order Lie-Trotter time-splitting and Lawson-type exponential integrator and the full-discretization by combining with the Fourier pseudospectral method in space.

\subsection{The Lie-Trotter time-splitting}
Splitting methods for the time integrator of \eqref{eq:GNLSE} are based on a decomposition of the flow. More precisely, let us define the flow $\Phi^t_A$ of the linear Schr\"odinger equation
\begin{equation}
\left\{
\begin{aligned}
&i\partial_t v(x, t) = -\Delta v(x, t), \quad x \in \Omega, \ t>0,	\\
&v(x, 0) =v_0(x),
\end{aligned}\right.
\label{eq:sub1}
\end{equation}
and the flow $\Phi^t_B$ for the differential equation
\begin{equation}
\left\{
\begin{aligned}
&i\partial_t w(x, t) =  \lambda f^{\eps}_{\rm Reg}(|w(x, t)|^2)w(x, t), \quad x \in \Omega, \ t>0,	\\
&w(x, 0) =w_0(x).
\end{aligned}\right.
\label{eq:sub2}
\end{equation}
The associated evolution operators are given by
\begin{equation}
v(x, t) = \Phi^t_A(v_0(x)) = e^{it\Delta}v_0(x), \quad w(x, t) = \Phi^t_B(w_0(x))= w_0(x) e^{-it \lambda f^{\eps}_{\rm Reg}(|w_0(x)|^2)}, \quad t \geq 0.	
\end{equation}
The idea of splitting methods is to approximate the flow of \eqref{eq:GNLSE} by combining the two flows $\Phi^t_A$ and $\Phi^t_B$ \cite{HFT,MQ}. Let $\tau = \Delta t > 0$ be the time step size and  $t_m = m\tau$ ($m = 0, 1, 2, \ldots$) as the time steps. Denote $\psi^{[\eps, m]} := \psi^{[\eps, m]}(x)$ as the approximation of $\psi^{\eps}(x, t_m)$, then the first-order Lie-Trotter time-splitting (LTTS) method can be written as
\begin{equation}
\psi^{[\eps, m+1]}	= \Phi^{\tau}(\psi^{[\eps, m]}) := \Phi^{\tau}_B\left(\Phi^{\tau}_A(\psi^{[\eps, m]})\right),\ m \geq 0; \quad \psi^{[\eps, 0]} = \psi_0.
\label{eq:LTTS}
\end{equation}

Then, we apply the Fourier pseudospectral method in space to derive a full-discretization.  Let $N$ be an positive integer, and define the spatial mesh size $h = (b-a)/N$, then the grid points are chosen as
\begin{equation*}
x_j := a + jh, \quad j \in \mathcal{T}^0_N	 = \{j~|~j = 0, 1, \ldots, N\}.
\end{equation*}
Denote $X_N := \{\psi= (\psi_0, \psi_1, \ldots, \psi_N)^T \in \mathbb{C}^{N+1} \ \vert \ \psi_0 = \psi_N\}$, $C_{\rm per}(\Omega)=\{\psi \in C(\overline \Omega) \ \vert \ \psi(a) = \psi(b)\}$ and
\begin{align*}
& Y_N := \text{span}\left\{e^{i\mu_l(x-a)},\ x \in \overline{\Omega}, \ l \in \mathcal{T}_N\right\},\; \mathcal{T}_N = \left\{l \ \vert \ l = -\frac{N}{2}, \ldots, \frac{N}{2}-1\right\},
\end{align*}
where $\mu_l=\frac{2\pi l}{b-a}$. For any $\psi(x) \in C_{\rm per}(\Omega)$ and a vector $\psi \in X_N$, let $P_N: L^2(\Omega) \to Y_N$ be the standard $L^2$-projection operator onto $Y_N$, $I_N : C_{\rm per}(\Omega) \to Y_N$ or $I_N : X_N \to Y_N$ be the trigonometric interpolation operator \cite{STL}, i.e.,

\begin{equation*}
P_N \psi  = \sum_{l \in \mathcal{T}_N} \widehat{\psi}_l e^{i\mu_l(x-a)},\quad I_N \psi  = \sum_{l \in \mathcal{T}_N} \widetilde{\psi}_l e^{i\mu_l(x-a)},\quad x \in \overline{\Omega},
\end{equation*}
where
\begin{equation*}
\widehat{\psi}_l = \frac{1}{b- a}\int^{b}_{a} \psi(x) e^{-i\mu_l (x-a)} dx, \quad \widetilde{\psi}_l = \frac{1}{N}\sum_{j=0}^{N-1} \psi_j e^{-i\mu_l (x_j-a)}, \quad l \in \mathcal{T}_N,
\end{equation*}
with $\psi_j$ interpreted as $\psi(x_j)$ when involved.

Let $\psi^{\eps, m}_j$ be the numerical approximation of $\psi^{\eps}(x_j, t_m)$ for $j \in \mathcal{T}^0_N	$, and $m \geq 0$. Denote $\psi^{\eps, m} = (\psi^{\eps, m}_0, \psi^{\eps, m}_1, \ldots, \psi^{\eps, m}_N)^T \in \mathbb{C}^{N+1}$ for $m = 0, 1, \ldots$. Then, the Lie-Trotter time-splitting Fourier pseudospectral (TSFP) discretization for the rNLSE \eqref{eq:GNLSE} can be written as
\begin{equation}
\label{eq:TSFP}
\begin{split}
&\psi^{(\eps, 1)}_j=\sum_{l \in \mathcal{T}_N} e^{-i\tau\mu^2_l} \widetilde{(\psi^{\eps, m})}_l e^{i\mu_l(x_j-a)}, \\
&\psi^{\eps, m+1}_j=e^{- i\tau \lambda f^{\eps}_{\rm Reg}(|\psi^{(\eps, 1)}_j|^2)}\psi^{(\eps, 1)}_j, \qquad j\in \mathcal{T}^0_N, \quad m \ge 0, \\
\end{split}
\end{equation}
where $\psi^{\eps, 0}_j = \psi_0(x_j)$ for $j\in \mathcal{T}^0_N$.

This full discretization is explicit, time-symmetric, time transverse invariant, and it conserves mass in the discretized level.

\subsection{The Lawson-type exponential integrator}
Exponential integrators for the flow \eqref{eq:GNLSE} are constructed by incorporating the exact propagator of the linear part in an appropriate way \cite{HO}. By Duhamel's formula (variation-of-constants formula), the exact solution of \eqref{eq:GNLSE} at $t = t_m+\tau$ is given by
\begin{equation}
\psi^{\eps}(t_m + \tau) = e^{it\Delta}\psi(t_m) - i\lambda \int^{\tau}_0 e^{i(\tau - s)\Delta}\left[f^{\eps}_{\rm Reg}(|\psi^{\eps}(t_m+s)|^2) \psi^{\eps}(t_m+s)\right] ds,
\end{equation}
where we denote by $\psi^{\eps}(t) = \psi^{\eps}(x, t)$ in short. By applying the approximation
\begin{equation}
\psi(t_m+s) \approx \psi(t_m)	
\end{equation}
in the integral terms and combining the first-order Lawson method \cite{JDL}, we obtain the semi-discretization via the Lawson-type exponential integrator (LTEI) method as
\begin{equation}
\psi^{[\eps, m+1]} =  e^{it\Delta}\psi^{[\eps, m]} - i \lambda \tau e^{it\Delta}\left[f^{\eps}_{\rm Reg}(|\psi^{[\eps, m]}|^2) \psi^{[\eps, m]}\right], \ m \geq 0; \quad  \psi^{[\eps, 0]} = \psi_0.
\label{eq:LTEI}
\end{equation}
Respectively, the exponential integrator Fourier pseudospectral (EIFP) method for the rNLSE \eqref{eq:GNLSE} is
\begin{equation}
\psi^{\eps, m+1}_j =  \sum_{l \in \mathcal{T}_N} e^{-i\tau\mu_l^2}\widetilde{(\psi^{\eps, m})}_l - i \lambda \tau  \sum_{l \in \mathcal{T}_N} e^{-i\tau\mu_l^2}\widetilde{g(\psi^{\eps, m})}_l, \qquad j\in \mathcal{T}^0_N, \quad n \ge 0,
\label{eq:EIFP}	
\end{equation}
where $g(\psi^{\eps, m}) = f^{\eps}_{\rm Reg}(|\psi^{\eps, m}|^2) \psi^{\eps, m}$ and $\psi^{\eps, 0}_j = \psi_0(x_j)$ for $j\in \mathcal{T}^0_N$.

Again, this full discretization is explicit, but it is not time-symmetric and time transverse invariant. Also, it does not conserve mass in the discretized level.

\begin{remark}
In this section, we apply the first-order Lie-Trotter time-splitting method and Lawson-type exponential integrator to discretize the regularized models \eqref{eq:GNLSE} in time. It can be extended to the second-order temporal semi-discretizations via the Strang time-splitting method \cite{MQ} and Lawson-type exponential integrator scheme \cite{JDL}.
\end{remark}

\section{Numerical results}
In this section, we first compare the convergence rate of the local energy regularized model \eqref{eq:ERNLSE} with the global regularized models \eqref{eq:RNLSE1} and \eqref{eq:RNLSE2}. Then, we test the order of accuracy of the TSFP method \eqref{eq:TSFP} and the EIFP method  \eqref{eq:EIFP}. We take $\lambda=1$ except stated otherwise. We consider the Gaussian initial data:
\begin{equation}
\psi_0(x) = 	\frac{1}{\pi^{1/4}} e^{-x^2/2}, \quad x \in \mathbb{R}.
\end{equation}
The sNLSE \eqref{eq:NLSE1} is numerically solved on the domain $\Omega = [-16, 16]$, which is large enough such that the truncation error is negligible. To quantify the numerical errors, we introduce the following error functions:
\begin{align*}
&\widehat{e}^{\eps}_{\rho}(t_m) := \rho(\cdot, t_m) - \rho^{\eps}(\cdot, t_m) = |\psi(\cdot, t_m)|^2- |\psi^{\eps}(\cdot, t_m)|^2\\
&\widehat{e}^{\eps}(t_m) := \psi(\cdot, t_m)- \psi^{\eps}(\cdot, t_m),  \quad \widetilde{e}^{\eps}(t_m) := \psi(\cdot, t_m)- \psi^{\eps, m},\\
&{e}^{\eps}(t_m) :=  \psi^{\eps}(\cdot, t_m)-\psi^{\eps, m}, \quad e^{\eps}_{\rm E} = |E(\psi_0)-E^{\eps}_{\rm Reg}(\psi_0)|,
\end{align*}
where $\psi$ and $\psi^{\eps}$ are the exact solutions of the sNLSE \eqref{eq:NLSE1} and the regularized equation \eqref{eq:GNLSE}, respectively, and $\psi^{\eps, n}$ is the numerical solution of the regularized equation \eqref{eq:GNLSE} by the TSFP method \eqref{eq:TSFP} or the EIFP method \eqref{eq:EIFP}. The 'exact' solution $\psi^{\eps}$ is numerically obtained by the TSFP scheme with a very fine mesh size $h = 1/16$ and a very small time step $\tau = 10^{-5}$. Similarly, the 'exact' solution $\psi$ is numerically computed by the TSFP scheme with a very fine mesh size and a very small time step as well as a very small regularization parameter $\eps = 10^{-12}$.

\begin{figure}[ht!]
\begin{minipage}{0.5\textwidth}
\centerline{\includegraphics[width=7.5cm,height=5cm]{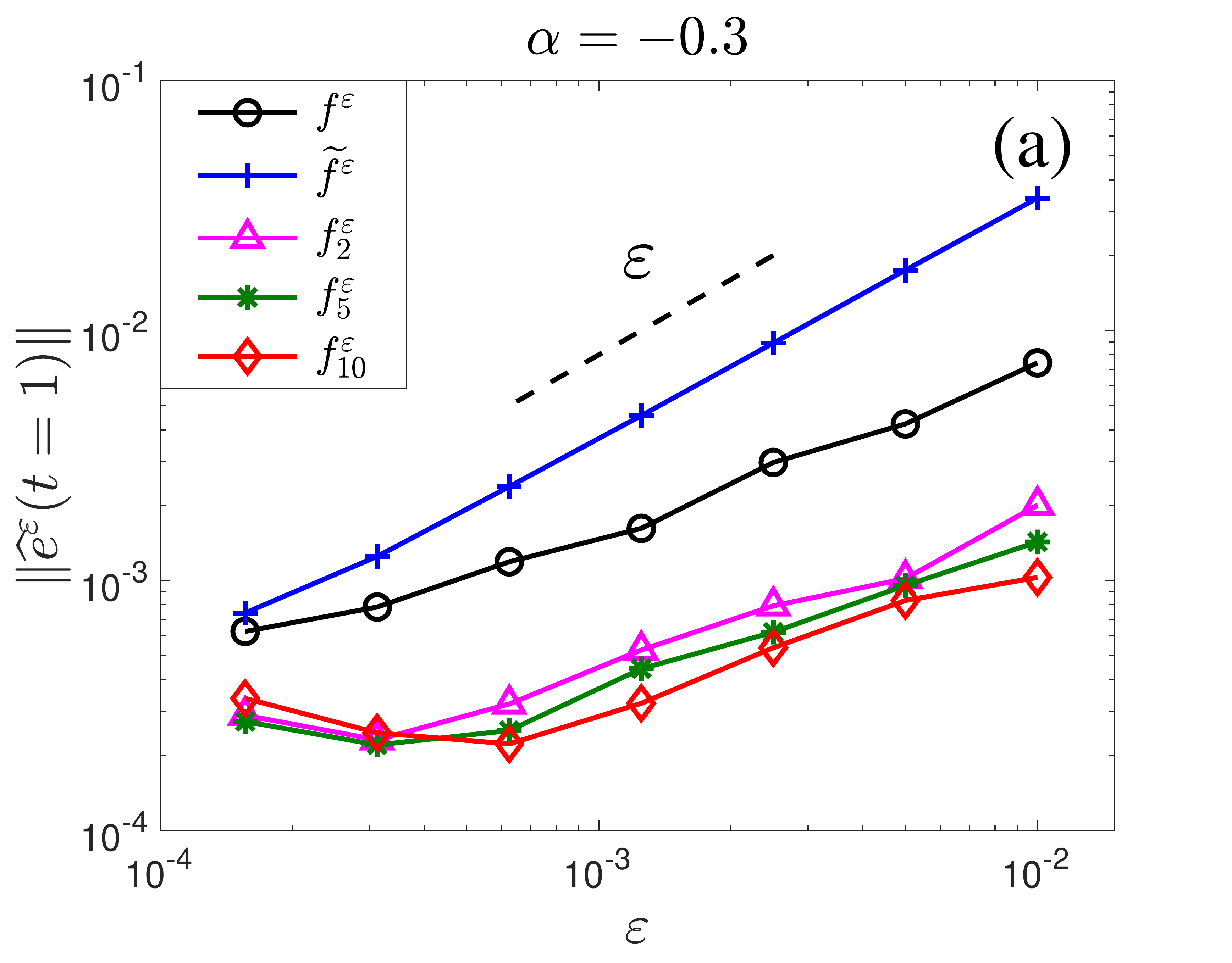}}
\end{minipage}
\begin{minipage}{0.5\textwidth}
\centerline{\includegraphics[width=7.5cm,height=5cm]{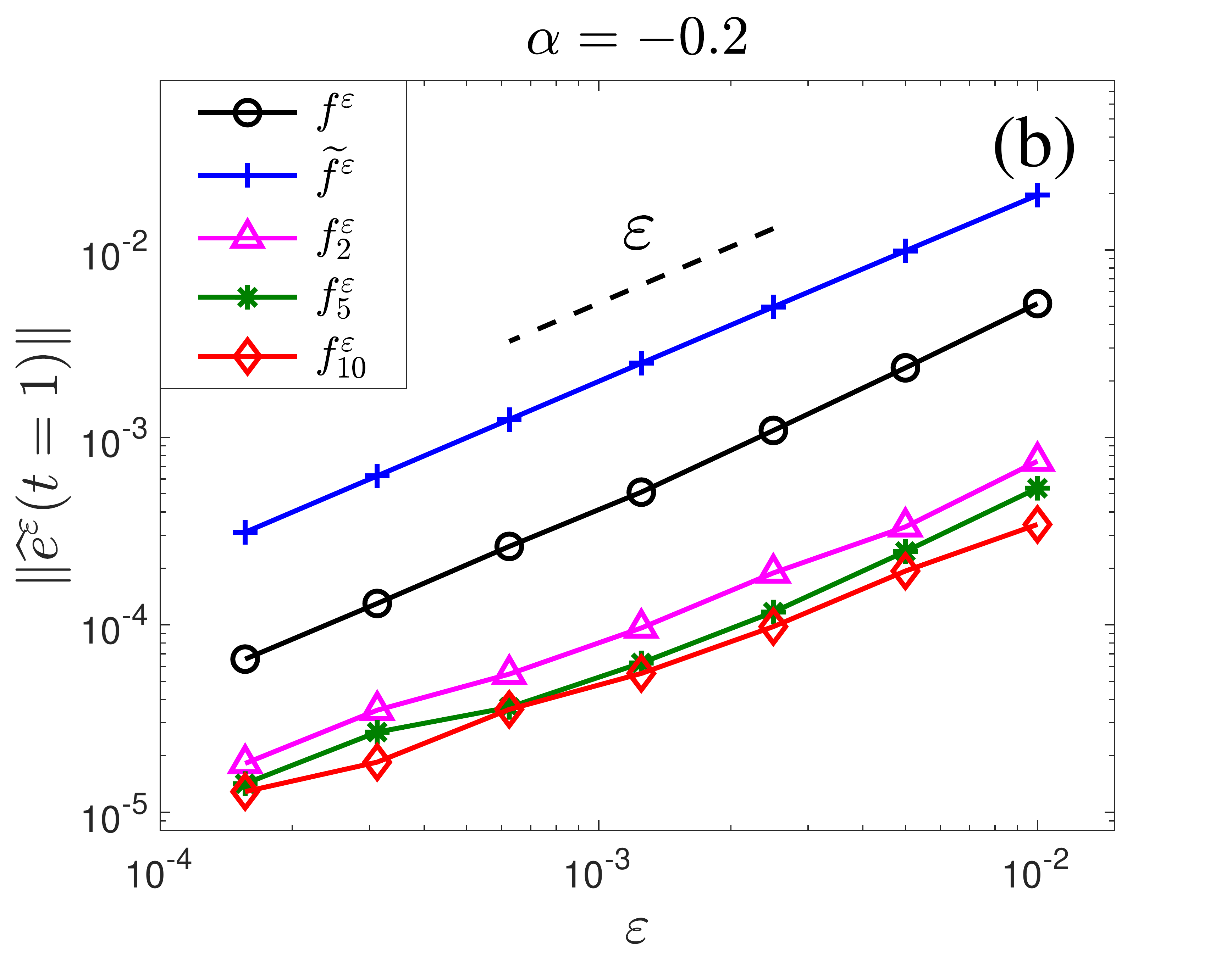}}
\end{minipage}
\caption{Convergence of the rNLSE \eqref{eq:GNLSE} with different regularized nonlinearity $f^{\eps}_{\rm Reg}$ to the sNLSE \eqref{eq:NLSE1}: errors $\left\|\widehat{e}^{\eps}(t=1)\right\|$ for (a) $\alpha = -0.3$; and (b) $\alpha = -0.2$.}
\label{fig:model}
\end{figure}

\begin{figure}[ht!]
\begin{minipage}{0.5\textwidth}
\centerline{\includegraphics[width=7.5cm,height=5cm]{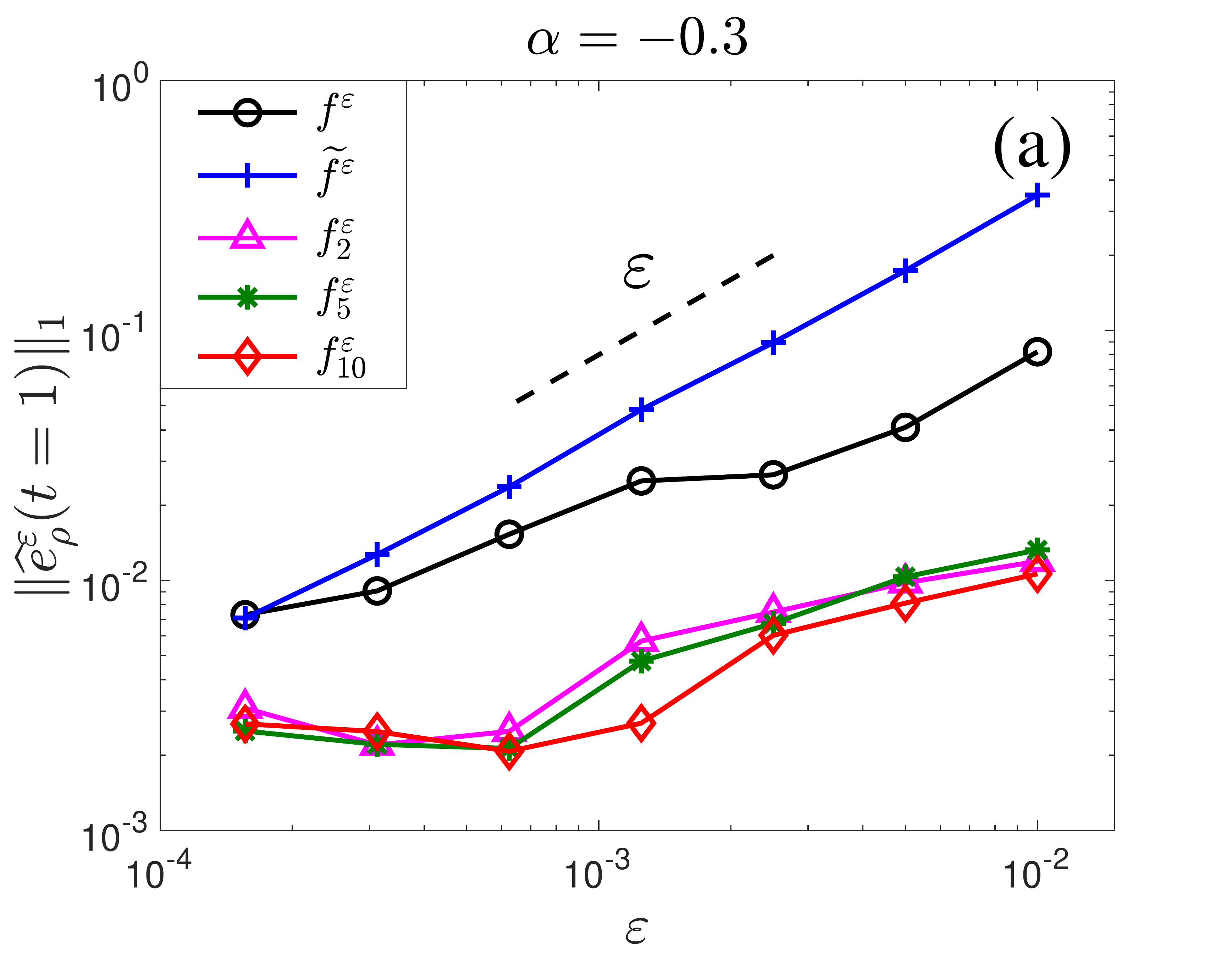}}
\end{minipage}
\begin{minipage}{0.5\textwidth}
\centerline{\includegraphics[width=7.5cm,height=5cm]{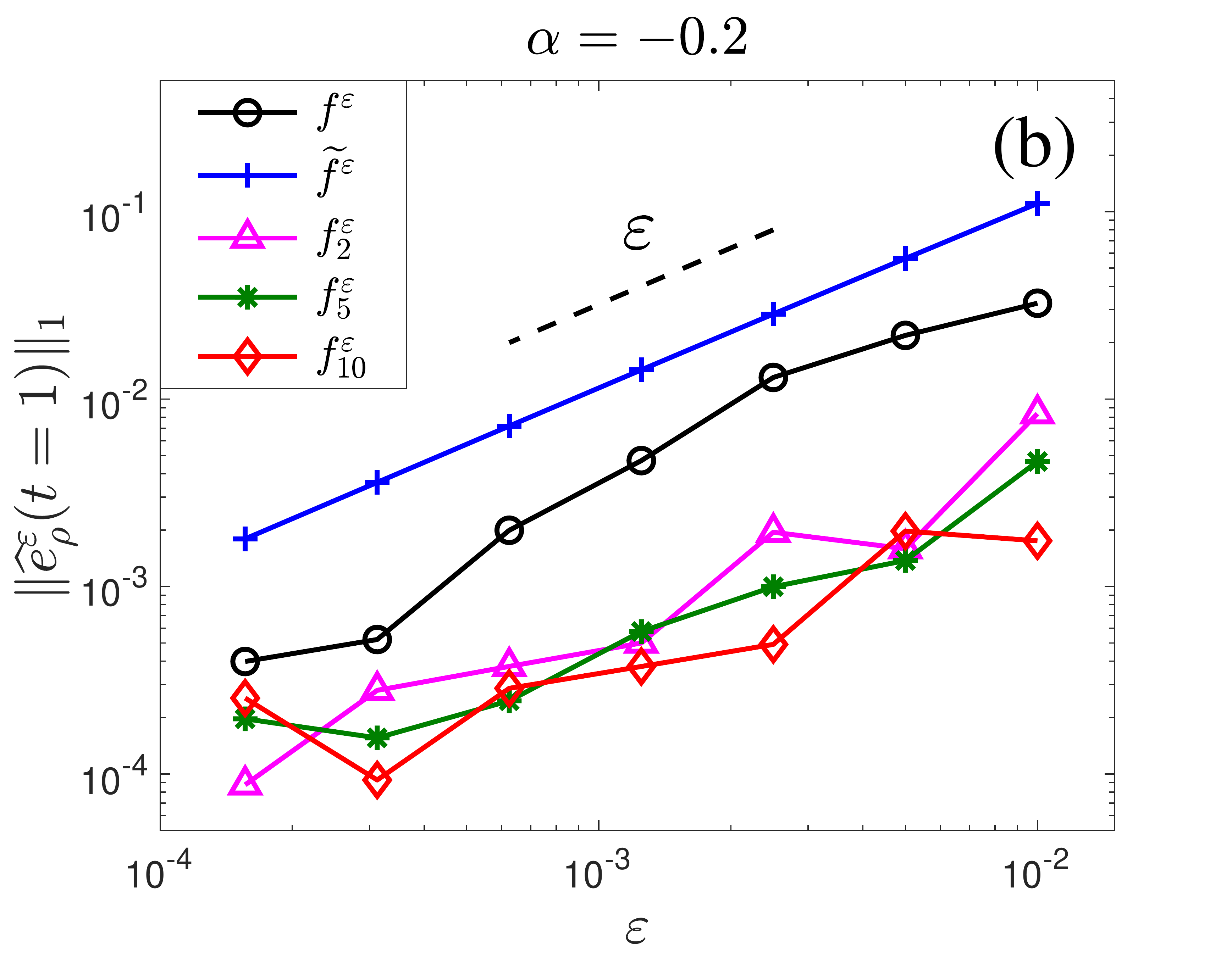}}
\end{minipage}
\caption{Convergence of the rNLSE \eqref{eq:GNLSE} with different regularized nonlinearity $f^{\eps}_{\rm Reg}$ to the sNLSE \eqref{eq:NLSE1}: errors $\left\|\widehat{e}^{\eps}_{\rho}(t=1)\right\|_1$ for (a) $\alpha = -0.3$; and (b) $\alpha = -0.2$.}
\label{fig:rho}
\end{figure}

\begin{figure}[ht!]
\begin{minipage}{0.5\textwidth}
\centerline{\includegraphics[width=7.5cm,height=5cm]{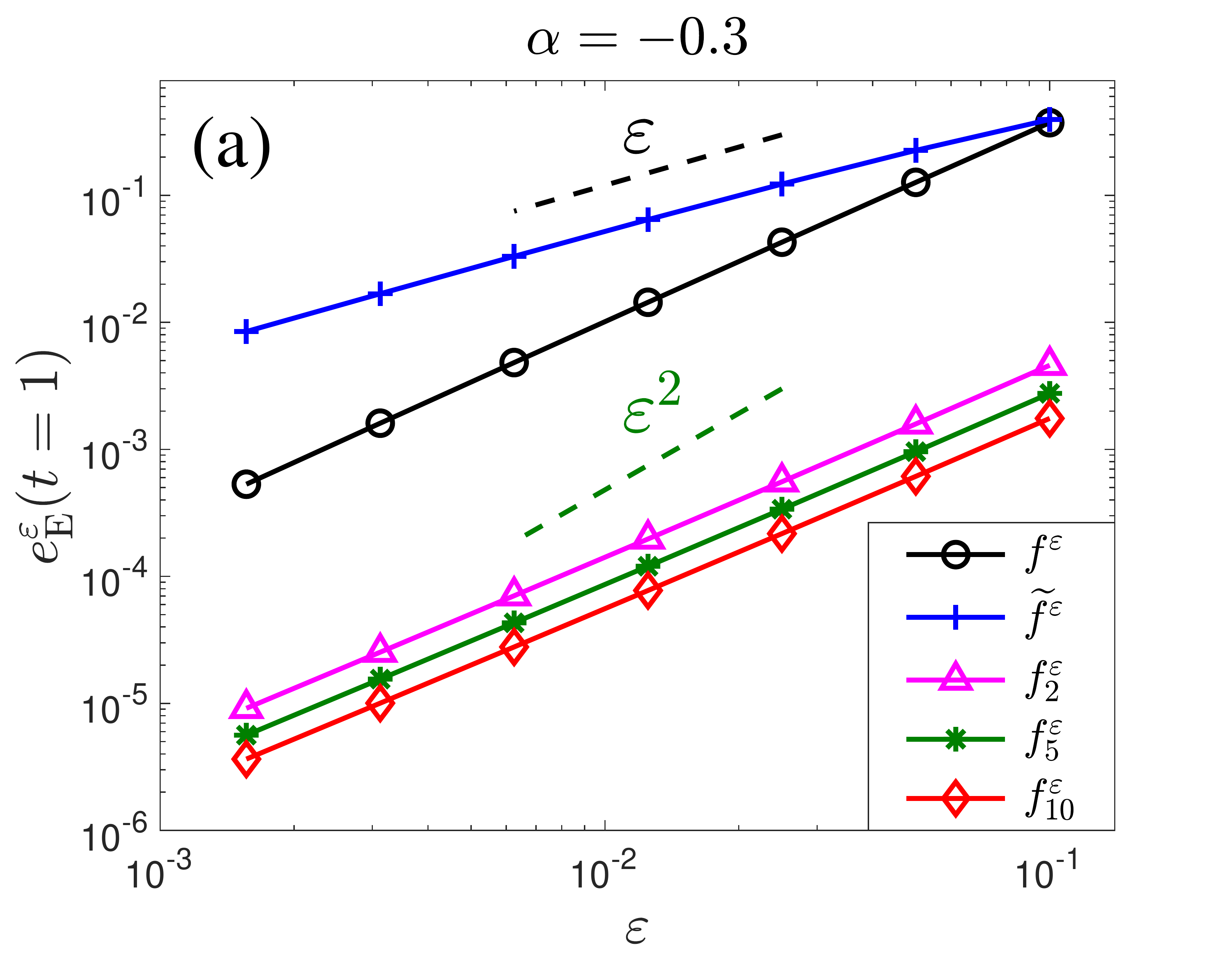}}
\end{minipage}
\begin{minipage}{0.5\textwidth}
\centerline{\includegraphics[width=7.5cm,height=5cm]{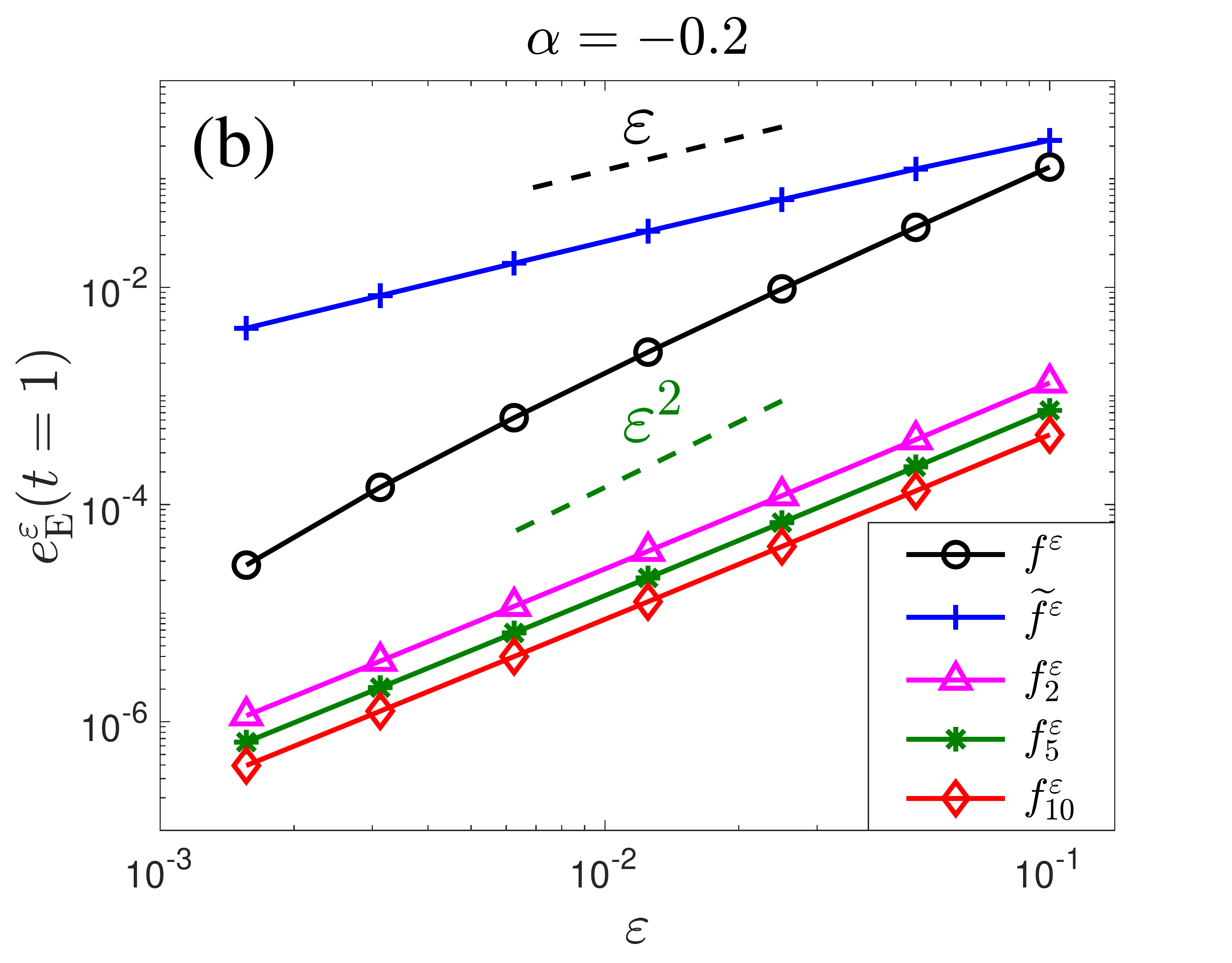}}
\end{minipage}
\caption{Convergence of the rNLSE \eqref{eq:GNLSE} with different regularized nonlinearity $f^{\eps}_{\rm Reg}$ to the sNLSE \eqref{eq:NLSE1}: errors $e^{\eps}_{\rm E}(t=1)$ for (a) $\alpha = -0.3$; and (b) $\alpha = -0.2$.}
\label{fig:Energy}
\end{figure}

\subsection{Convergence rates of different regularized  models}
First, we test the errors between the solutions of the sNLSE  \eqref{eq:NLSE1} and the regularized model \eqref{eq:GNLSE} with different choices of regularized nonlinearity $f^{\eps}_{\rm Reg}$. Figures \ref{fig:model}, \ref{fig:rho} and \ref{fig:Energy} show the errors $\left\|\widehat{e}^{\eps}(t=1)\right\|$,   $\left\|\widehat{e}^{\eps}_{\rho}(t=1)\right\|_1$ and $e^{\eps}_{\rm E}(t=1)$, respectively.

From these figures and additional numerical results not shown here for brevity, we have the following observations: (i) The solution of the regularized NLSE \eqref{eq:RNLSE2} converges linearly to that of the sNLSE \eqref{eq:NLSE1} in terms of $\eps$, while the regularized NLSE \eqref{eq:RNLSE1} and the energy regularized model \eqref{eq:ERNLSE} converge much slower and the convergence rate depends on the parameter $\alpha$. The amplitudes of the errors from the local energy regularization are much smaller than the two global regularizations. (ii) For errors of the density in $L^1$-norm, the density of the rNLSE \eqref{eq:RNLSE2} converges linearly to that of the sNLSE \eqref{eq:NLSE1} in terms of $\eps$, while the convergence rate of the rNLSE \eqref{eq:RNLSE1} and the erNLSE \eqref{eq:ERNLSE} depends on the parameter $\alpha$. (iii) For the energy of different regularizations, the energy $\widetilde E^{\eps}$ converges linearly, while the convergence rate of $E^{\eps}$ and $E^{\eps}_n$ (for any $n \geq 1$) depending on $\alpha$ is between linearly and quadratically. In addition, all these three regularized energies are smaller when $\alpha$ is larger. (iv) The proposed energy regularized model (i.e., $f^{\eps}_{\rm Reg} = f^{\eps}_n$ in the rNLSE \eqref{eq:GNLSE}) performs better than the other two regularizations  (i.e., $f^{\eps}_{\rm Reg} = f^{\eps}$ and $f^{\eps}_{\rm Reg} = \widetilde{f}^{\eps}$) in the sense that its corresponding errors in wave function, density and energy are much smaller.

\begin{figure}[t!]
\begin{minipage}{0.5\textwidth}
\centerline{\includegraphics[width=7.5cm,height=5cm]{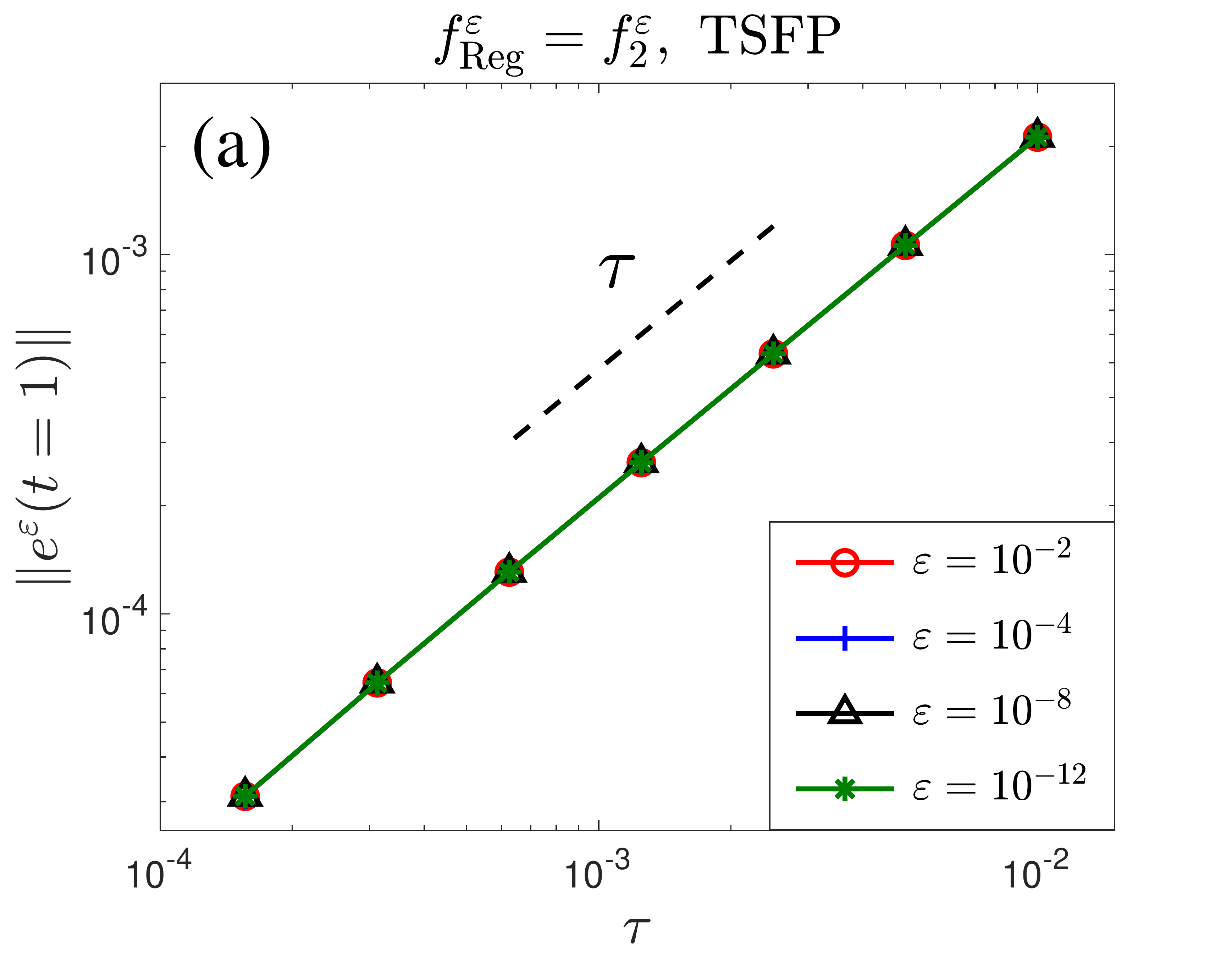}}
\end{minipage}
\begin{minipage}{0.5\textwidth}
\centerline{\includegraphics[width=7.5cm,height=5cm]{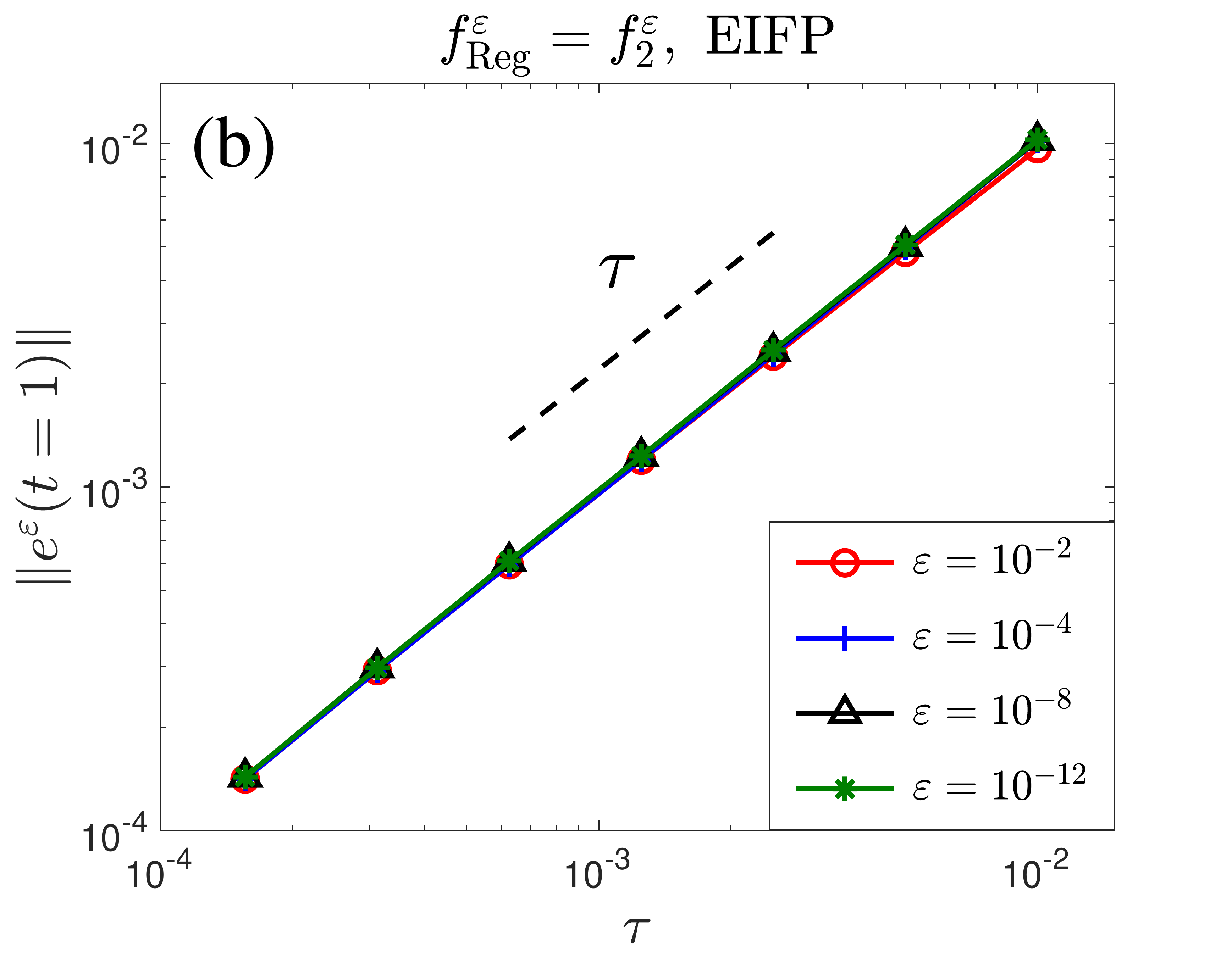}}
\end{minipage}
\begin{minipage}{0.5\textwidth}
\centerline{\includegraphics[width=7.5cm,height=5cm]{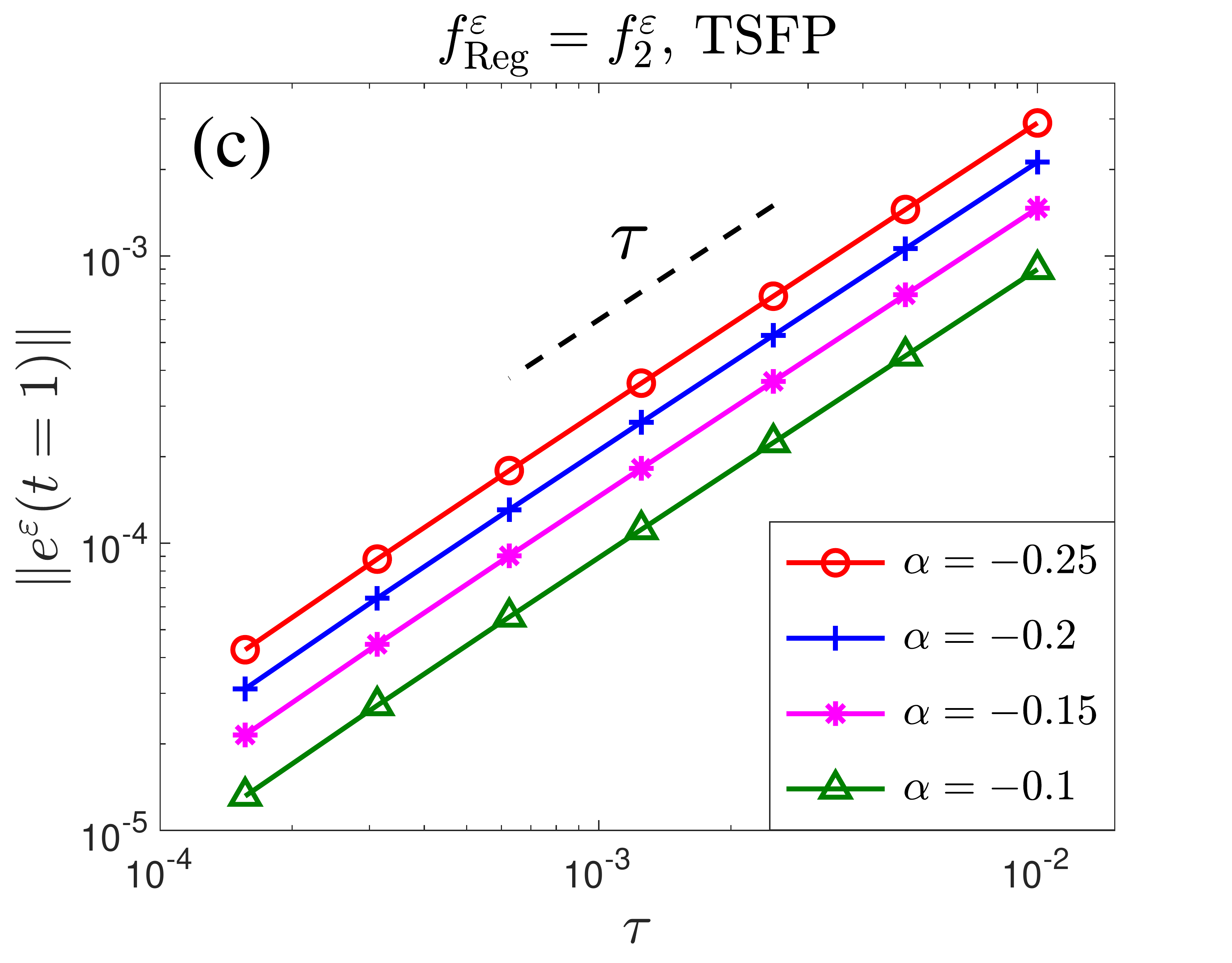}}
\end{minipage}
\begin{minipage}{0.5\textwidth}
\centerline{\includegraphics[width=7.5cm,height=5cm]{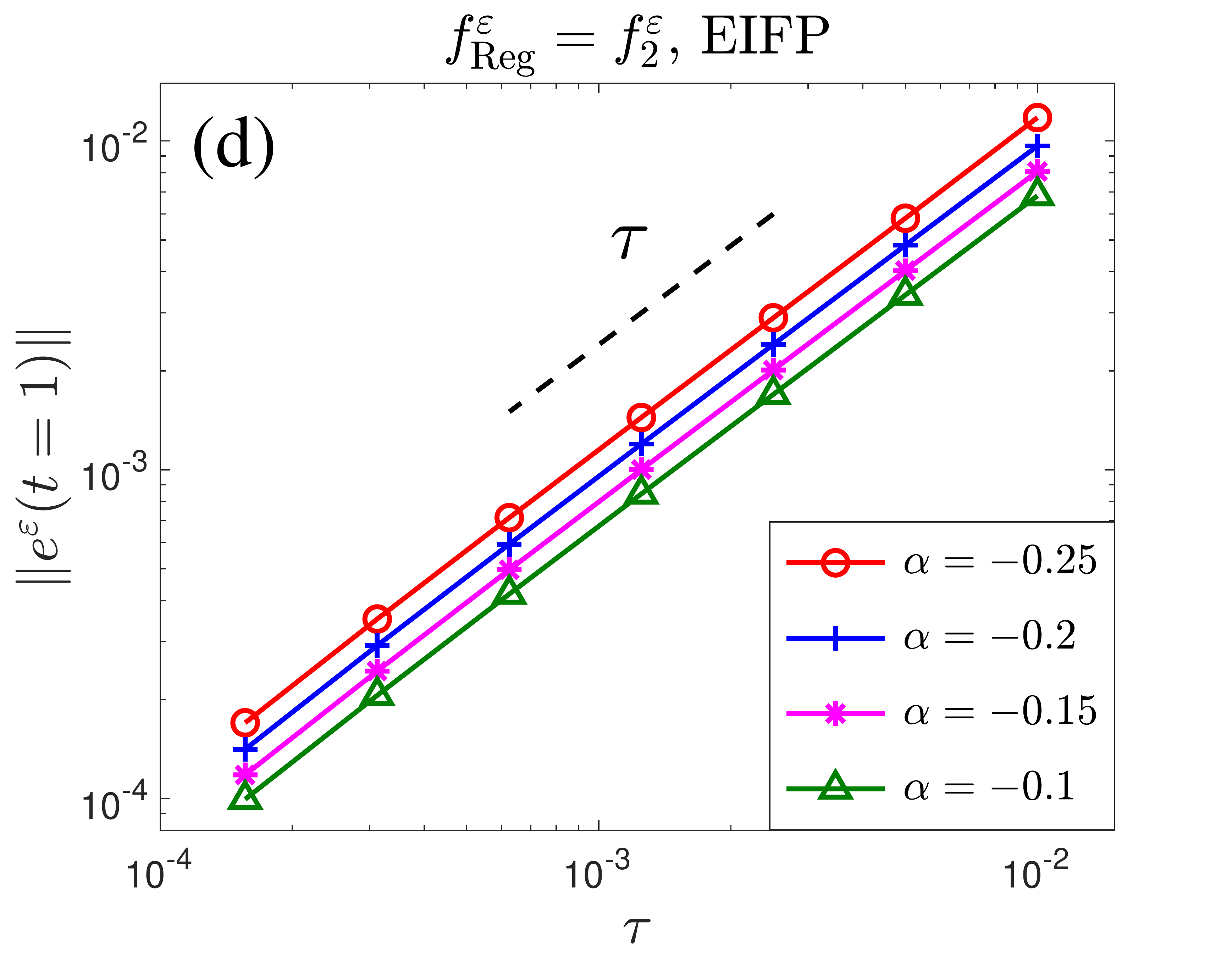}}
\end{minipage}
\begin{minipage}{0.5\textwidth}
\centerline{\includegraphics[width=7.5cm,height=5cm]{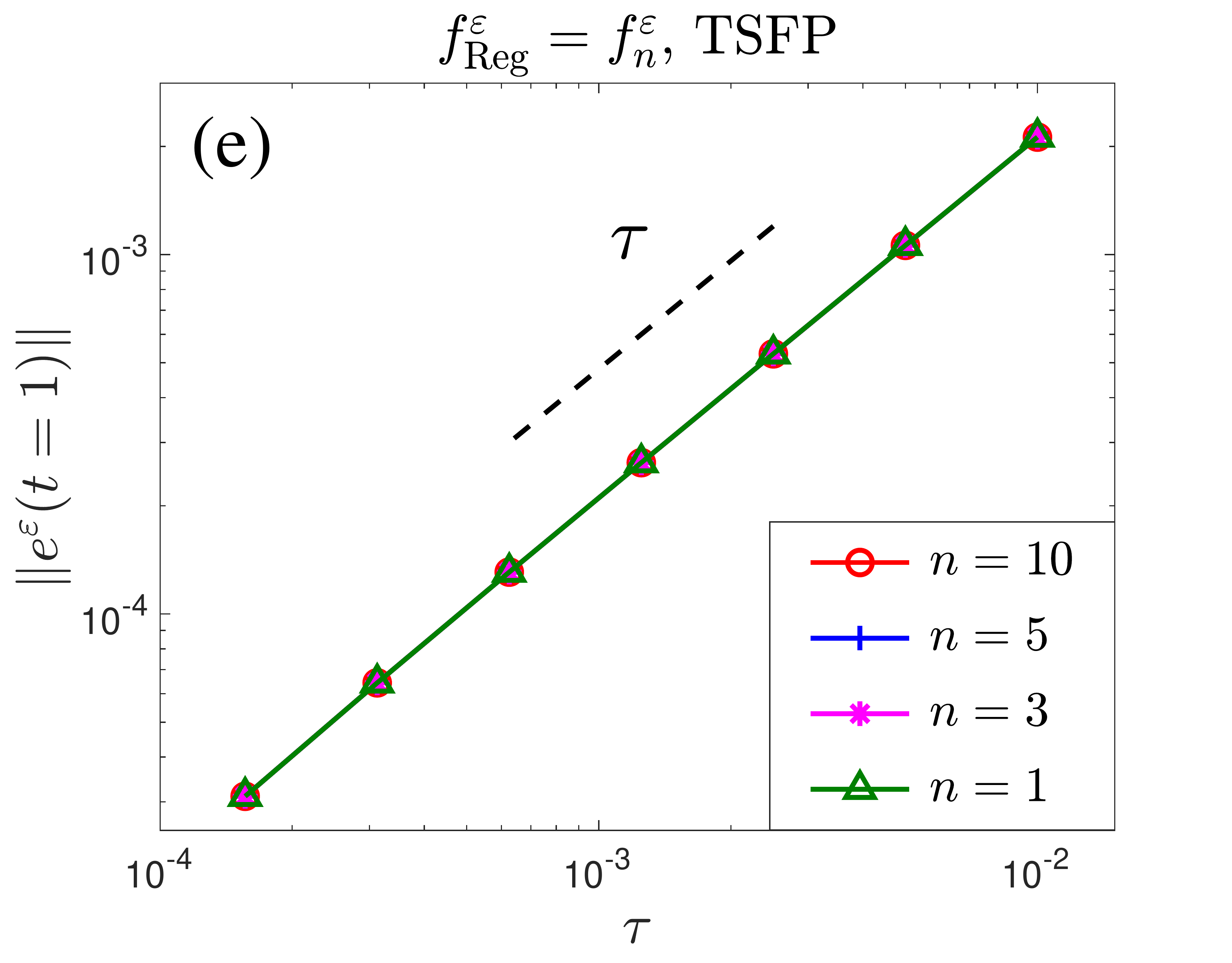}}
\end{minipage}
\begin{minipage}{0.5\textwidth}
\centerline{\includegraphics[width=7.5cm,height=5cm]{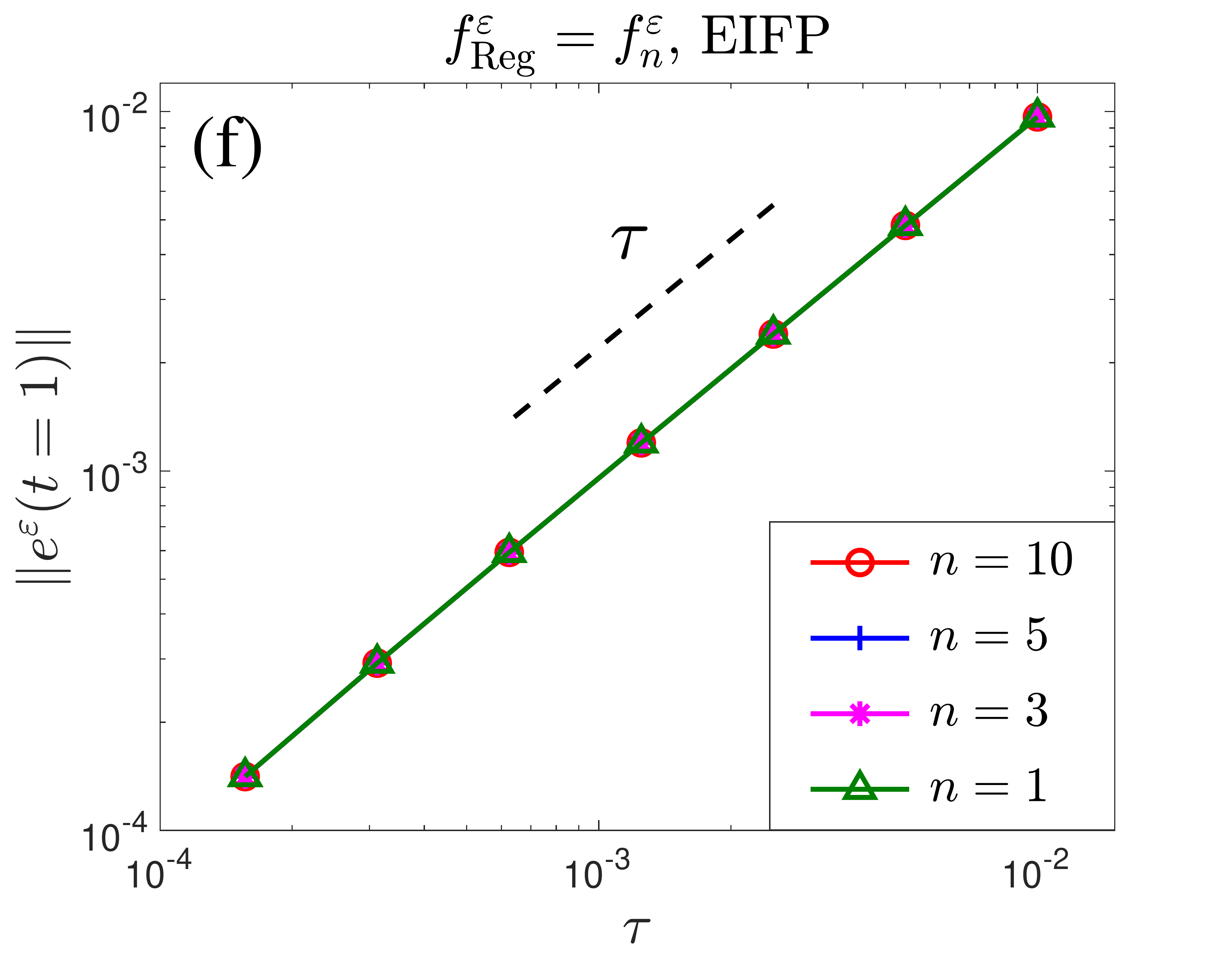}}
\end{minipage}
\caption{Convergence of the TSFP method \eqref{eq:TSFP} (left column) and EIFP method \eqref{eq:EIFP} (right column) to the erNLSE \eqref{eq:ERNLSE} for different $\eps$ (top row), $\alpha$ (middle row) and degree $n$ (bottom row), i.e., errors $\left\|e^{\eps}(t=1)\right\|$ versus the time step $\tau$.}
\label{fig:ef}
\end{figure}

\begin{figure}[h!]
\begin{minipage}{0.5\textwidth}
\centerline{\includegraphics[width=7.5cm,height=5cm]{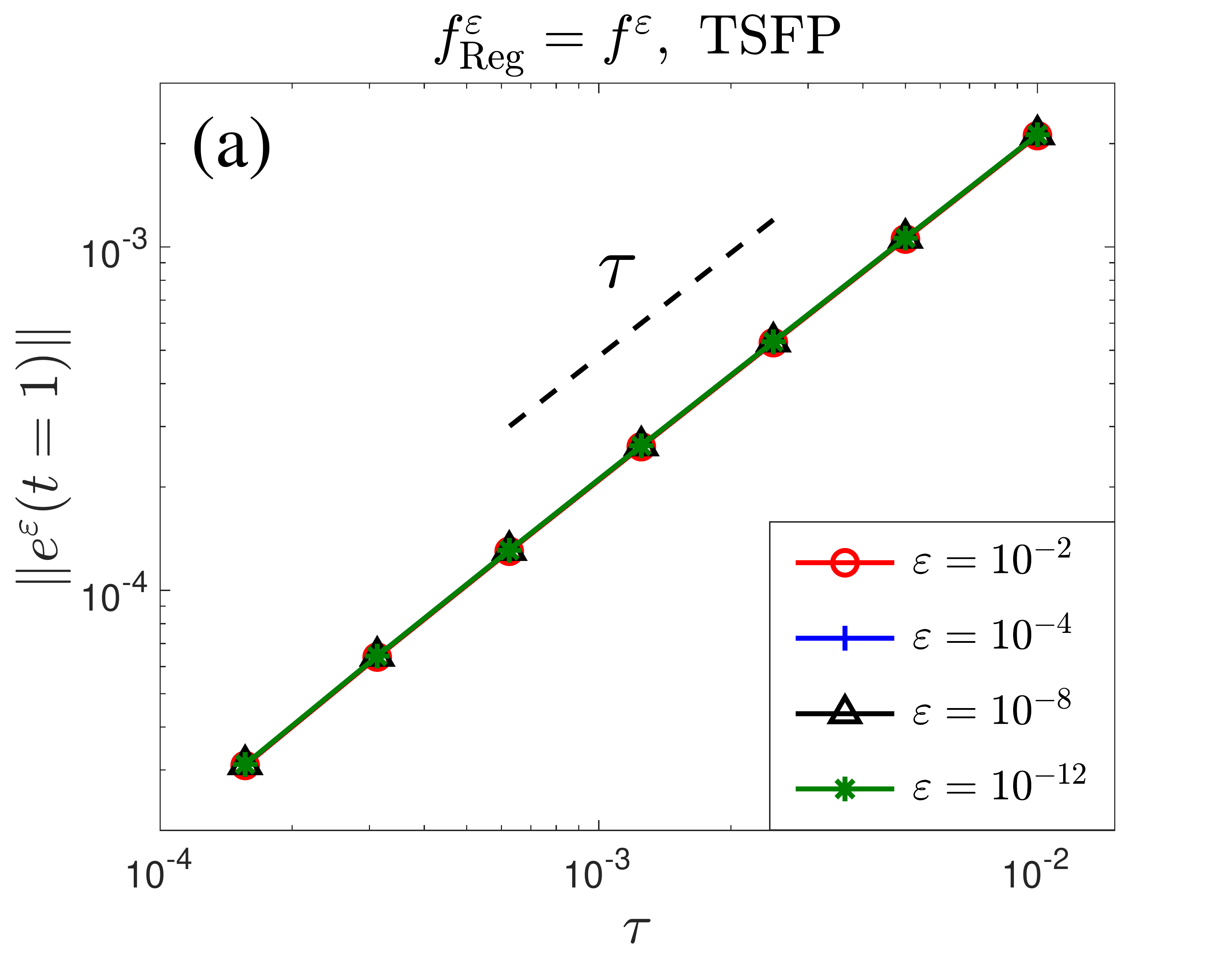}}
\end{minipage}
\begin{minipage}{0.5\textwidth}
\centerline{\includegraphics[width=7.5cm,height=5cm]{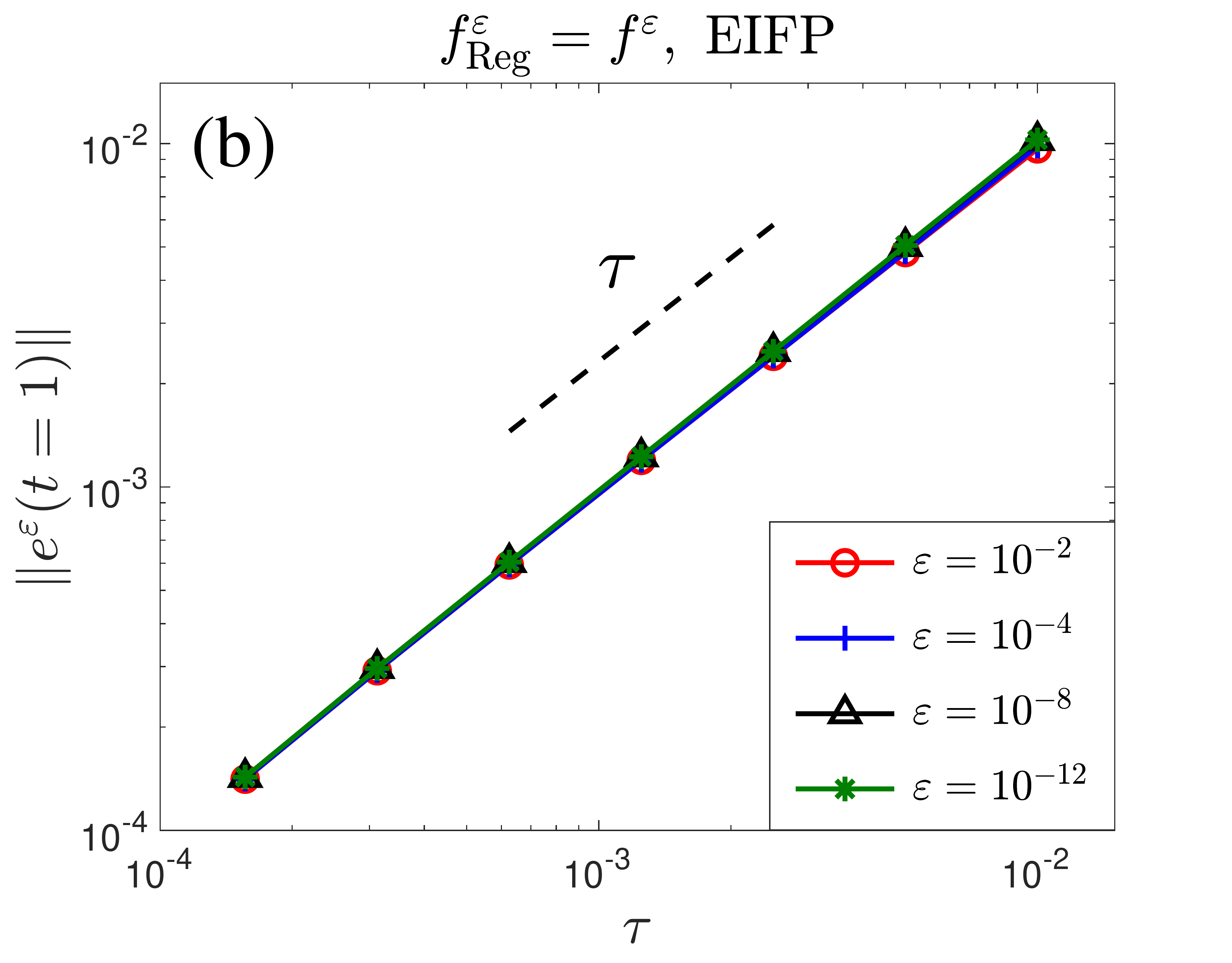}}
\end{minipage}
\begin{minipage}{0.5\textwidth}
\centerline{\includegraphics[width=7.5cm,height=5cm]{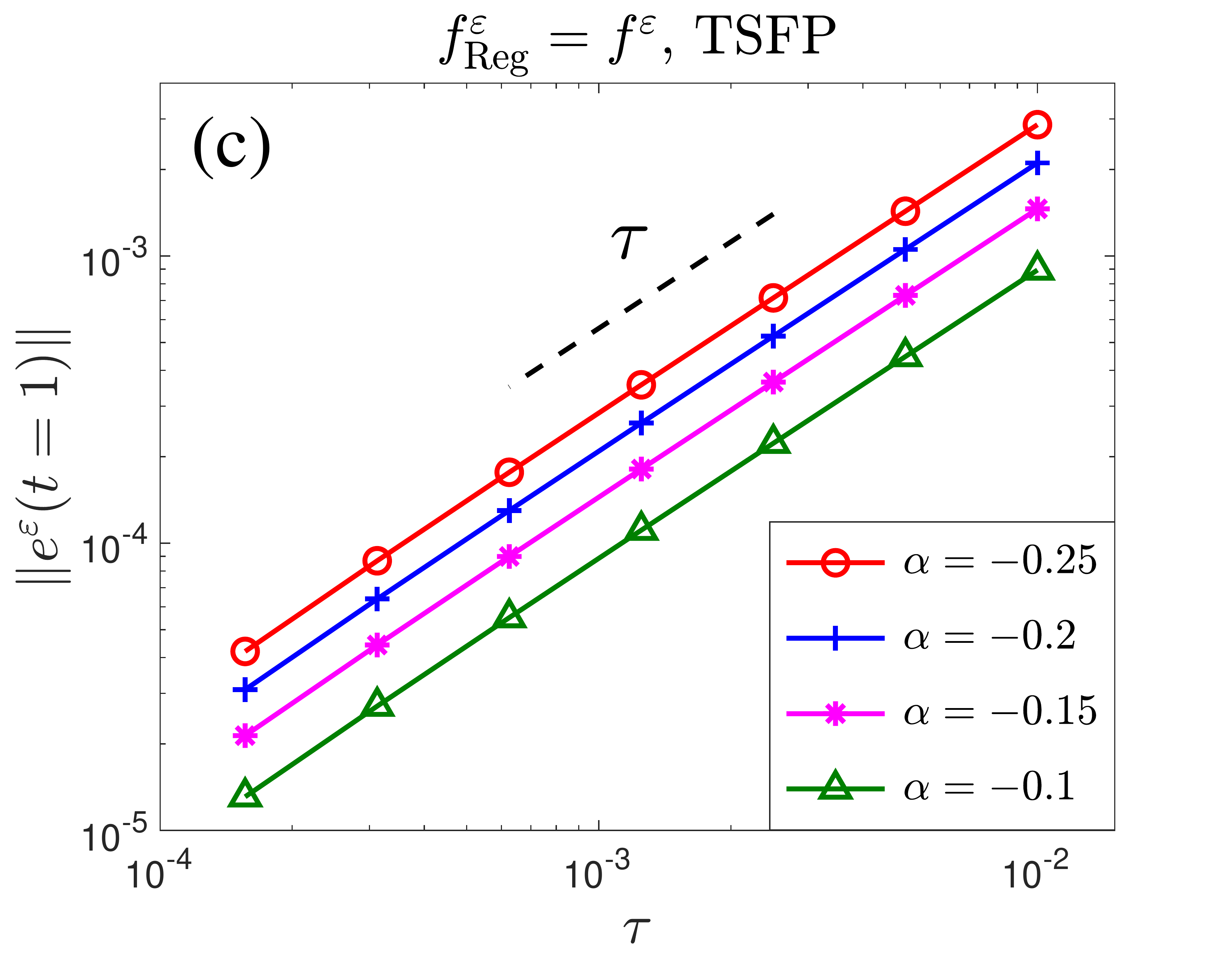}}
\end{minipage}
\begin{minipage}{0.5\textwidth}
\centerline{\includegraphics[width=7.5cm,height=5cm]{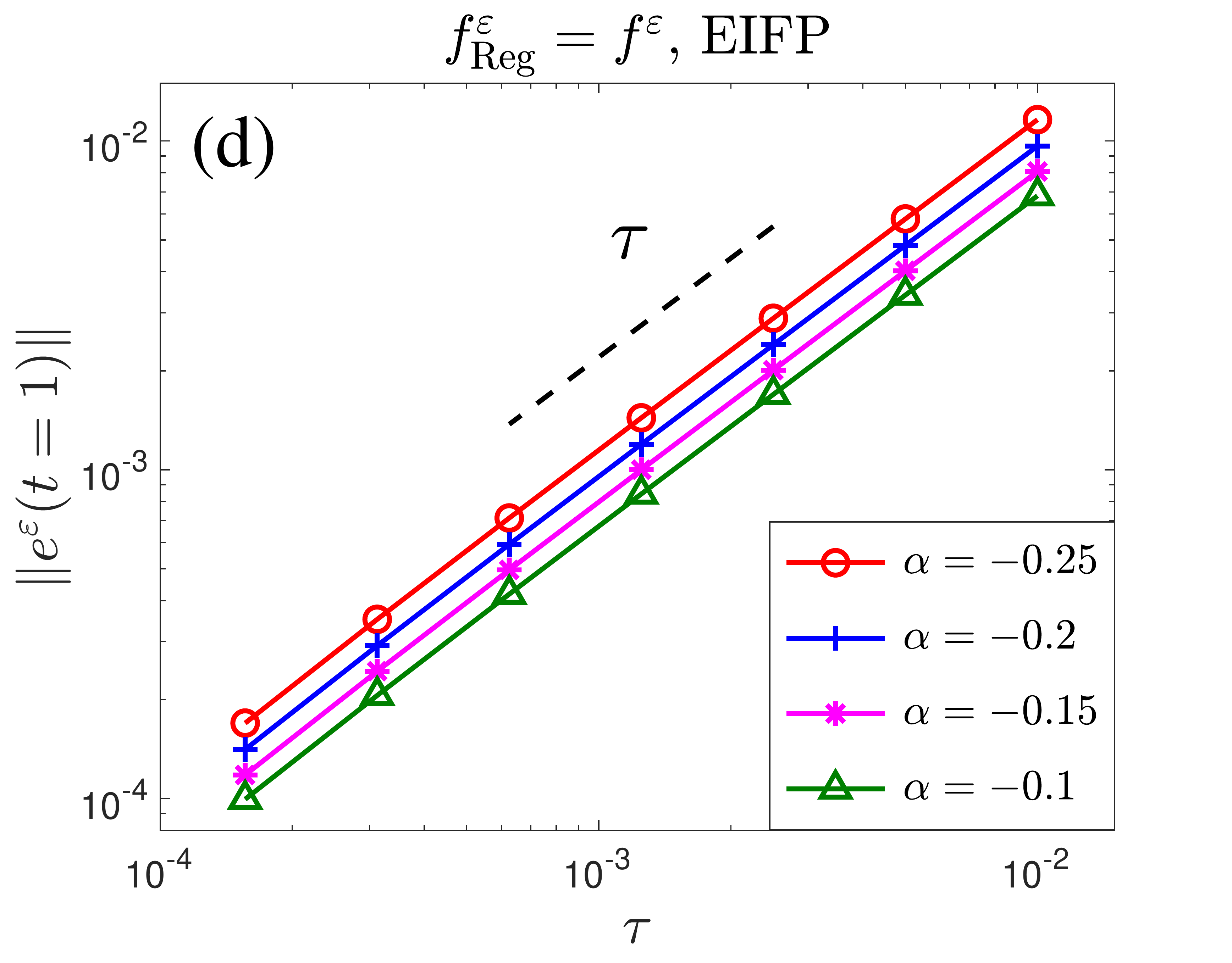}}
\end{minipage}
\caption{Convergence of the TSFP method \eqref{eq:TSFP} (left column) and EIFP method \eqref{eq:EIFP} (right column) to the rNLSE \eqref{eq:RNLSE1} for different $\eps$ (upper row) and $\alpha$ (lower row), i.e., errors $\left\|e^{\eps}(t=1)\right\|$ versus the time step $\tau$.}
\label{fig:f1}
\end{figure}

\subsection{Convergence rates of TSFP and EIFP methods}
\begin{figure}[ht!]
\begin{minipage}{0.5\textwidth}
\centerline{\includegraphics[width=7.5cm,height=5cm]{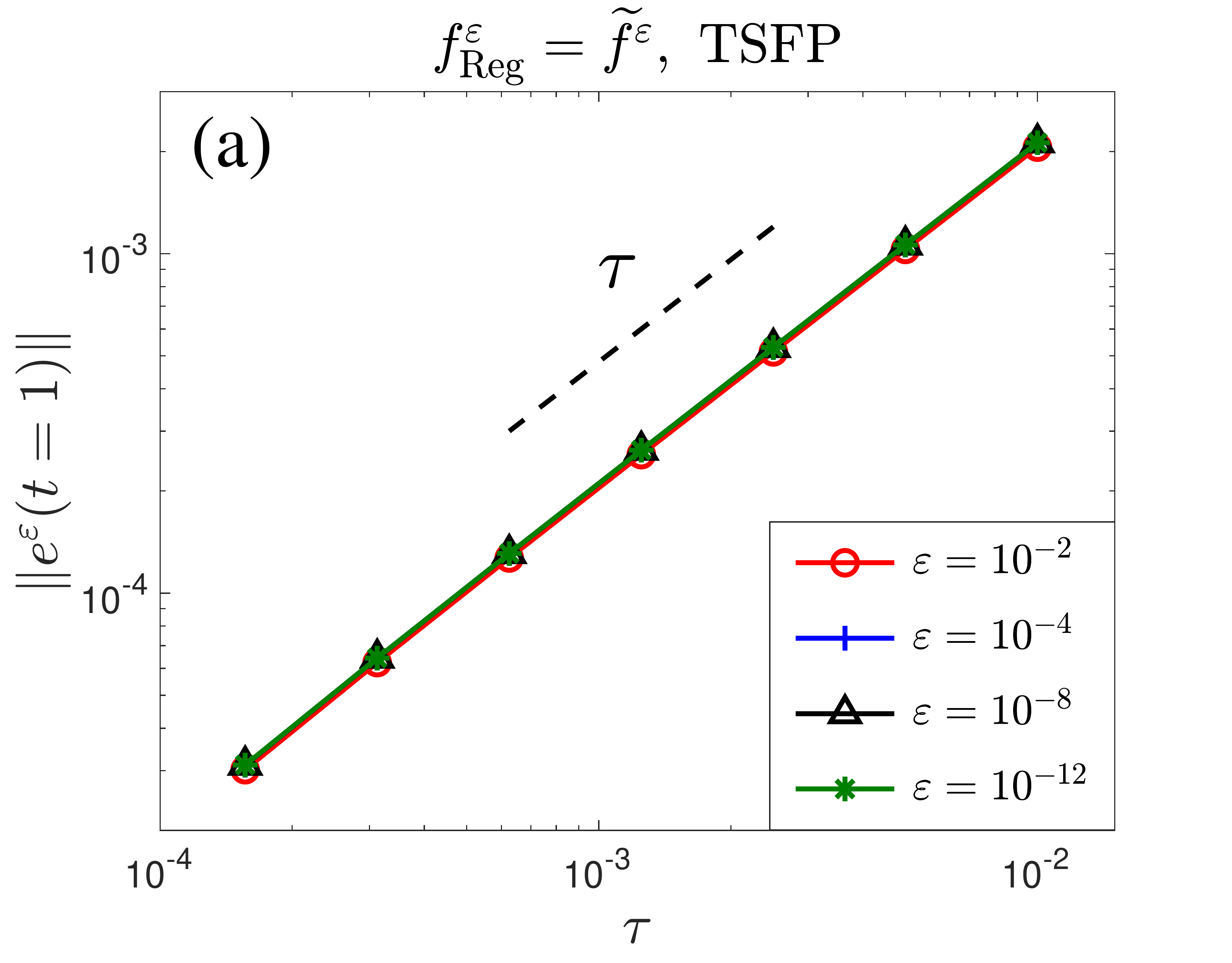}}
\end{minipage}
\begin{minipage}{0.5\textwidth}
\centerline{\includegraphics[width=7.5cm,height=5cm]{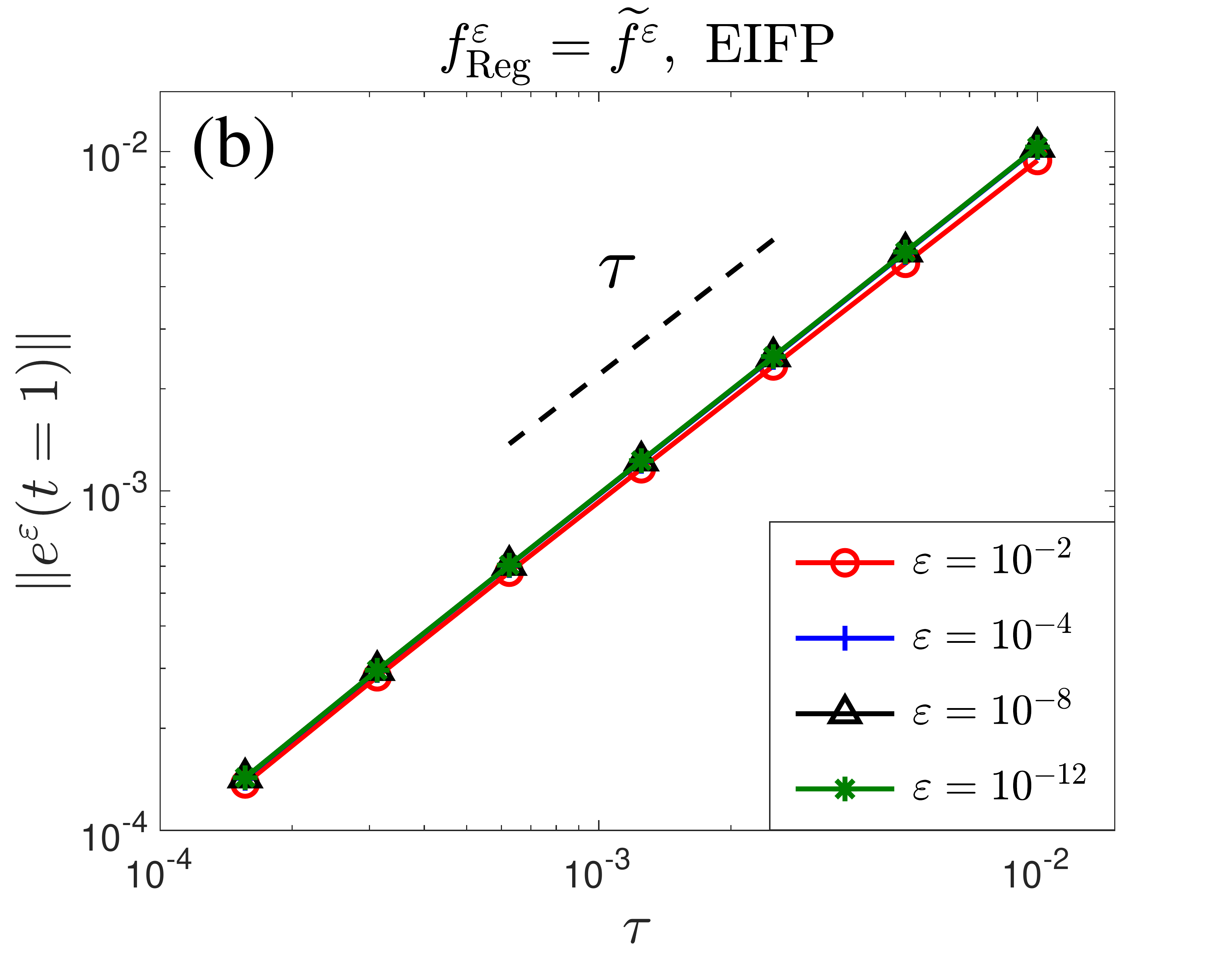}}
\end{minipage}
\begin{minipage}{0.5\textwidth}
\centerline{\includegraphics[width=7.5cm,height=5cm]{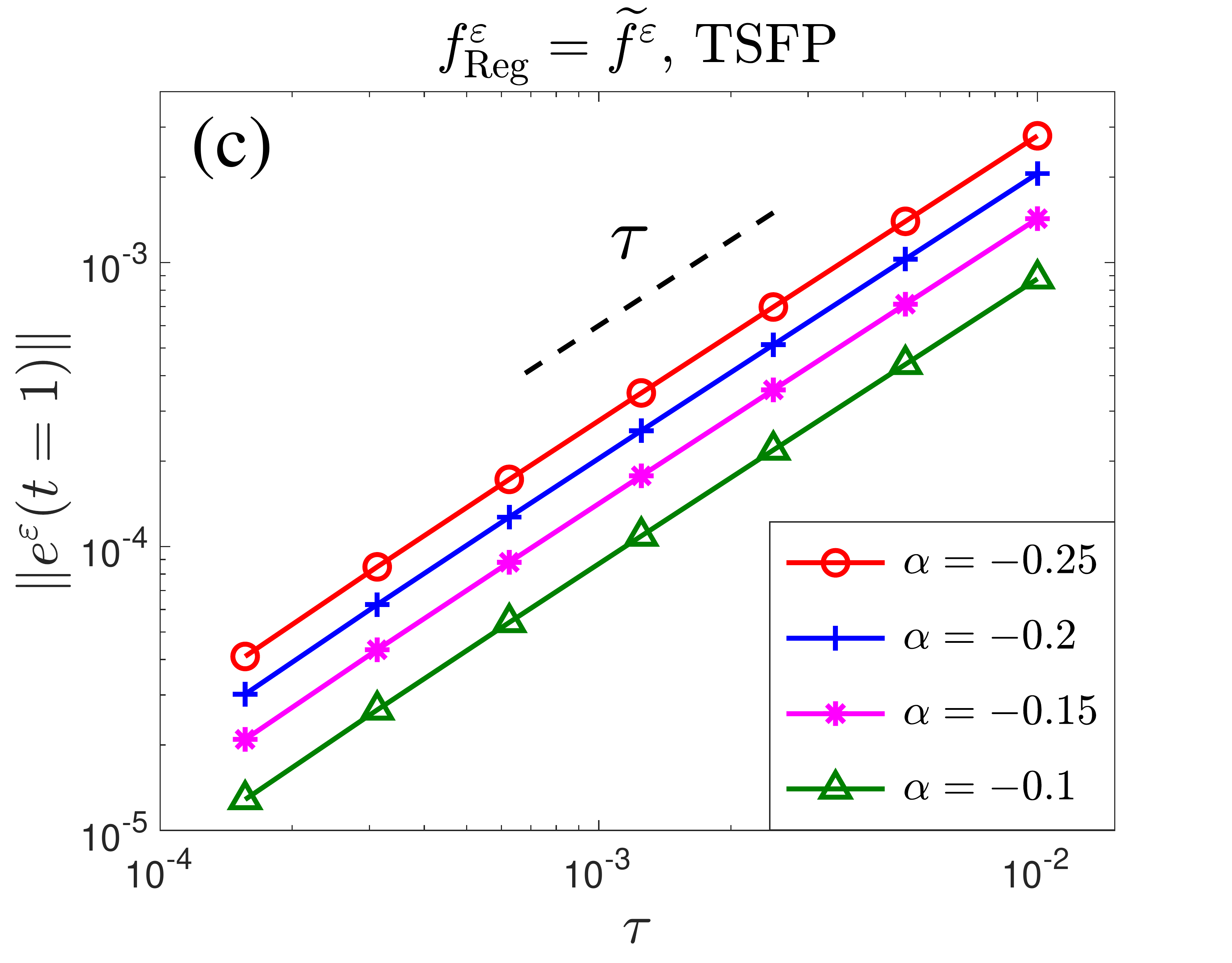}}
\end{minipage}
\begin{minipage}{0.5\textwidth}
\centerline{\includegraphics[width=7.5cm,height=5cm]{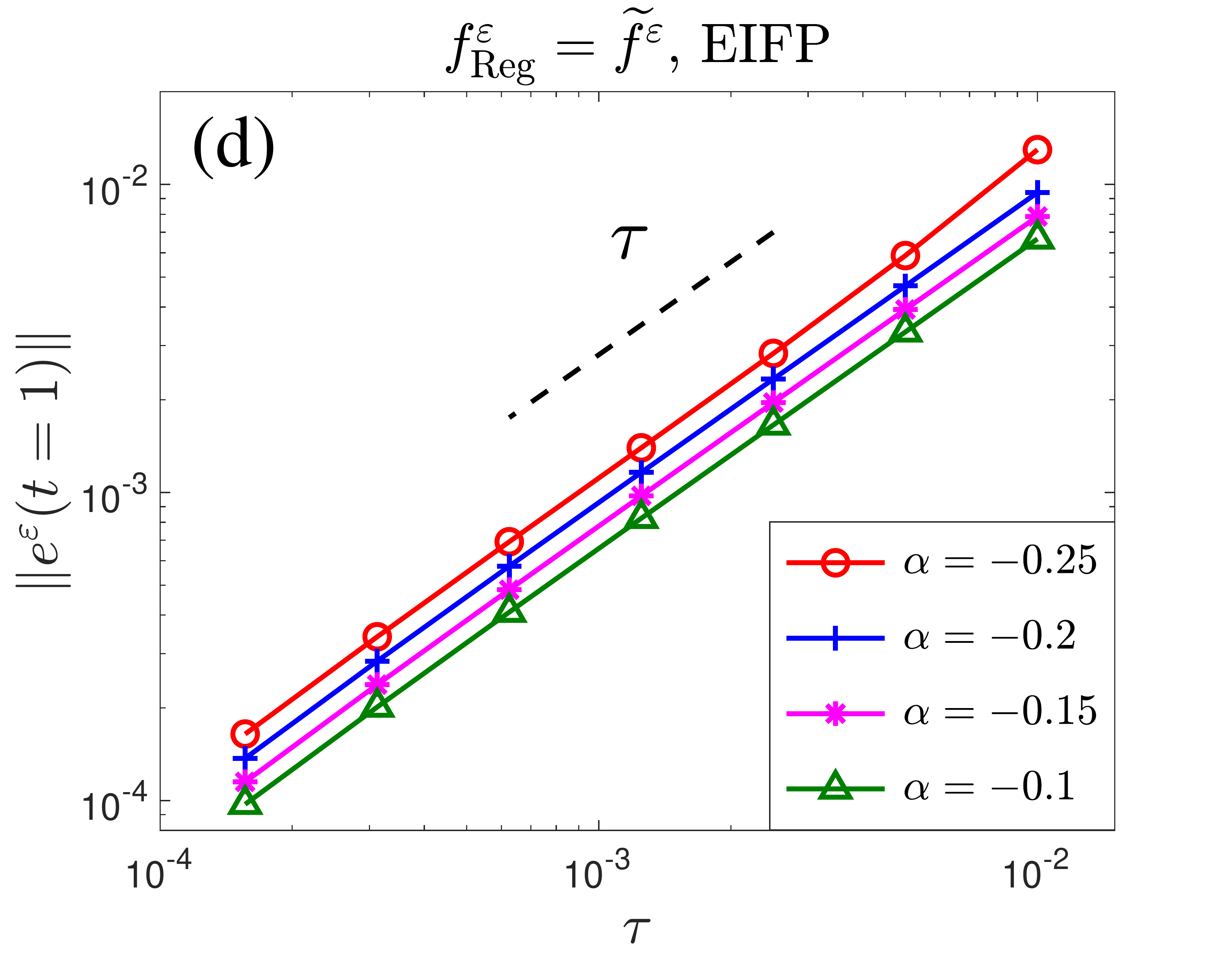}}
\end{minipage}
\caption{Convergence of the TSFP method \eqref{eq:TSFP} (left column) and EIFP method \eqref{eq:EIFP} (right column) to the rNLSE \eqref{eq:RNLSE2} for different $\eps$ (upper row) and $\alpha$ (lower row), i.e., errors $\left\|e^{\eps}(t=1)\right\|$ versus the time step $\tau$.}
\label{fig:f2}
\end{figure}

\begin{figure}[ht!]
\begin{minipage}{0.5\textwidth}
\centerline{\includegraphics[width=7.5cm,height=5cm]{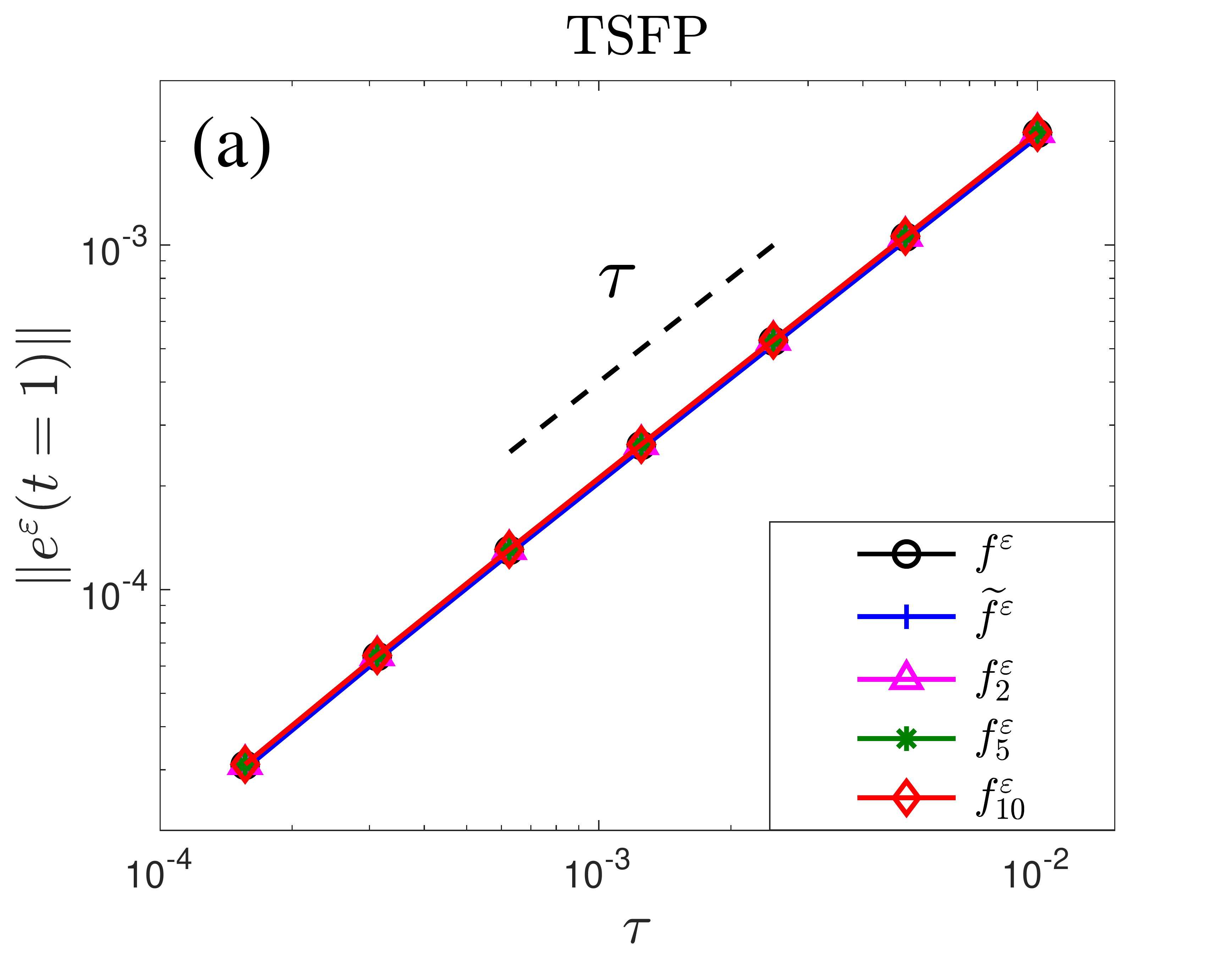}}
\end{minipage}
\begin{minipage}{0.5\textwidth}
\centerline{\includegraphics[width=7.5cm,height=5cm]{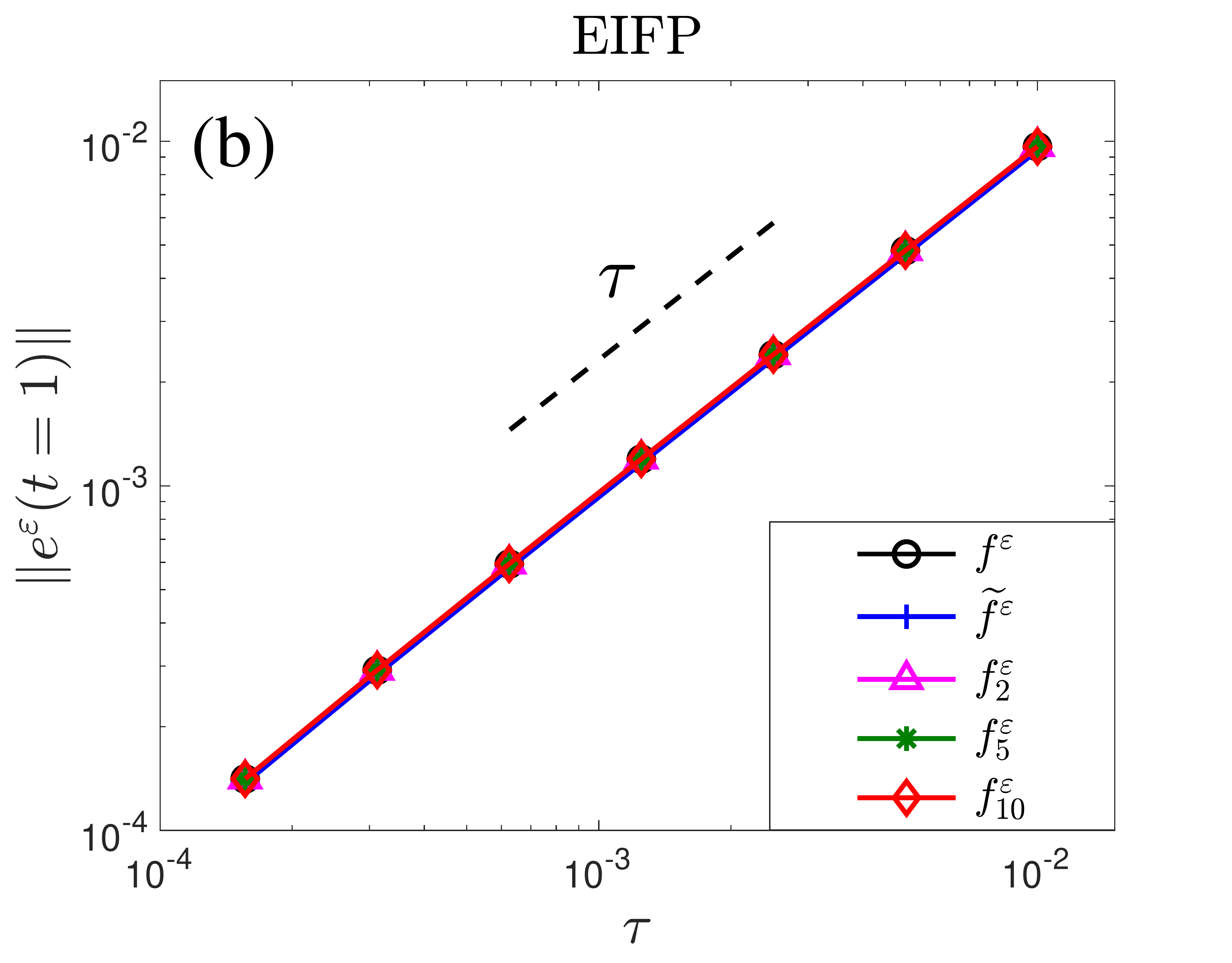}}
\end{minipage}
\caption{Convergence of the TSFP method \eqref{eq:TSFP} (left) and EIFP method \eqref{eq:EIFP} (right) for different regularizations, i.e., errors $\left\|e^{\eps}(t=1)\right\|$ versus the time step $\tau$.}
\label{fig:df}
\end{figure}

We are going to test the convergence rates of the TSFP method \eqref{eq:TSFP} and the EIFP method \eqref{eq:EIFP} to numerically solve the rNLSE \eqref{eq:GNLSE} in terms of the time step $\tau$ for fixed $0<\eps\ll1$. Figures \ref{fig:ef}, \ref{fig:f1} and \ref{fig:f2} display the convergence of the TSFP method \eqref{eq:TSFP} and the EIFP method \eqref{eq:EIFP} to solve these three regularized models, respectively, for different $\eps$ and $\alpha$ with different degree $n$ for the local energy regularization. Figure \ref{fig:df} compares the errors of the TSFP method \eqref{eq:TSFP} and the EIFP method \eqref{eq:EIFP} for different regularizations.

 From these figures and additional numerical results not shown here for brevity, we can clearly see that: (i) For all the regularizations with the fixed $0<\eps\ll1$  and $\alpha \in (-1/3, 0)$, the TSFP method \eqref{eq:TSFP} and the EIFP method \eqref{eq:EIFP} converge linearly in time. (ii) The errors of the TSFP method are much smaller than those of the EIFP method for the same $\eps$ and $\alpha$ by choosing the same time step $\tau$. (iii) For different regularizations, the errors of the each numerical scheme are close. The errors of the TSFP method and EIFP method are independent of $\eps$ and the degree $n$ for the local energy regularization, but depend on the the parameter $\alpha$. When $\alpha$ is close to $0^-$, the errors of these two numerical schemes become smaller and the change of the TSFP method is larger than that of the EIFP method. We remark here, the independence of the errors of the TSFP method and EIFP method in terms of $\eps$ is just from the view of numerical simulations, while it deserves to show analytically and it is our ongoing work.

\subsection{Convergence rates of regularized numerical methods}
By the definition of the error functions, we explain the errors in the following diagram
{\Large\begin{equation*}
  \begin{tikzcd}[row sep=large, column sep=large]
      \psi^{\eps, m} \arrow{r}{e^{\eps}(t_m)} \arrow[dashed,swap]{dr}{\widetilde{e}^{\eps}(t_m)} & [7em] \psi^{\eps}(\cdot, t_m)  \arrow{d}{\widehat{e}^{\eps}(t_m)} \\
     & \psi(\cdot, t_m)
  \end{tikzcd}
\end{equation*}}

\begin{figure}[ht!]
\begin{minipage}{0.5\textwidth}
\centerline{\includegraphics[width=7.5cm,height=5cm]{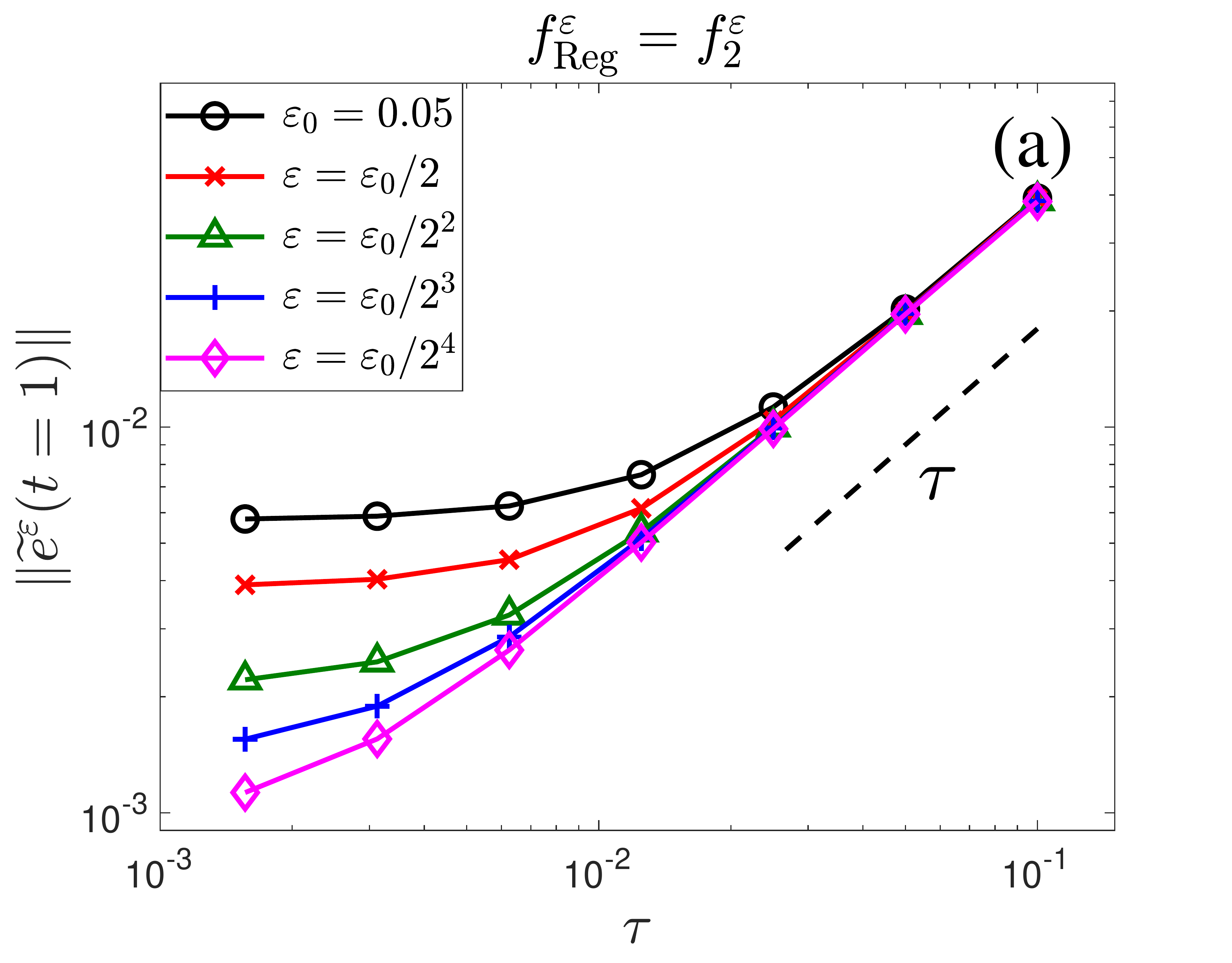}}
\end{minipage}
\begin{minipage}{0.5\textwidth}
\centerline{\includegraphics[width=7.5cm,height=5cm]{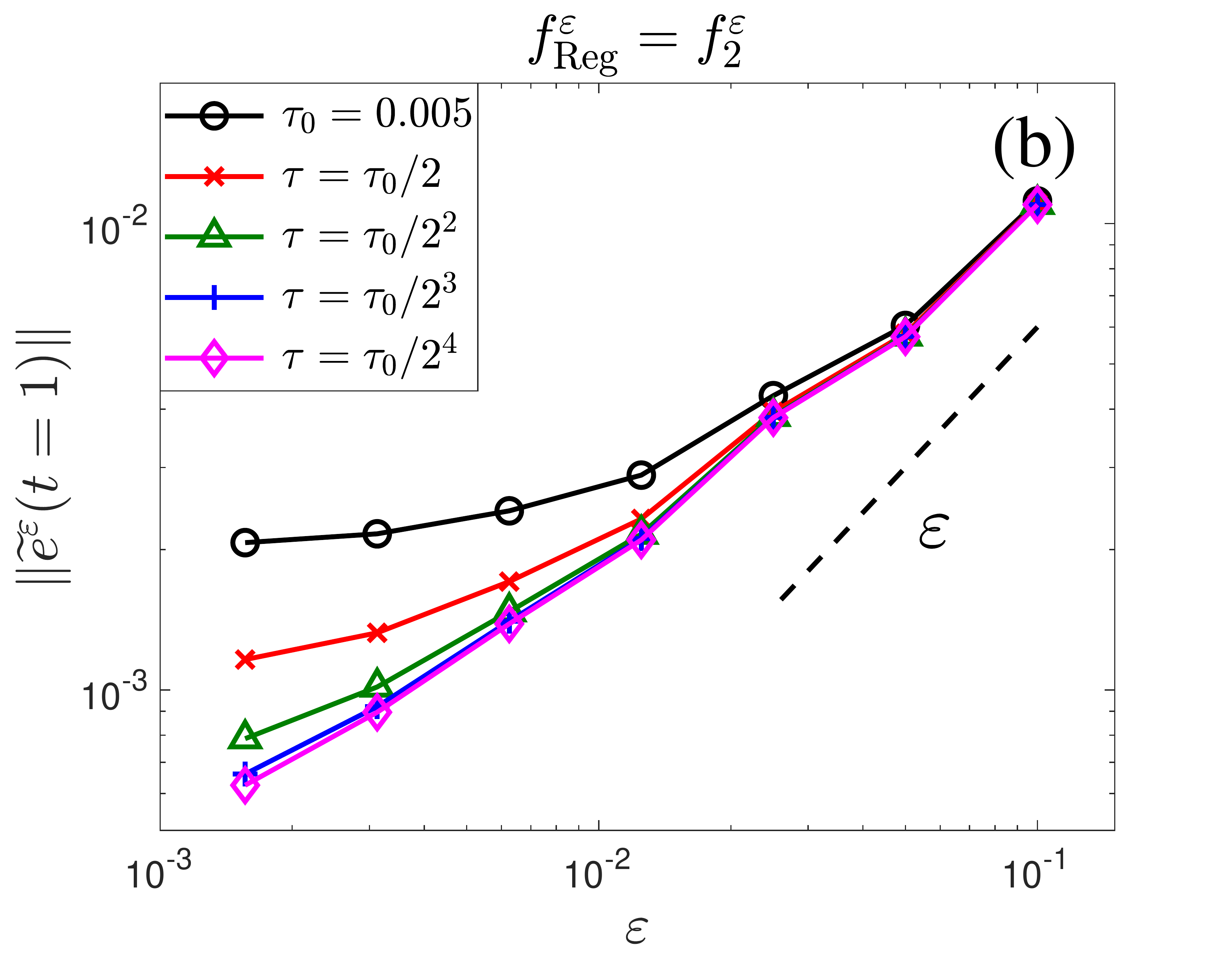}}
\end{minipage}
\caption{Convergence of the TSFP method \eqref{eq:TSFP} for the erNLSE \eqref{eq:ERNLSE} to the sNLSE \eqref{eq:NLSE1} for (a) different $\eps$; and  (b) different $\tau$.}
\label{fig:tee}
\end{figure}

\begin{figure}[ht!]
\begin{minipage}{0.5\textwidth}
\centerline{\includegraphics[width=7.5cm,height=5cm]{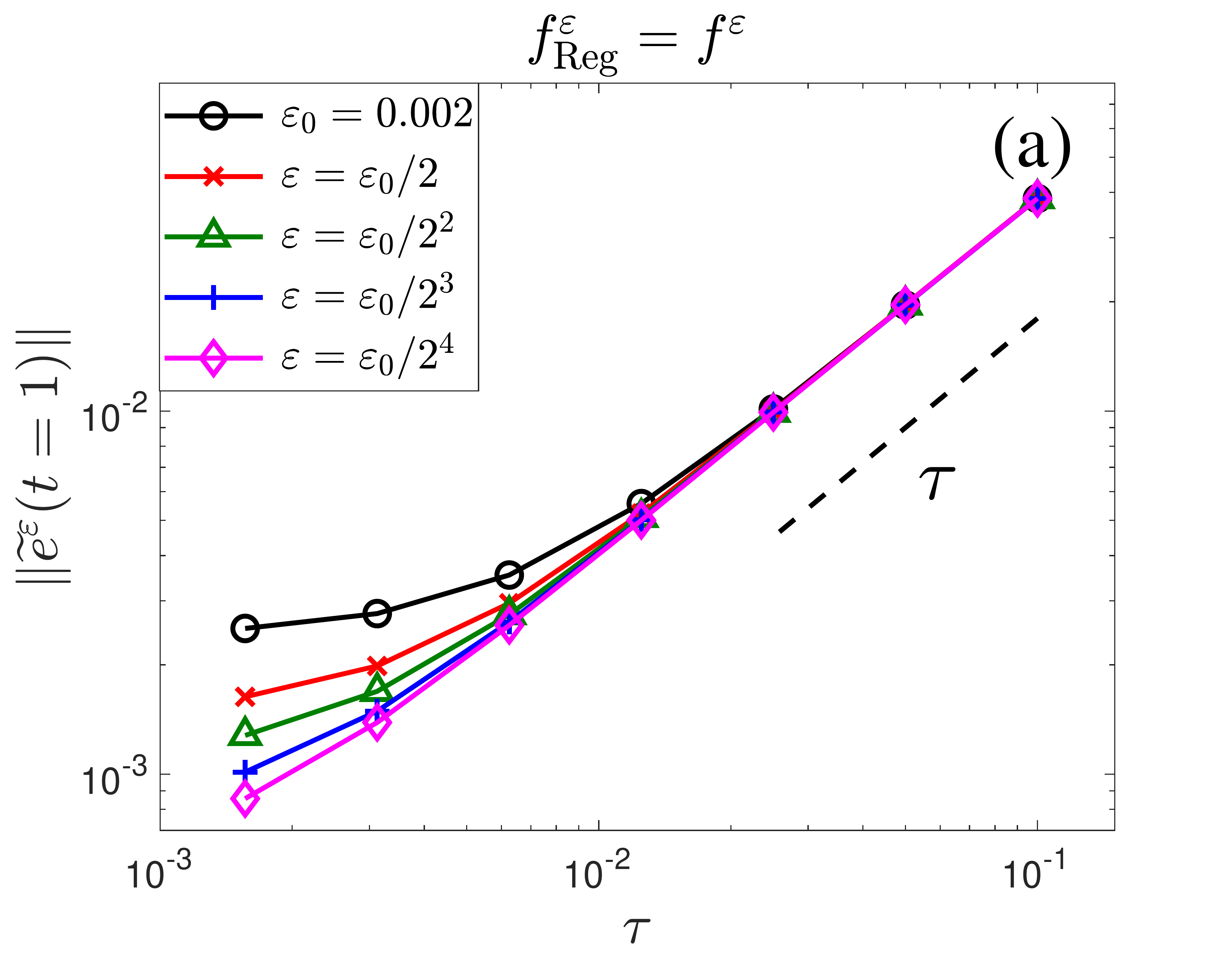}}
\end{minipage}
\begin{minipage}{0.5\textwidth}
\centerline{\includegraphics[width=7.5cm,height=5cm]{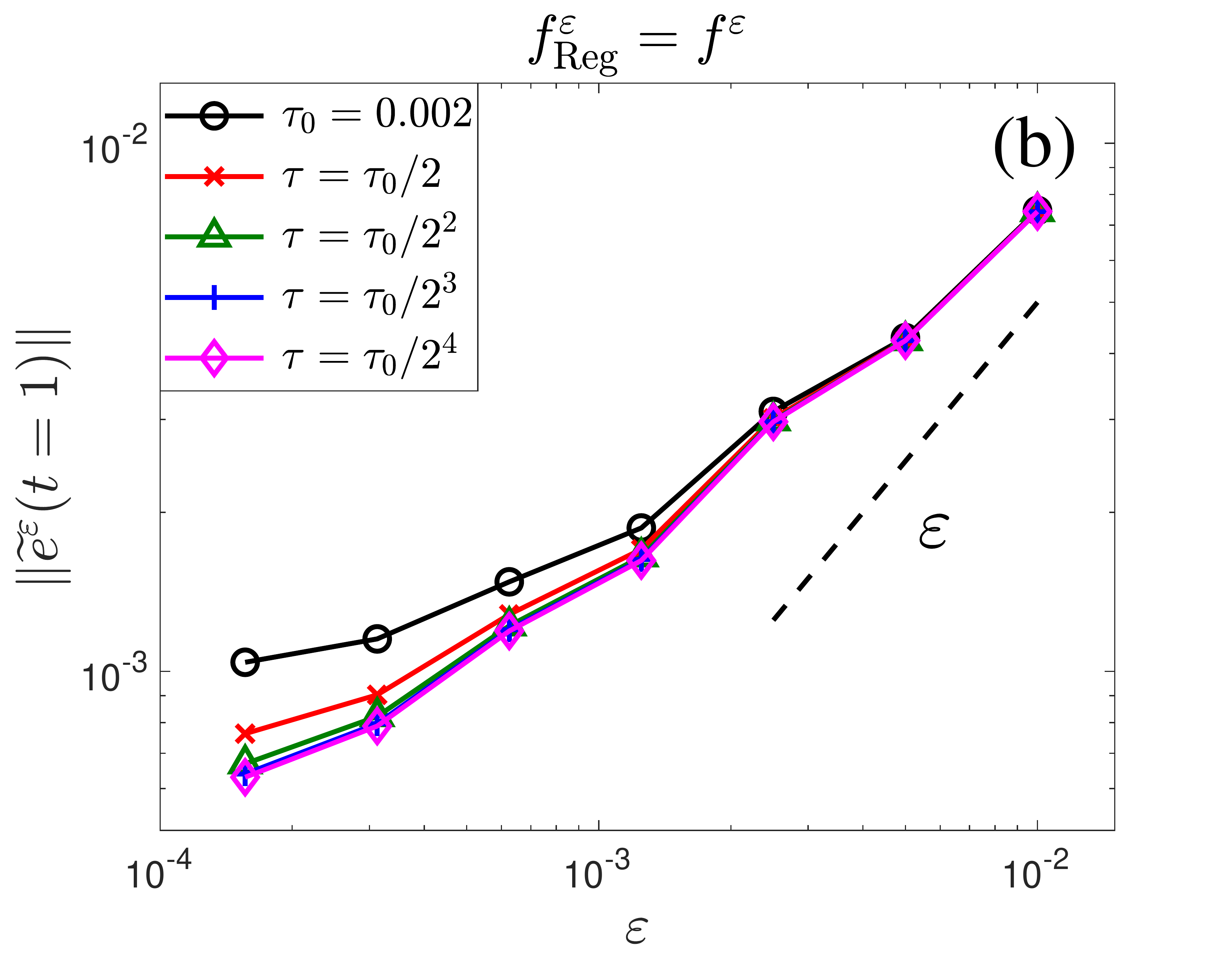}}
\end{minipage}
\caption{Convergence of the TSFP method \eqref{eq:TSFP} for the rNLSE \eqref{eq:RNLSE1} to the sNLSE \eqref{eq:NLSE1} for (a) different $\eps$; and (b) different $\tau$.}
\label{fig:te1}
\end{figure}

\begin{figure}[ht!]
\begin{minipage}{0.5\textwidth}
\centerline{\includegraphics[width=7.5cm,height=5cm]{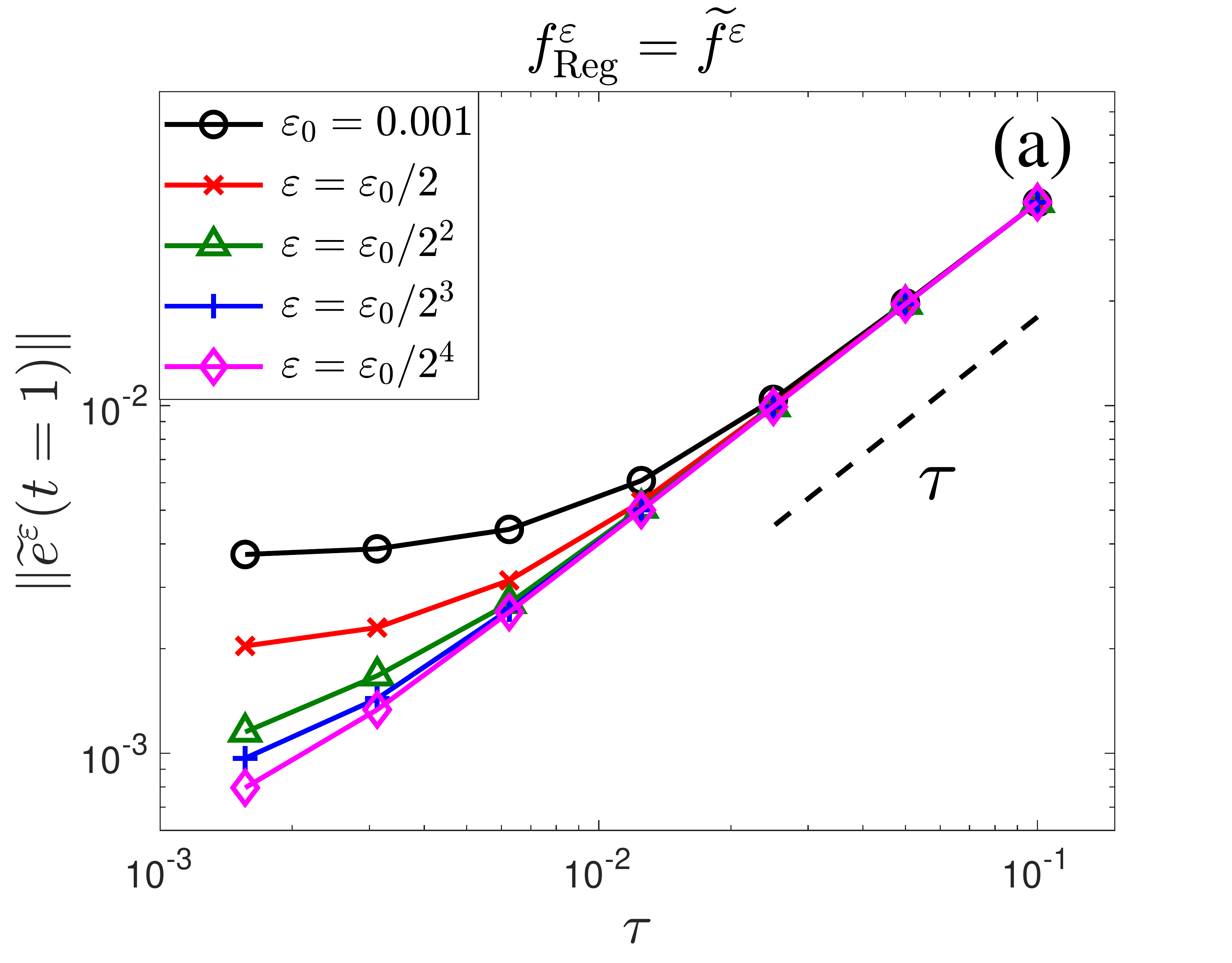}}
\end{minipage}
\begin{minipage}{0.5\textwidth}
\centerline{\includegraphics[width=7.5cm,height=5cm]{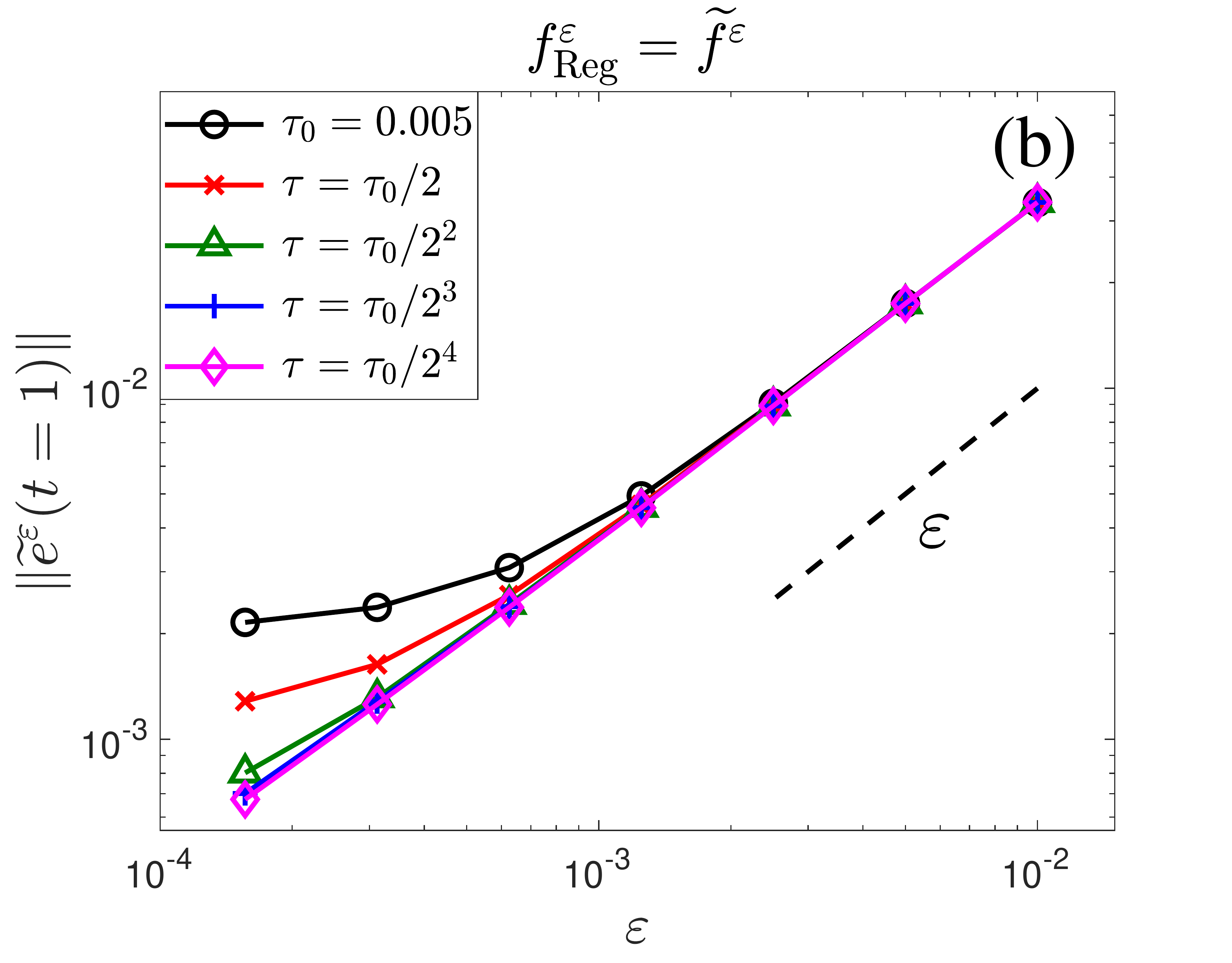}}
\end{minipage}
\caption{Convergence of the TSFP method \eqref{eq:TSFP} for the rNLSE \eqref{eq:RNLSE2} to the sNLSE \eqref{eq:NLSE1} for (a) different $\eps$; and (b) different $\tau$.}
\label{fig:te2}
\end{figure}

\noindent In order to get an accurate approximations for the solution  $\psi$ of sNLSE \eqref{eq:NLSE1}, we need to choose a suitable regularized model and an accurate numerical scheme to solve the regularized model. Then, we will test the converge rate of the error $\widetilde{e}^{\eps}(t_m)$. Figures \ref{fig:tee}, \ref{fig:te1} and \ref{fig:te2} depict the errors $\left\|\widetilde{e}^{\eps}(t=1)\right\|$ for the TSFP method \eqref{eq:TSFP} to solve these three regularized models.

 From these figures and additional numerical results not shown here for brevity, we can observe: (i) When the time step $\tau$ become smaller, the first-order convergence can be observed when $\eps$ is smaller. (ii) When the time step $\tau$ is small enough, the rNLSE \eqref{eq:RNLSE2} converges linearly and the convergence of the other two regularizations are slower.

\subsection{Application for the dynamics in 2D}
In this section, we compare the dynamics of a single Gausson, a vortex pair and a vortex dipole in two dimensions (2D) for different exponent $\alpha$ in the nonlinearity \eqref{nonlinear}, e.g., $\alpha = -0.3$, $\alpha = -0.1$ and  $\alpha = 1$ (the cubic NLSE). Here, the initial data in the three cases are given as \cite{BLS,KJZB,Zhang} :

Case I. A single Gausson, i.e.,
\begin{equation*}
\psi_0(x, y) = e^{-(x^2+y^2)}.
\end{equation*}

Case II. A vortex pair, i.e.,
\begin{equation*}
\psi_0(x, y) = ((x-0.5)+iy)((x+0.5)+iy)e^{-(x^2+y^2)}.
\end{equation*}

Case III. A vortex dipole, i.e.,
\begin{equation*}
\psi_0(x, y) = ((x-0.5)+iy)((x+0.5)-iy)e^{-(x^2+y^2)}.
\end{equation*}

\begin{figure}[ht!]
\begin{minipage}{0.33\textwidth}
\centerline{\includegraphics[width=4.2cm,height=3.5cm]{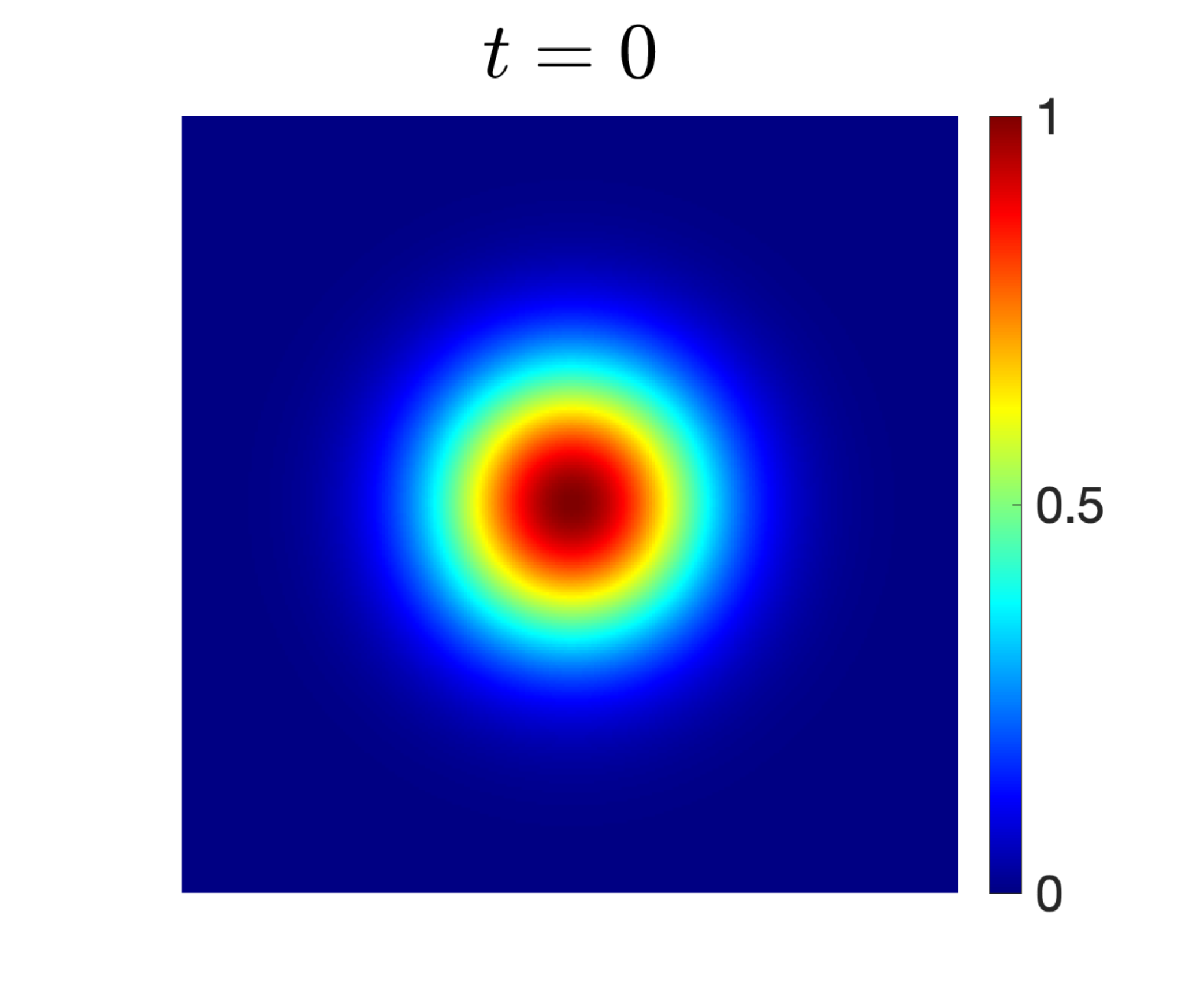}}
\end{minipage}
\begin{minipage}{0.33\textwidth}
\centerline{\includegraphics[width=4.2cm,height=3.5cm]{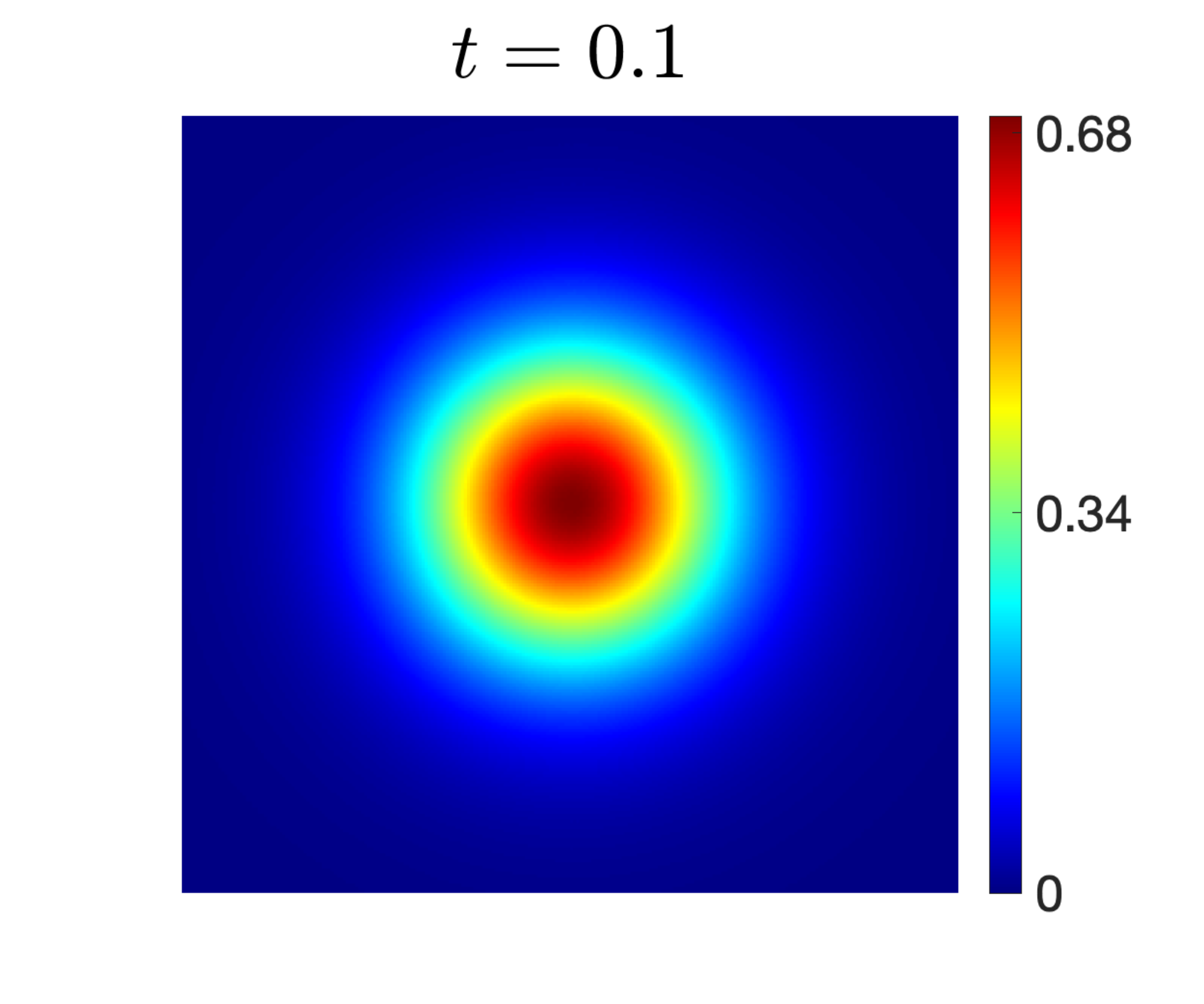}}
\end{minipage}
\begin{minipage}{0.33\textwidth}
\centerline{\includegraphics[width=4.2cm,height=3.5cm]{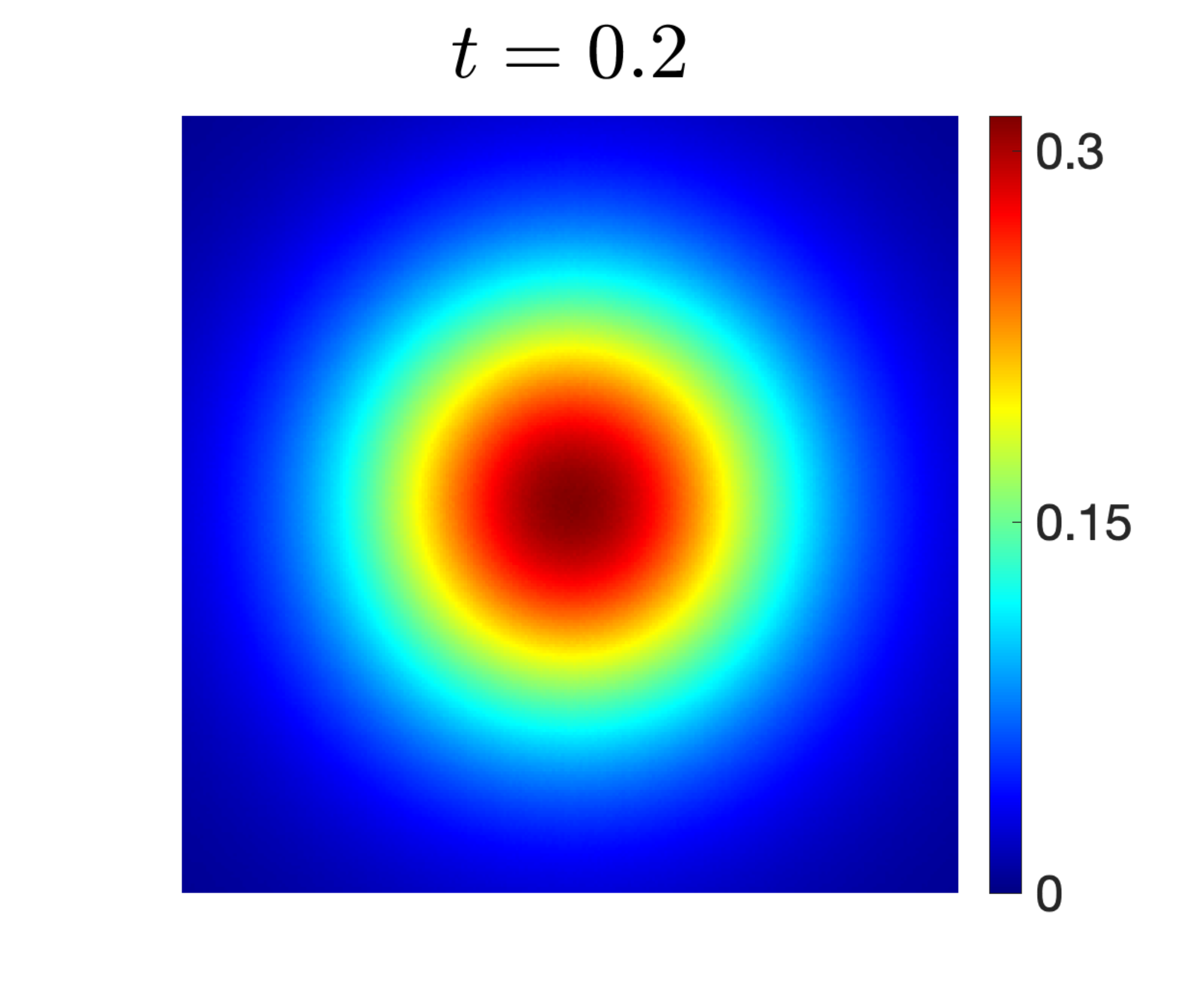}}
\end{minipage}
\begin{minipage}{0.33\textwidth}
\centerline{\includegraphics[width=4.2cm,height=3.5cm]{case1_in_den.eps}}
\end{minipage}
\begin{minipage}{0.33\textwidth}
\centerline{\includegraphics[width=4.2cm,height=3.5cm]{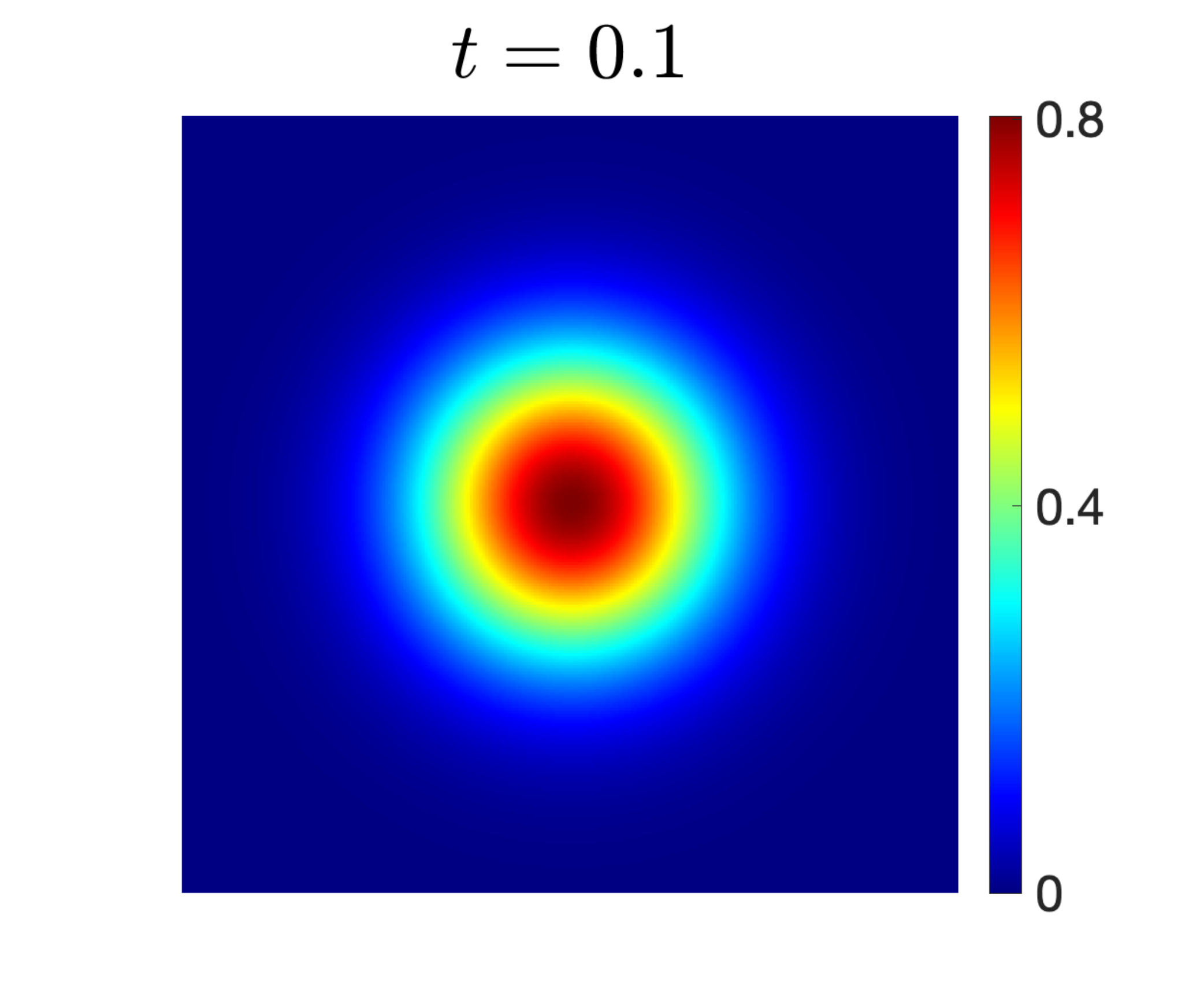}}
\end{minipage}
\begin{minipage}{0.33\textwidth}
\centerline{\includegraphics[width=4.2cm,height=3.5cm]{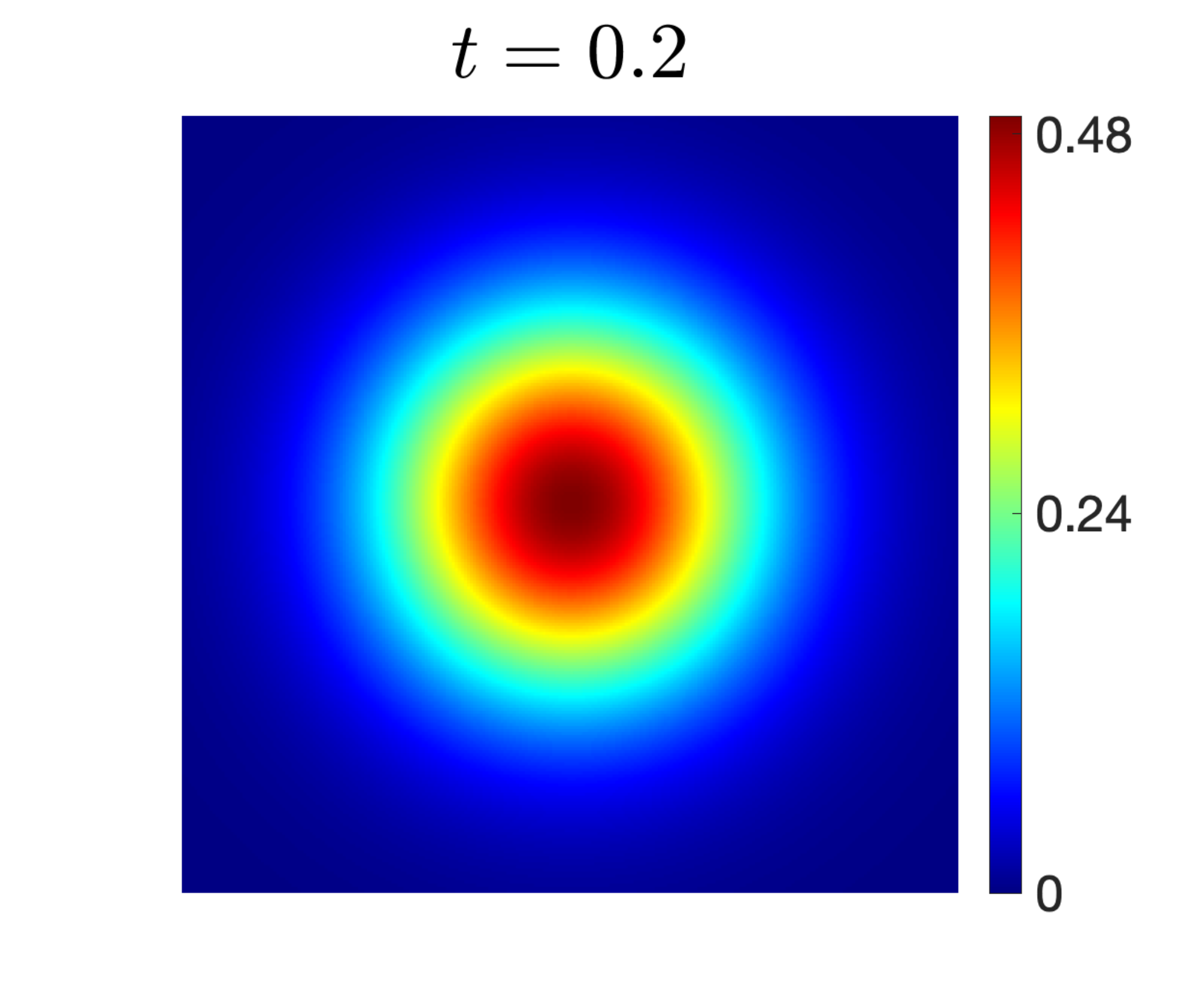}}
\end{minipage}
\begin{minipage}{0.33\textwidth}
\centerline{\includegraphics[width=4.2cm,height=3.5cm]{case1_in_den.eps}}
\end{minipage}
\begin{minipage}{0.33\textwidth}
\centerline{\includegraphics[width=4.2cm,height=3.5cm]{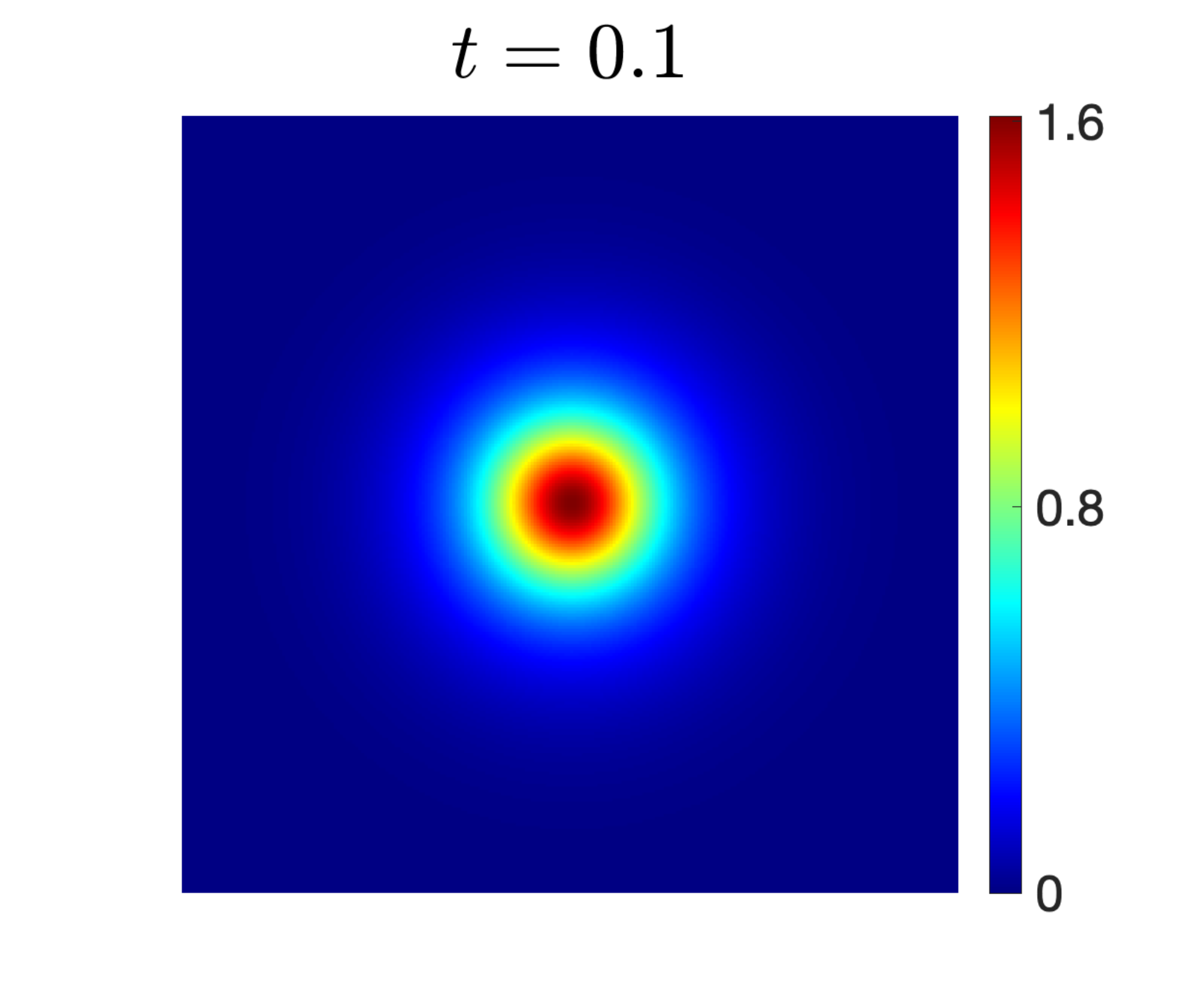}}
\end{minipage}
\begin{minipage}{0.33\textwidth}
\centerline{\includegraphics[width=4.2cm,height=3.5cm]{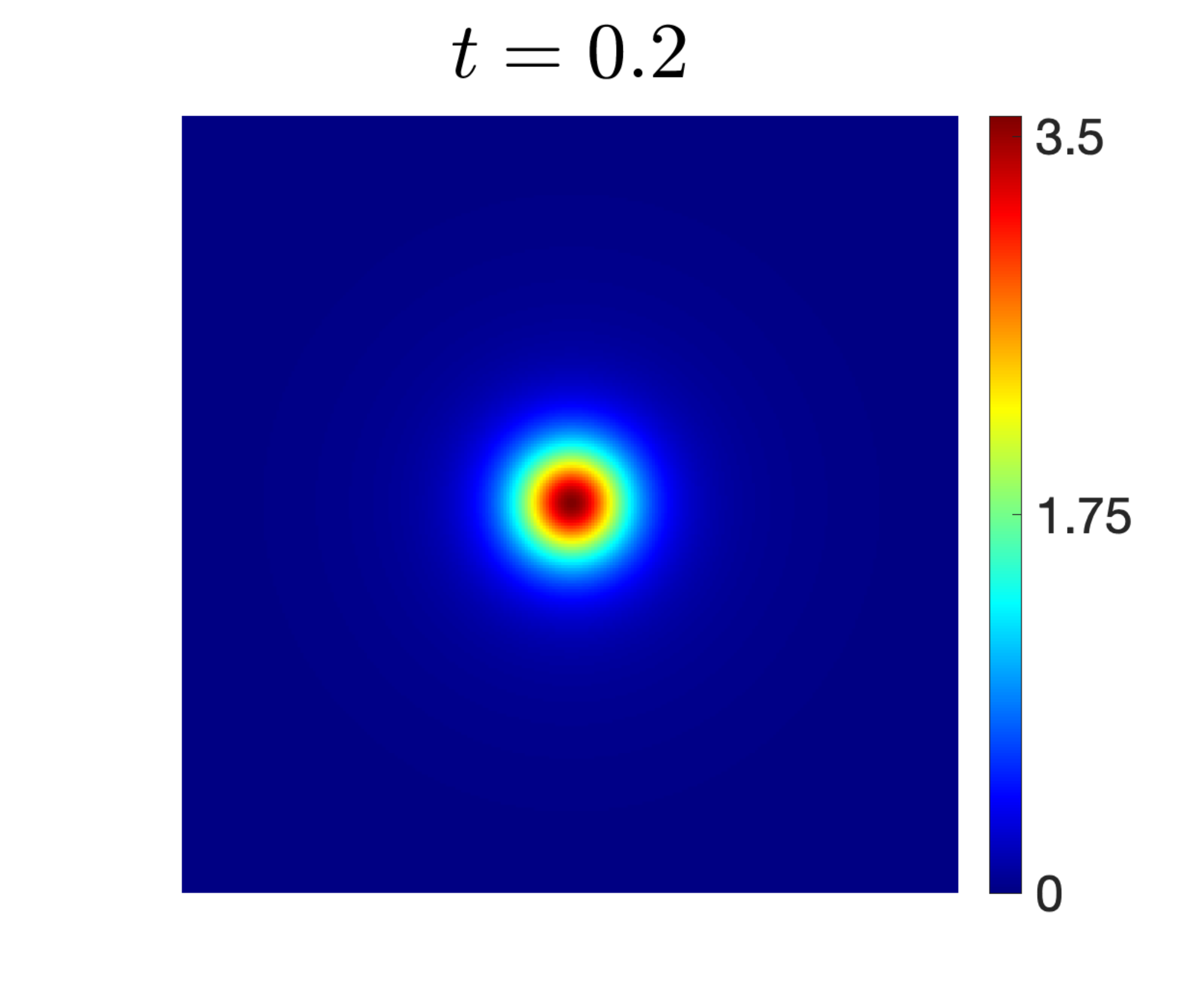}}
\end{minipage}
\caption{Plots of the density $\rho(x, y, t)$ at different times in the region $[-2, 2]^2$ for a single Gausson with the initial data in Case I for different $\alpha$: $\alpha = -0.3$ (top row);  $\alpha = -0.1$ (middle row); and  $\alpha = 1$ (bottom  row).}
\label{fig:2D_case1}
\end{figure}

\begin{figure}[ht!]
\begin{minipage}{0.33\textwidth}
\centerline{\includegraphics[width=4.2cm,height=3.5cm]{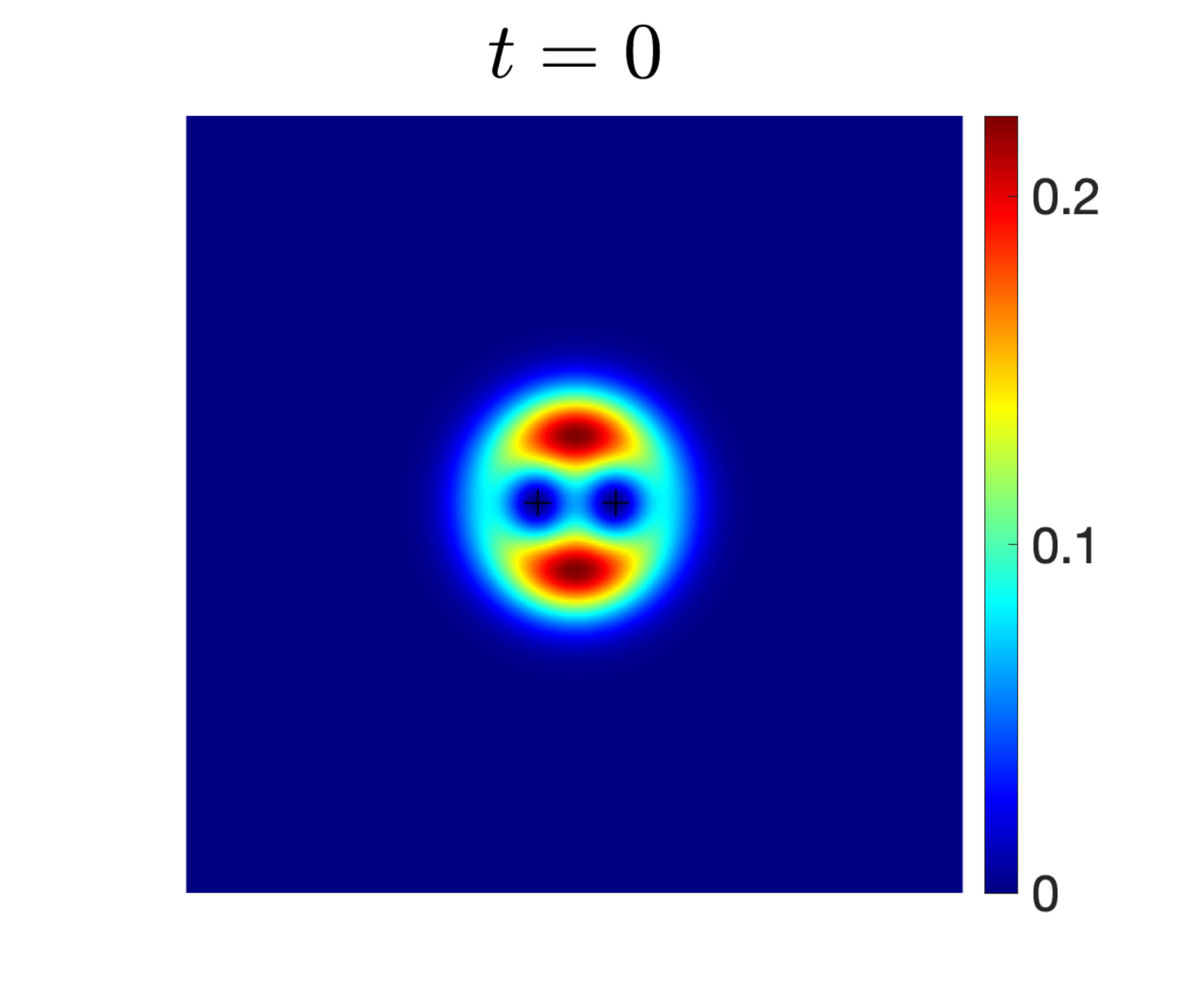}}
\end{minipage}
\begin{minipage}{0.33\textwidth}
\centerline{\includegraphics[width=4.2cm,height=3.5cm]{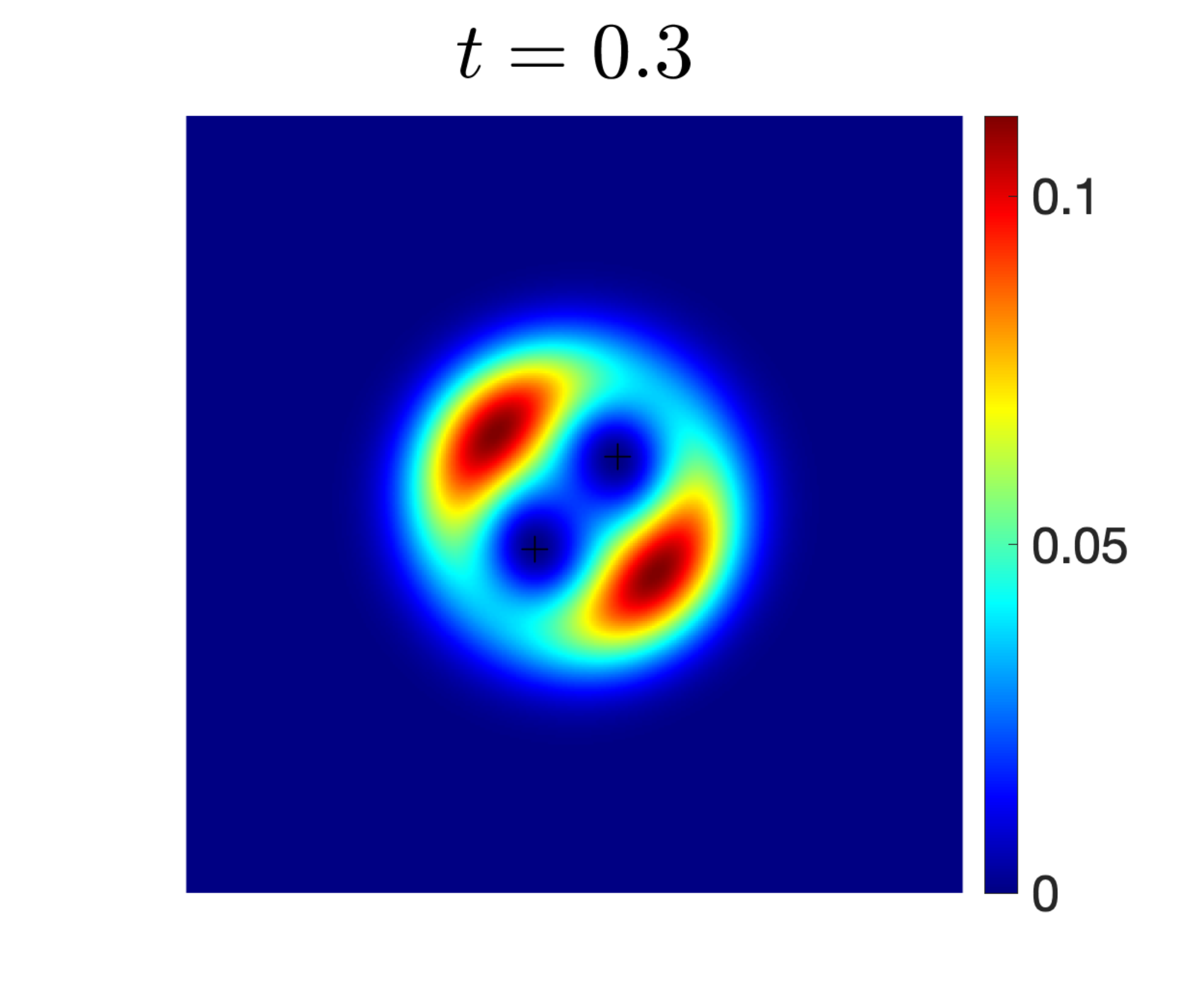}}
\end{minipage}
\begin{minipage}{0.33\textwidth}
\centerline{\includegraphics[width=4.2cm,height=3.5cm]{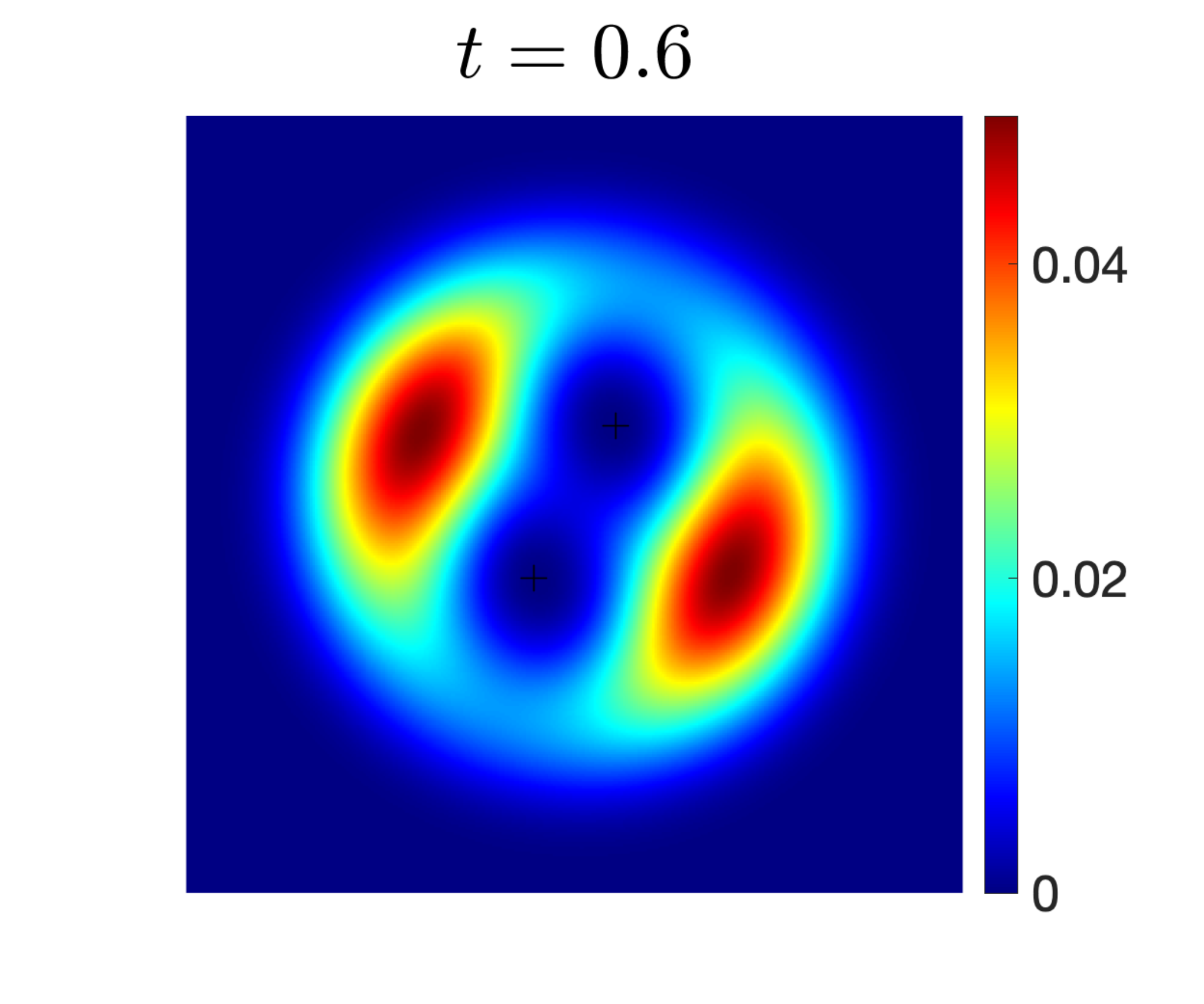}}
\end{minipage}
\begin{minipage}{0.33\textwidth}
\centerline{\includegraphics[width=4.2cm,height=3.5cm]{case2_in_den.eps}}
\end{minipage}
\begin{minipage}{0.33\textwidth}
\centerline{\includegraphics[width=4.2cm,height=3.5cm]{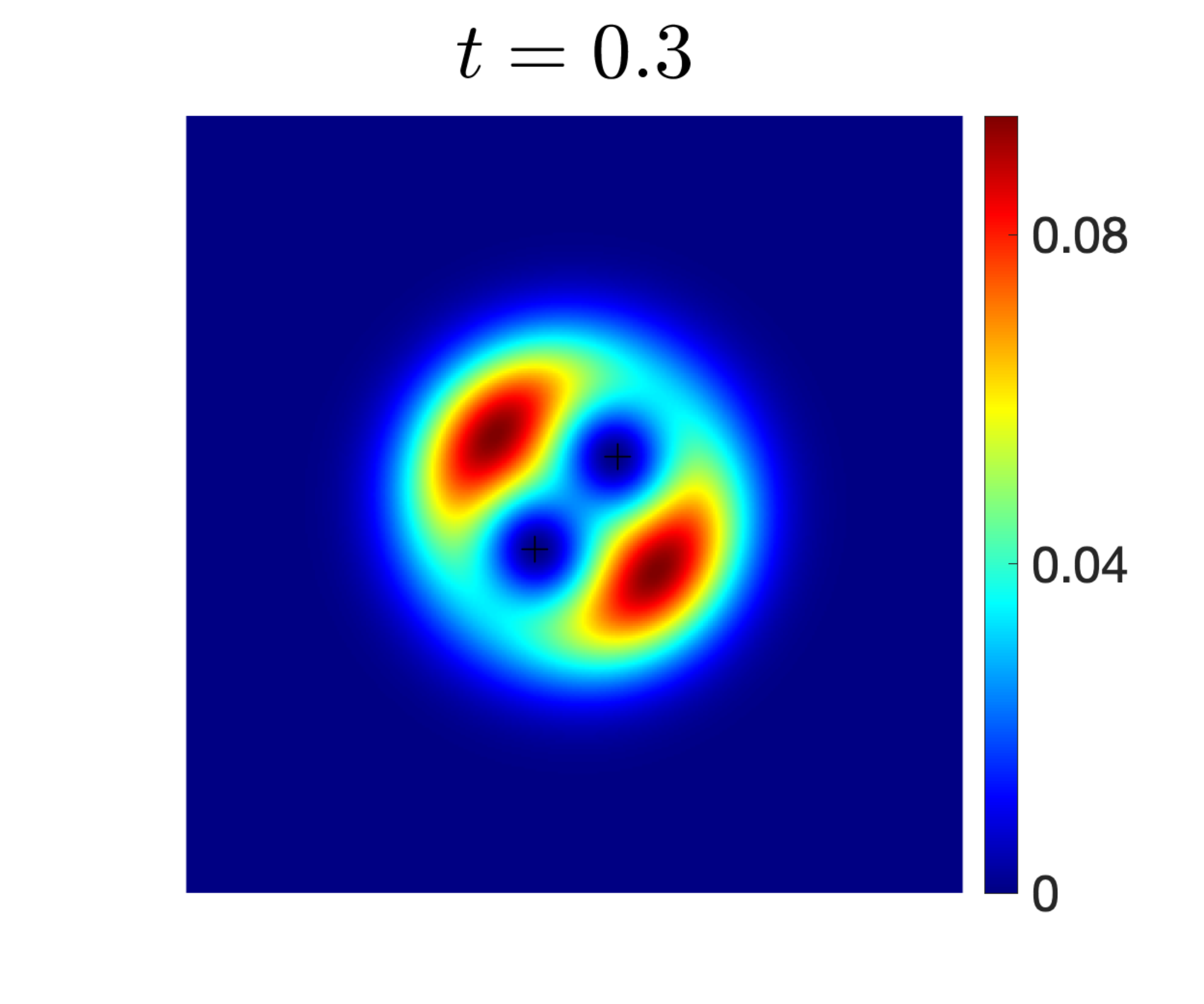}}
\end{minipage}
\begin{minipage}{0.33\textwidth}
\centerline{\includegraphics[width=4.2cm,height=3.5cm]{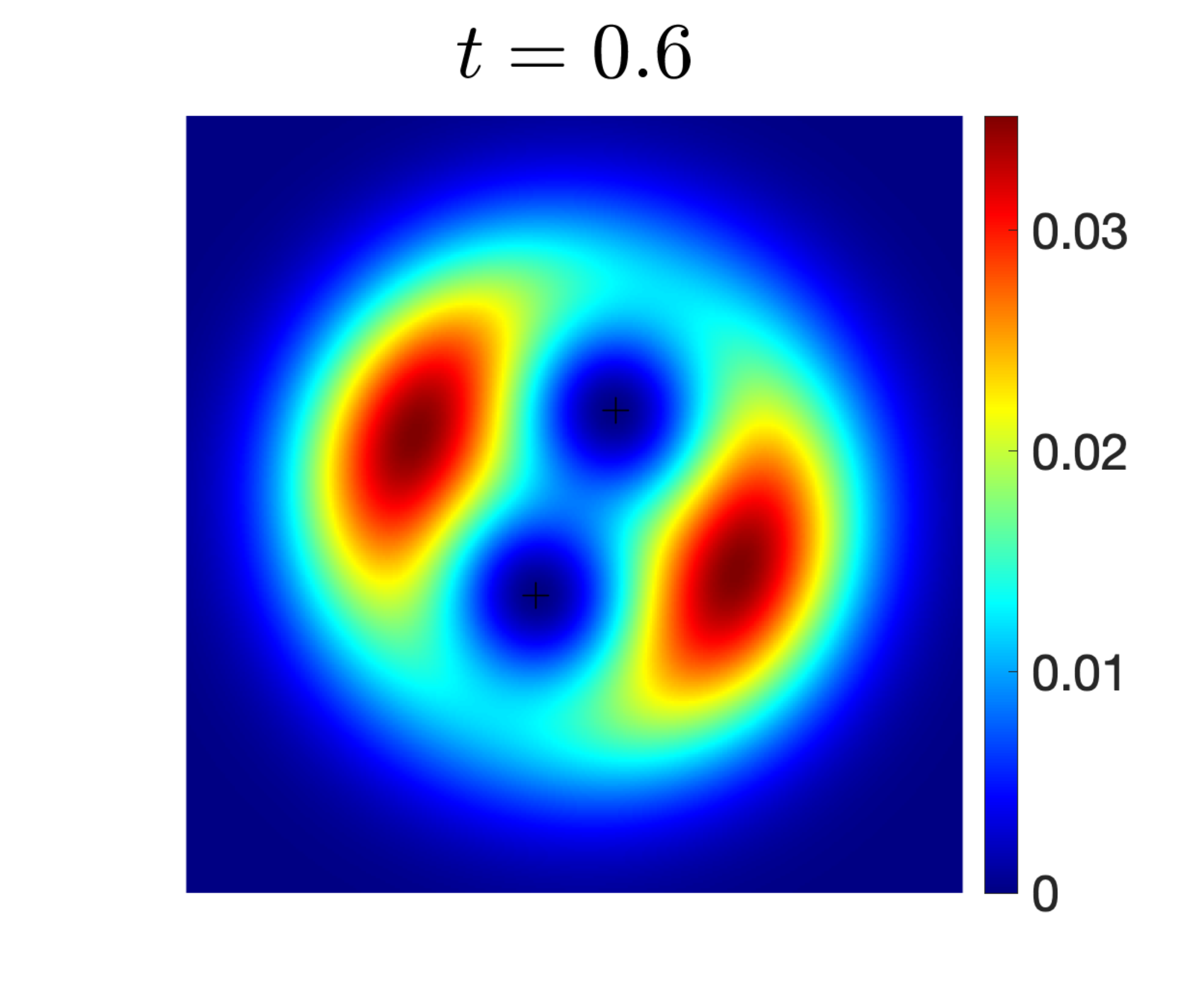}}
\end{minipage}
\begin{minipage}{0.33\textwidth}
\centerline{\includegraphics[width=4.2cm,height=3.5cm]{case2_in_den.eps}}
\end{minipage}
\begin{minipage}{0.33\textwidth}
\centerline{\includegraphics[width=4.2cm,height=3.5cm]{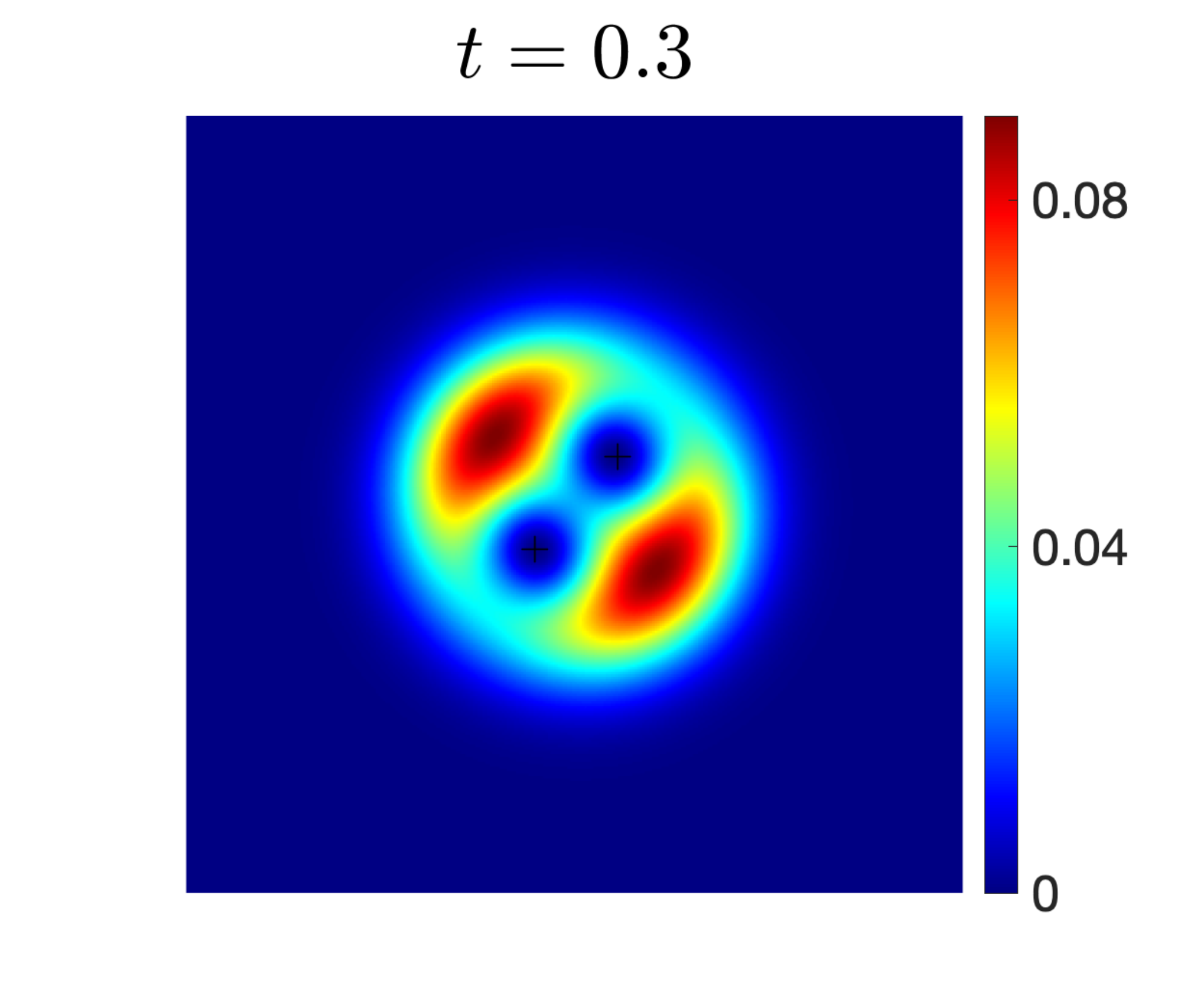}}
\end{minipage}
\begin{minipage}{0.33\textwidth}
\centerline{\includegraphics[width=4.2cm,height=3.5cm]{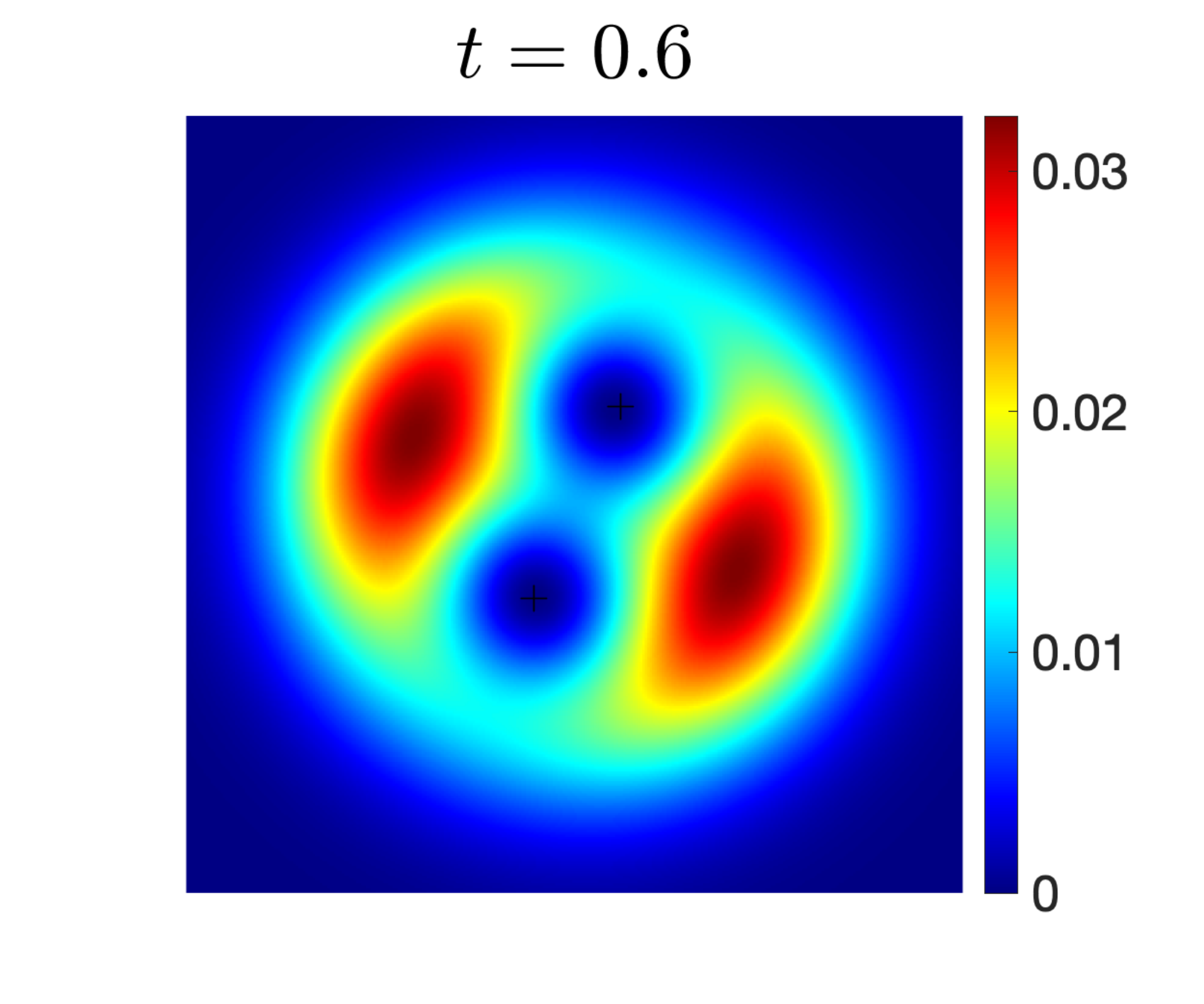}}
\end{minipage}
\caption{Plots of the density $\rho(x, y, t)$ at different times in the region $[-5, 5]^2$ for a vortex pair with the initial data in Case II for different $\alpha$: $\alpha = -0.3$ (top row);  $\alpha = -0.1$ (middle row); and  $\alpha = 1$ (bottom row).}
\label{fig:2D_case2}
\end{figure}

\begin{figure}[ht!]
\begin{minipage}{0.33\textwidth}
\centerline{\includegraphics[width=4.5cm,height=3.5cm]{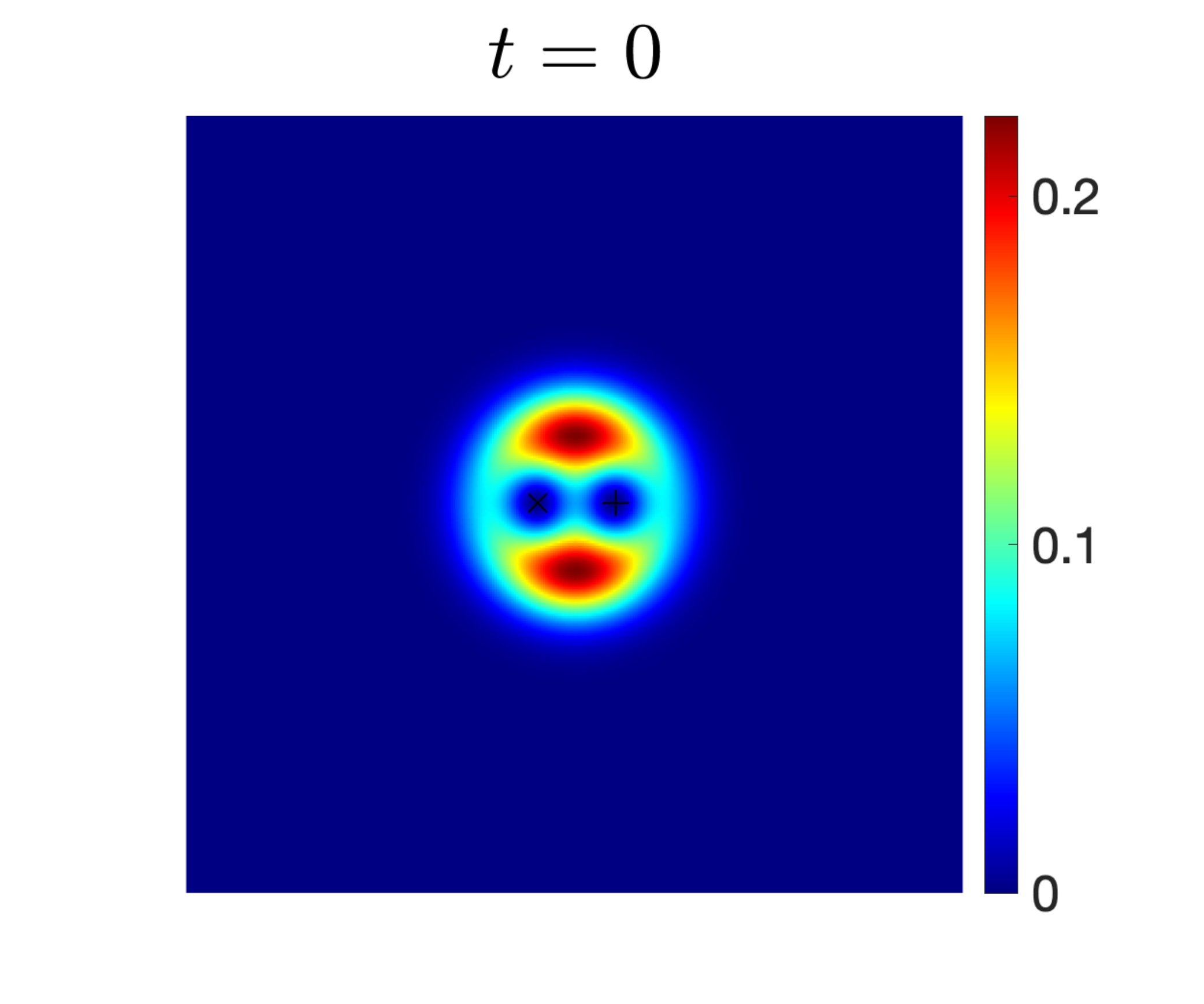}}
\end{minipage}
\begin{minipage}{0.33\textwidth}
\centerline{\includegraphics[width=4.5cm,height=3.5cm]{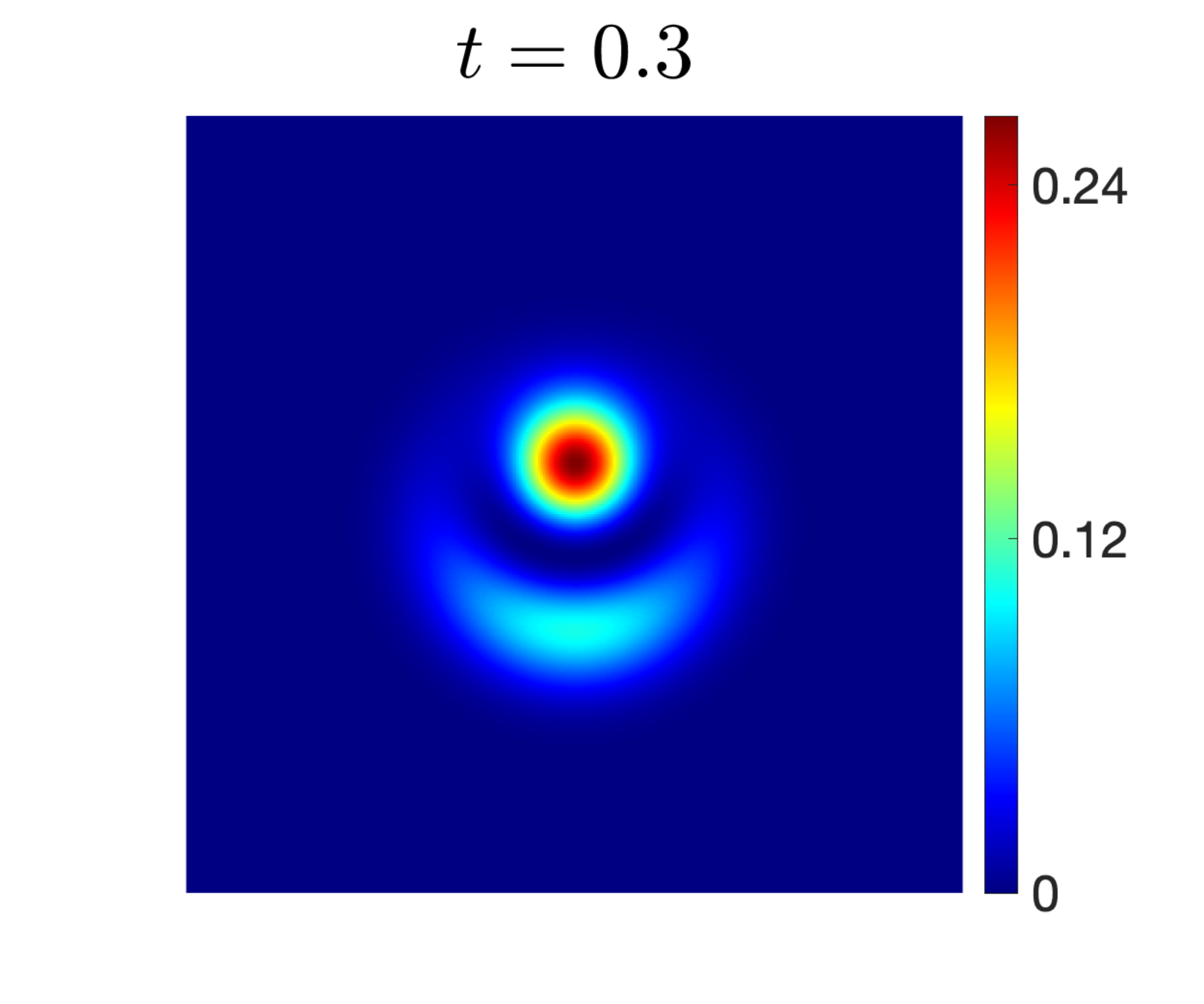}}
\end{minipage}
\begin{minipage}{0.33\textwidth}
\centerline{\includegraphics[width=4.5cm,height=3.5cm]{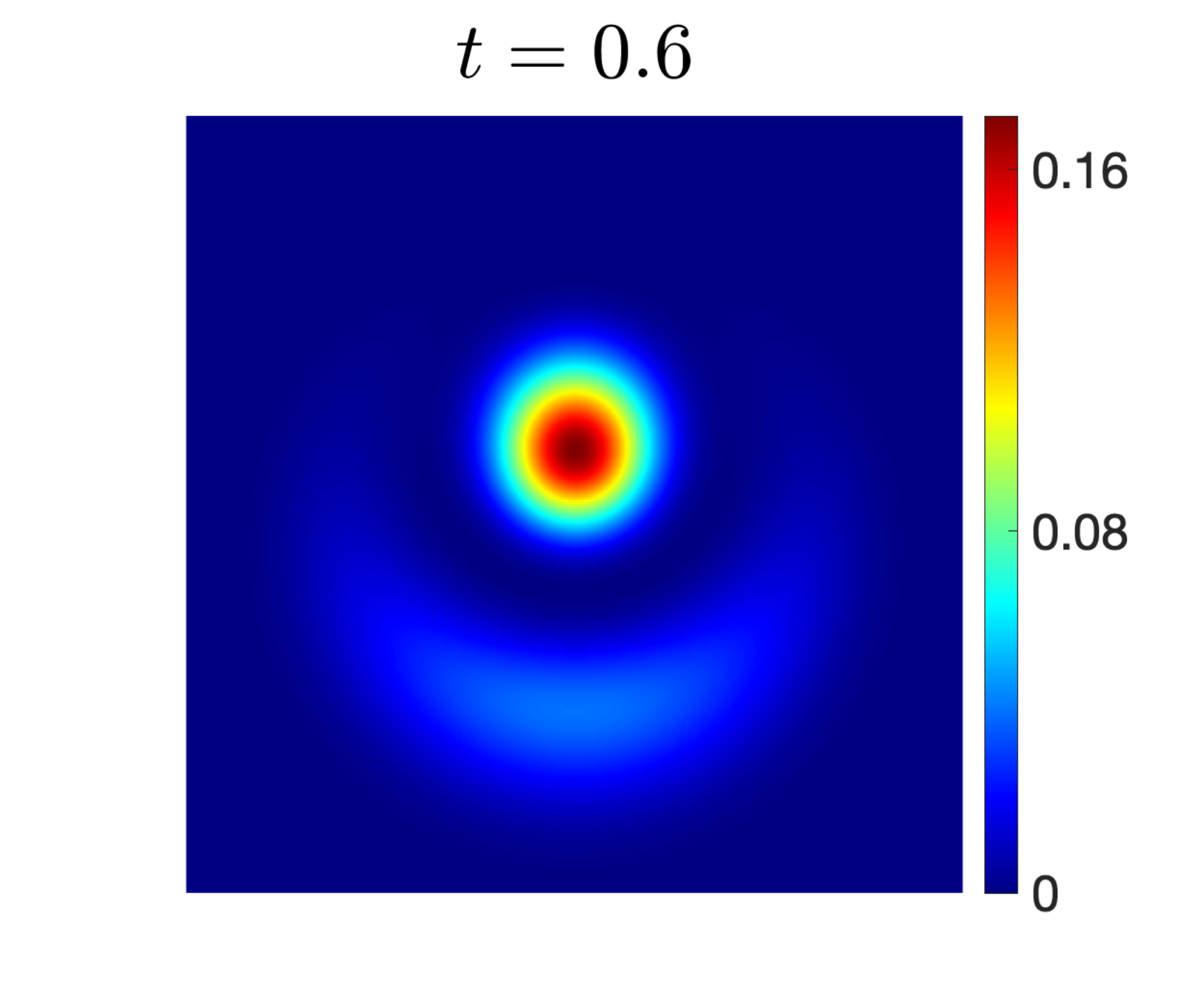}}
\end{minipage}
\begin{minipage}{0.33\textwidth}
\centerline{\includegraphics[width=4.5cm,height=3.5cm]{case3_in_den.eps}}
\end{minipage}
\begin{minipage}{0.33\textwidth}
\centerline{\includegraphics[width=4.5cm,height=3.5cm]{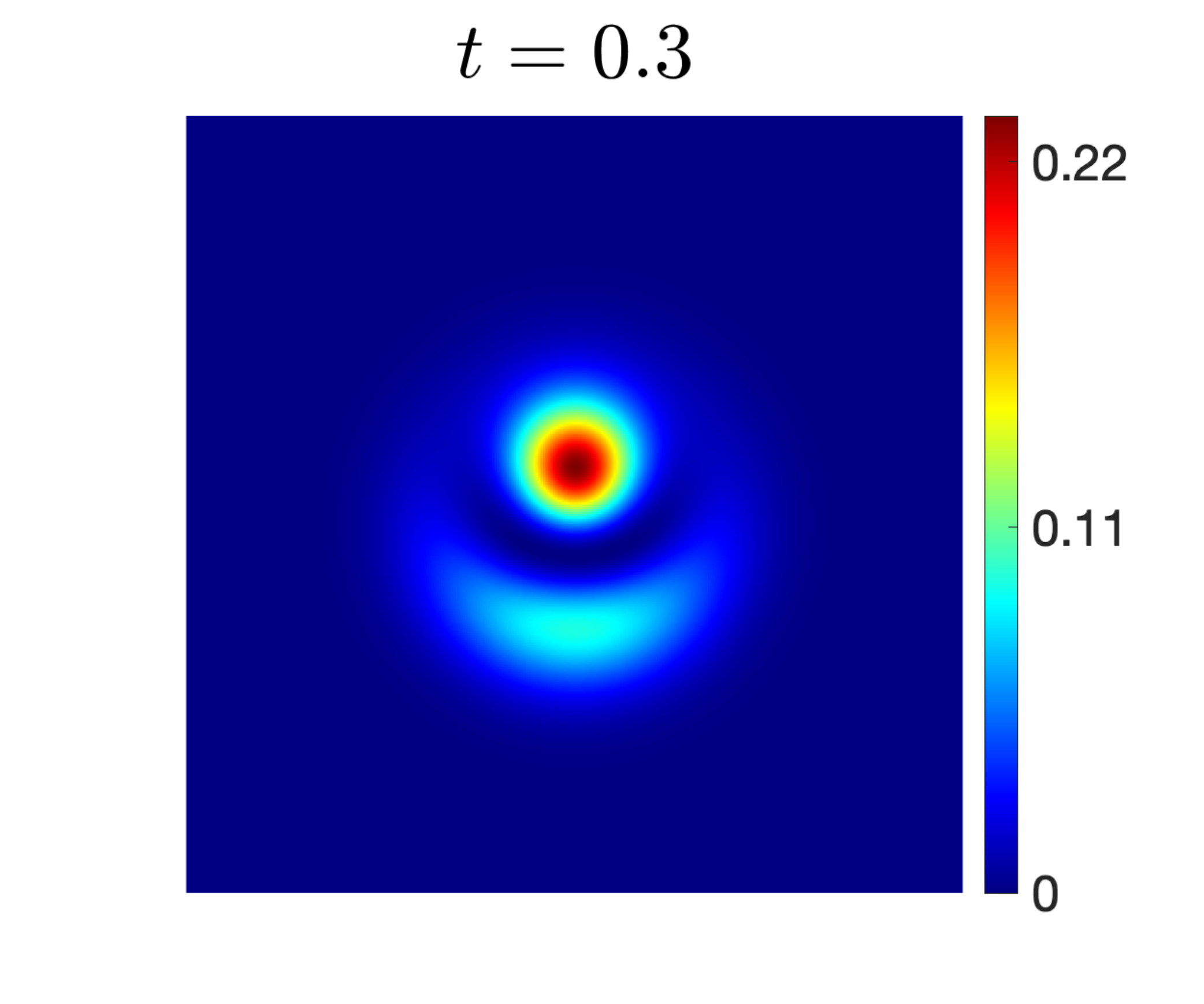}}
\end{minipage}
\begin{minipage}{0.33\textwidth}
\centerline{\includegraphics[width=4.5cm,height=3.5cm]{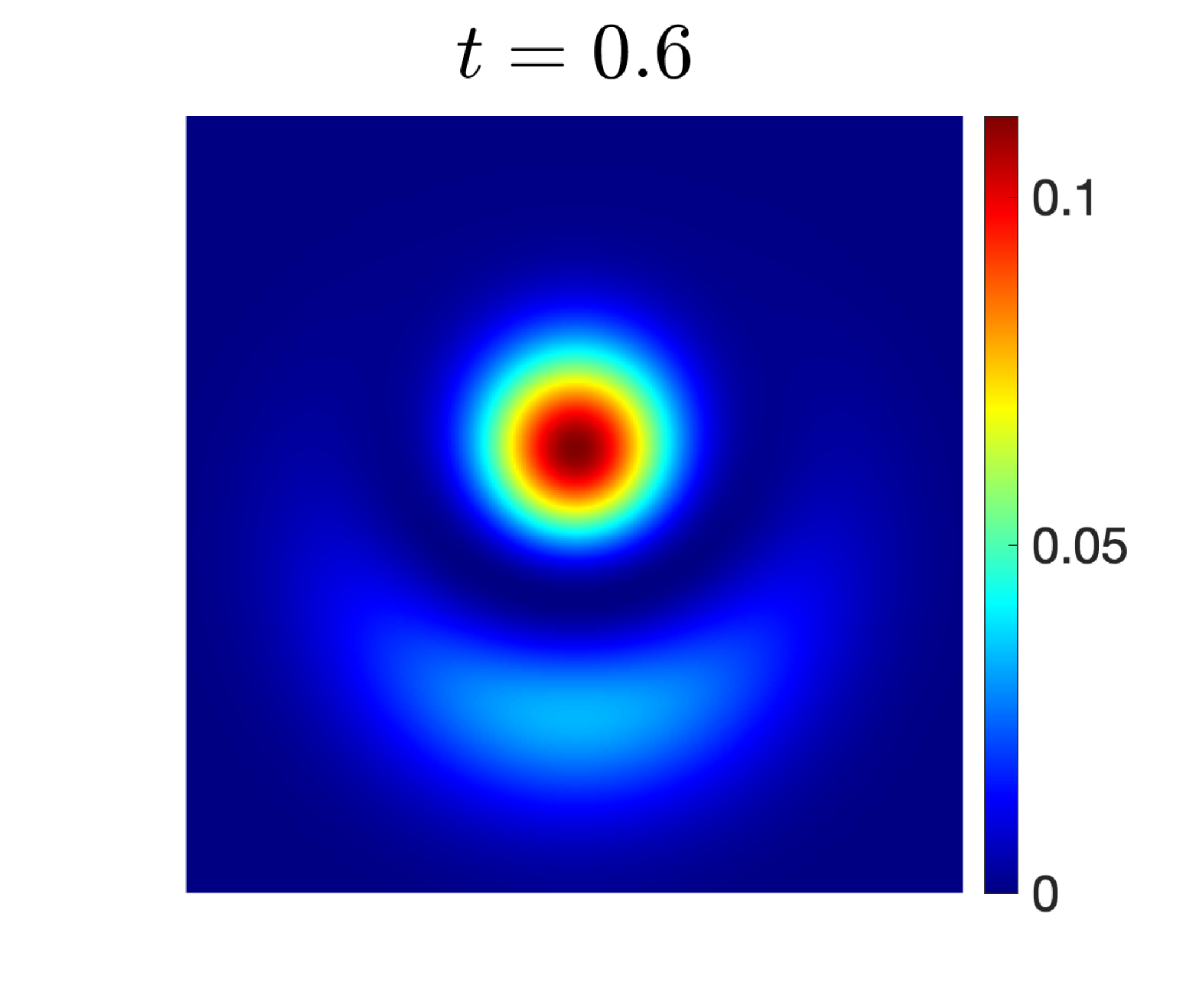}}
\end{minipage}
\begin{minipage}{0.33\textwidth}
\centerline{\includegraphics[width=4.5cm,height=3.5cm]{case3_in_den.eps}}
\end{minipage}
\begin{minipage}{0.33\textwidth}
\centerline{\includegraphics[width=4.5cm,height=3.5cm]{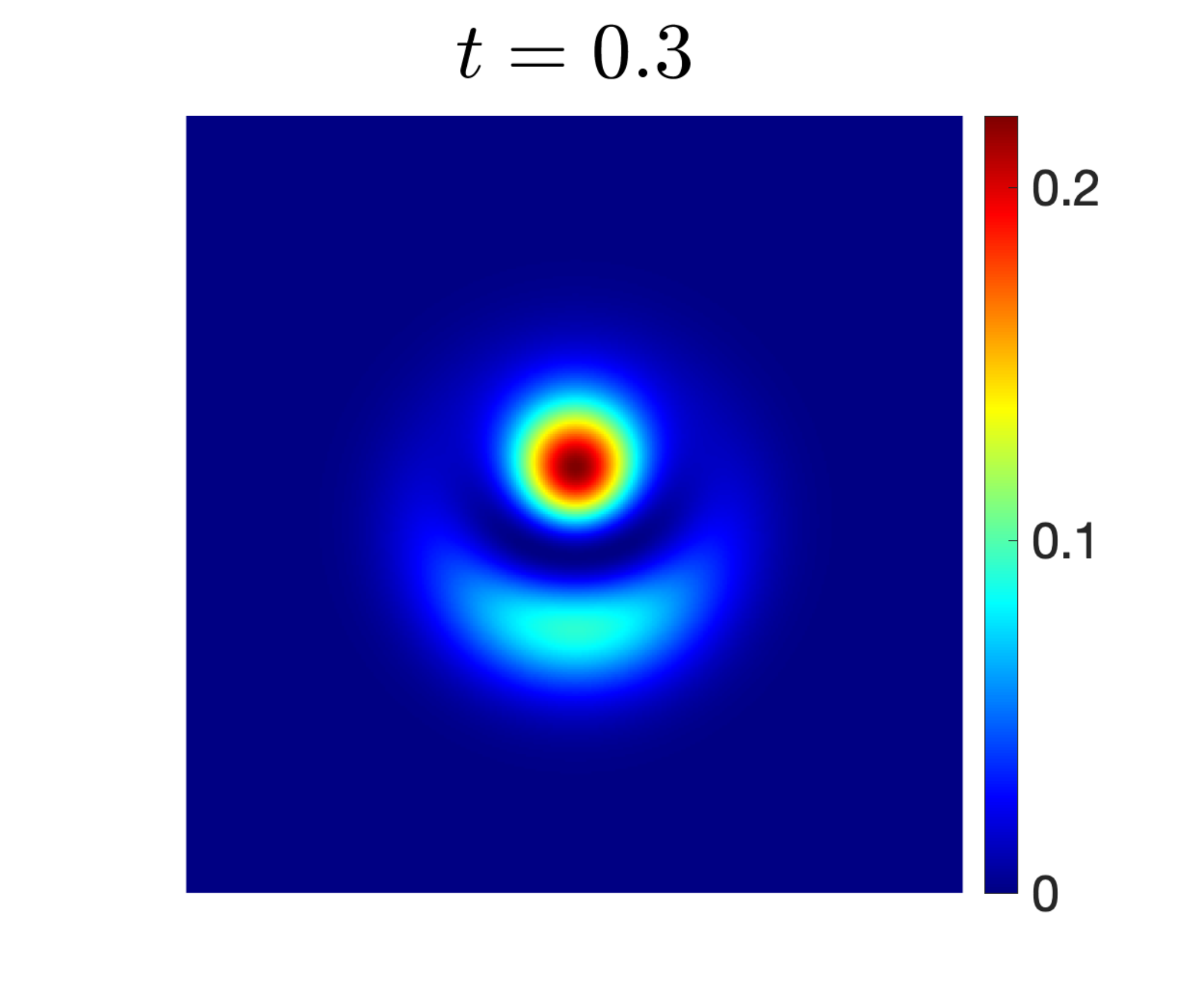}}
\end{minipage}
\begin{minipage}{0.33\textwidth}
\centerline{\includegraphics[width=4.5cm,height=3.5cm]{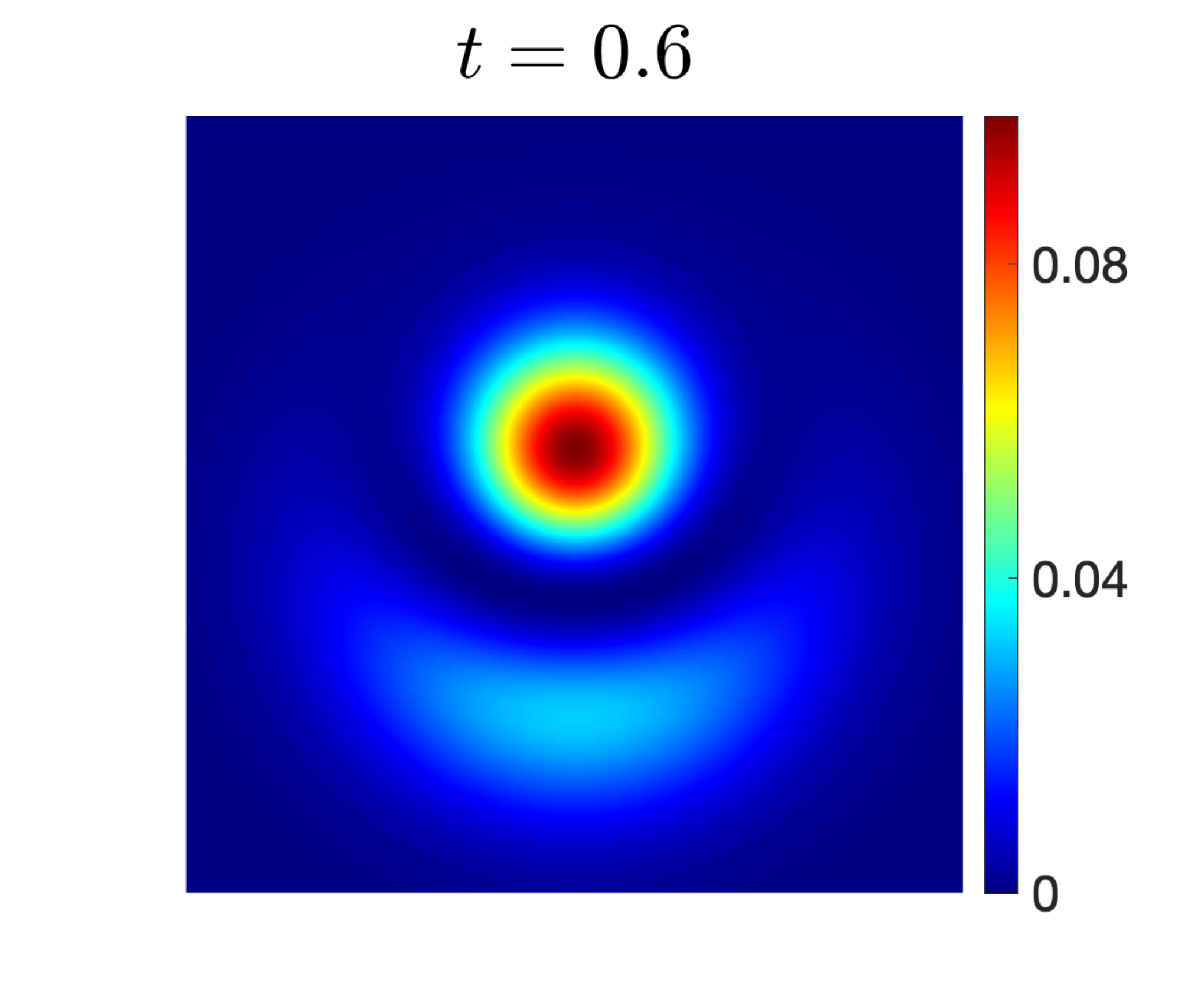}}
\end{minipage}
\caption{Plots of the density $\rho(x, y, t)$ at different times in the region $[-5, 5]^2$ for a vortex dipole with the initial data in Case III for different $\alpha$: $\alpha = -0.3$ (top row);  $\alpha = -0.1$ (middle row); and  $\alpha = 1$ (bottom row).}
\label{fig:2D_case3}
\end{figure}

We solve the problem by the TSFP method \eqref{eq:TSFP} with $\eps = 10^{-12}$, $\tau  = 0.001$. We take $\lambda = -10$, $h_x = h_y = 1/64$ and the computational domain $\Omega = [-8, 8]^2$ for Case I and  $\lambda = 1$, $h_x = h_y = 1/32$ and the computational domain $\Omega = [-16, 16]^2$ for Case II and Case III. Figures \ref{fig:2D_case1}, \ref{fig:2D_case2} and \ref{fig:2D_case3} show the density $\rho(x, y, t)$ at different times for Case I, Case II and Case III, respectively. In the figures, we use `$+$' and `$\times$' to represent the position of the vortex which has the winding number $m = 1$ and $m = -1$, respectively.

From these figures, we can draw the following conclusions:

(i) In Case I, we initially have a single Gausson and consider the focusing case, i.e., $\lambda < 0$.When the time $t$ evolves, it is still a single Gausson for all the choices of $\alpha$. However, when $\alpha = 1$, i.e., the cubic NLSE, it is concentrated and the peak value becomes larger with the time evolution. For $\alpha < 0$, it is decentralized with the time evolution and the density $\rho(x, y, t)$ depends on the exponent $\alpha$  (cf. Figure \ref{fig:2D_case1}).

(ii) In Case II, we initially have a vortex pair located at $(\pm 0.5, 0)$ with winding number $m = 1$. When the time $t$ evolves, the vortex pair rotates with each other and they never collide and annihilate (cf. Figure \ref{fig:2D_case2}). The vortex centers are independent of the exponent $\alpha$, while the size of the vortex core is much larger when $\alpha$ is smaller.

(iii) In Case III, we initially have a vortex dipole located at $(\pm 0.5, 0)$ with winding number $m = \pm 1$. When the time $t$ evolves, the two vortices start moving together,  and then they collide and disappear on the $y$-axis within a short time. In addition, the density $\rho(x, y, t)$ also depends on the exponent $\alpha$ (cf. Figure \ref{fig:2D_case3}).

\section{Conclusion}
We proposed two types of regularizations for the nonlinear Schr\"odinger equation with singular nonlinearity by introducing a small regularization parameter $0 < \eps \ll 1$ to overcome the singularity of the nonlinear term $f(\rho) = \rho^{\alpha}$ with $-1/3 < \alpha < 0$. One is by adapting the local energy regularization, which regularizes the energy density $F(\rho) = \frac{1}{\alpha+1}\rho^{\alpha+1}$ locally near $\rho = 0^+$ with a polynomial approximation. The other type is by using the global regularization, which directly regularizes the singular nonlinearity $f(\rho)$ and changes the nonlinearity for any $\rho \geq 0$. Then we presented the first-order Lie-Trotter time-splitting Fourier pseudospectral method and Lawson-type exponential integrator Fourier pseudospectral method to numerically solve the regularized models. Extensive numerical results are presented to show the convergence rates of different regularizations and to compare the performance of the numerical schemes as well as to illustrate rich dynamics of the nonlinear Schr\"odinger equation with singular nonlinearity.

\section*{Acknowledgments}
This work was partially supported by the Ministry of Education of Singapore grant A-0004265-00-00  (MOE2019-T2-1-063) (W. Bao), the European Research Council (ERC) under the European Union's Horizon 2020 research and innovation programme (grant agreement No. 850941) (Y. Feng) and the Research Foundation for Beijing University of Technology New Faculty (grant No. 006000514122521) (Y. Ma).


\end{document}